\documentclass{gen-p-l}
\usepackage[dvips]{graphicx}
\usepackage{ifthen,amsfonts,amsmath,amssymb,epic,eepic}  
\usepackage{epsfig,color}  

\theoremstyle{plain}
\newtheorem{theorem}{Theorem}[section]
\newtheorem{lemma}[theorem]{Lemma}
\newtheorem{corollary}[theorem]{Corollary}
\newtheorem{proposition}[theorem]{Proposition}

\newtheorem{question}{Question}

\theoremstyle{remark}

\newtheorem{remark}[theorem]{\bf Remark}
\theoremstyle{definition}
\newtheorem{definition}[theorem]{\bf Definition}
\newtheorem{example}[theorem]{\bf Example}

\newcommand{\dist}{{\text{dist}}}
\newcommand{\Angle}{{\text{Angle}}}
\newcommand{\Energy}{{\text{Energy}}}
\newcommand{\rp}{{\text{Re }}}
\newcommand{\ip}{{\text{Im }}}

\newcommand{\mint}{-\!\!\!\!\!\!\int}

\newcommand{\zz}{{\zeta}}
\def\ZZ{{\bold Z}}
\def\RP{{\bold RP}}
\def\RR{{\bold R}}
\def\SS{{\bold S}}
\newcommand{\V}{{\text {V}}}
\newcommand{\nn}{{\text {n}}}

\newcommand{\dv}{{\text{div }}}

\newcommand{\Vol}{{\text {Vol}}}
\newcommand{\Area}{{\text {Area}}}
\newcommand{\Ric}{{\text {Ric}}}

\def\CC{{\bold C }}
\newcommand{\eqr}[1]{(\ref{#1})}

\newcommand{\Hess}{{\text {Hess}}}
\newcommand{\e}{{\text {e}}}
\newcommand{\cB}{{\mathcal {B}}}
\newcommand{\cP}{{\mathcal {P}}}
\newcommand{\cS}{{\mathcal {S}}}
\newcommand{\cH}{{\mathcal {H}}}
\newcommand{\cL}{{\mathcal {L}}}

\newcommand{\cM}{{\mathcal {M}}}
\newcommand{\cF}{{\mathcal {F}}}
\newcommand{\cT}{{\mathcal {T}}}
\newcommand{\Conv}{\operatorname{Conv}}

\numberwithin{equation}{section}

\begin{document}

\title[An excursion into geometric analysis]
{An excursion into geometric analysis}

\author{Tobias H. Colding}%
\address{Courant Institute of Mathematical Sciences\\
251 Mercer Street\\ New York, NY 10012}
\author{William P. Minicozzi II}%
\address{Department of Mathematics\\
Johns Hopkins University\\
3400 N. Charles St.\\
Baltimore, MD 21218}
\thanks{The authors were partially supported by NSF Grants DMS
0104453 and DMS 0104187}

\email{colding@cims.nyu.edu and minicozz@math.jhu.edu}






\maketitle \tableofcontents

\section{Introduction}

This is a guided tour through some selected topics in geometric
analysis.  Most of the results here can be found in the literature
but some are new and do  not appear elsewhere. We have chosen to
illustrate many of the basic ideas as they apply to the theory of
minimal surfaces.  This is, in part, because minimal surfaces is,
if not the oldest, then certainly one of the oldest areas of
geometric analysis dating back to Euler's work in the 1740's and
in fact many of the basic ideas in geometric analysis originated
in the study of minimal surfaces.  In any case, the ideas apply to
a variety of different fields and we will mention some of these as
we go along.

Part \ref{p:1} reviews some of the classical ideas and results in
geometric analysis.  We begin with the definition  and basic
results for minimal surfaces, including the first variation
formula and maximum principle in Section \ref{s:2}.  Section
\ref{s:3} gives some applications of the Bochner formula to
comparison theorems, vanishing theorems (such as the famous
Bochner theorem), and harmonic functions.  We turn next to the
monotonicity formula and mean value inequality in Section
\ref{s:4}. These play a fundamental role in many areas of
geometric analysis, however, we have chosen to focus on the
special cases of minimal surfaces and manifolds with non--negative
Ricci curvature. Section \ref{s:5} recalls the Bernstein theorem
(for entire solutions) and Bers' theorem (for exterior solutions)
of the minimal surface equation. These results illustrate an
interesting rigidity for solutions of the minimal surface equation
which comes from the nonlinearity (and does not occur for
solutions of linear equations).
 In Sections \ref{s:6} and
\ref{s:7}, we briefly review the basic facts for mean curvature
flow and Ricci flow.  The next two sections discuss some
fundamental a priori estimates in pde.  First, Section \ref{s:8}
gives various  gradient estimates for  linear/nonlinear,
elliptic/parabolic equations, all based on the maximum principle.
The importance of this fundamental estimate has been
well-understood since the work of Bernstein in the early 1900's.
Section \ref{s:9}  recalls a much more recent tool for a priori
estimates, namely,  Simons' inequality, and illustrates its
usefulness for proving a priori estimates.  The original
inequality of Simons was for the Laplacian of the norm squared of
the second fundamental form of a minimal hypersurface, but
variations of this inequality appear in a surprising number of
fields (Einstein manifolds, harmonic maps, Yang-Mills
connections, various parabolic equations, etc.).  Finally, in
Section \ref{s:10} we derive the basic estimates for minimal
annuli with small total curvature, including a quantitative form
of Bers' theorem.  This last section also sets the stage for some
estimates in Part \ref{p:2} for multi--valued graphs.

In Part \ref{p:2}, we turn our attention to embedded multi--valued
minimal graphs (the basic example is half of the helicoid).  These
are graphs of multi--valued functions and should be thought of as
``spiral staircases.'' The analysis of these has played a major
role in recent developments in minimal surface theory.  The first
two sections, \ref{s:11} and \ref{s:12}, prove the fundamental
estimates on the separation and curvature.  Section \ref{s:13}
extends the Bers' theorem from Part \ref{p:1} to this setting.
Finally, Section \ref{s:14} proves some original results,
including a representation formula showing that an embedded
multi--valued minimal graph can be written as a sum of a helicoid,
a catenoid, and a small perturbation.

In Part \ref{p:3}, we survey some of the key  ideas in classical
regularity theory, recent developments on embedded minimal disks,
and some global results for minimal surfaces in $\RR^3$. Sections
\ref{s:16} and \ref{s:17} focus on Reifenberg type conditions,
where a set is assumed to be close to a plane at all points and at
all scales (``close'' is in the Hausdorff or Gromov--Hausdorff
sense and is defined in Section \ref{s:16}).  This condition
automatically gives H\"older regularity (and hence higher
regularity if the set is also a weak solution to a natural
equation).  Section \ref{s:18} surveys the role of monotonicity
and scaling in regularity theory, including $\epsilon$-regularity
theorems (such as Allard's theorem) and tangent cone analysis
(such as Almgren's refinement of Federer's dimension reducing).
Section \ref{s:19} briefly reviews recent results of the authors
for embedded minimal disks, developing a regularity theory in a
setting where the classical methods cannot be applied and in
particular where there is no monotonicity.  The estimates and
ideas discussed in Section \ref{s:19} have applications to the
global theory of minimal surfaces in $\RR^3$.  In Section
\ref{s:20}, we give a quick tour of some recent results in this
classical, but rapidly developing, area.

Thus far, we have mainly dealt with regularity and a priori
estimates but have ignored questions of existence.  Part \ref{p:4}
surveys some of the most useful existence results for minimal
surfaces and gives an application to Ricci flow.  Section
\ref{s:21} recalls the classical Weierstrass representation,
including a few modern applications, and the Kapouleas
desingularization method.  Section \ref{s:22} deals with area
minimizing surfaces (whether for fixed boundary, fixed homotopy
class, etc.) and questions of embeddedness.  The next section
discusses unstable (hence {\it not} minimizing) surfaces and the
corresponding questions for geodesics, concentrating on whether
the Morse index can be bounded uniformly.  Section \ref{s:24}
recalls the min--max construction for producing unstable minimal
surfaces and, in particular, doing so while controlling the
topology and guaranteeing embeddedness.  Finally, Section
\ref{s:25} discusses a recent application of min--max surfaces to
bound the extinction time for Ricci flow, answering a question of
Perelman.

Finally, in Part \ref{p:5}, we discuss some global results for
harmonic functions and a few applications of function theory. We
begin by reviewing the basic theory of harmonic functions on
Euclidean space.  This starts with the Liouville theorem and the
relationship between polynomial growth harmonic functions and
eigenfunctions on the sphere; see Section \ref{s:26}.  In
Section \ref{s:27} we sketch the proof that the spaces of harmonic
functions of polynomial growth are finite dimensional on manifolds
with non--negative Ricci curvature. Section \ref{s:28} gives a
version of this for minimal submanifolds and a geometric
application of this.  Finally, Section \ref{s:29} discusses two
estimates related to nodal  sets of eigenfunctions.

\part{Classical and almost classical results}       \label{p:1}

\section{Minimal surfaces}  \label{s:2}

Let $\Sigma\subset \RR^3$ be a \underline{smooth} orientable
surface (possibly with boundary) with unit normal $\nn_{\Sigma}$.
Given a function $\phi$ in the space $C^{\infty}_0(\Sigma)$ of
infinitely differentiable (i.e., smooth), compactly supported
functions on $\Sigma$, consider the one--parameter variation
\begin{equation}
\Sigma_{t,\phi}=\{x+t\,\phi (x)\,\nn_{\Sigma}(x) | x\in \Sigma\}\,
.
\end{equation}
The so called first variation formula of area is the equation
(integration is with respect to $d\text{area}$)
\begin{equation}  \label{e:frstvar}
\left.\frac{d}{dt} \right|_{t=0}\Area (\Sigma_{t,\phi})
=\int_{\Sigma}\phi\,H\, ,
\end{equation}
where $H$ is the mean curvature of $\Sigma$.  (When $\Sigma$ is
noncompact, then $\Sigma_{t,\phi}$ in \eqr{e:frstvar} is replaced
by $\Gamma_{t,\phi}$, where $\Gamma$ is any compact set containing
the support of $\phi$.) The surface $\Sigma$ is said to be a {\it
minimal} surface (or just minimal) if
\begin{equation}
\left.\frac{d}{dt} \right|_{t=0}\Area (\Sigma_{t,\phi})=0
\,\,\,\,\,\,\,\,\,\,\,\text{ for all } \phi\in
C^{\infty}_0(\Sigma)
\end{equation}
 or, equivalently by \eqr{e:frstvar}, if the
mean curvature $H$ is identically zero.  Thus $\Sigma$ is minimal
if and only if it is a critical point for the area functional.
(Since a critical point is not necessarily a minimum the term
``minimal'' is misleading, but it is time honored.  The equation
for a critical point is also sometimes called the Euler--Lagrange
equation.) Moreover, a computation shows that if $\Sigma$ is
minimal, then
\begin{equation}
\left. \frac{d^2}{dt^2} \right|_{t=0}\Area (\Sigma_{t,\phi})
=-\int_{\Sigma}\phi\,L_{\Sigma}\phi\, ,
\,\,\,\,\,\,\,\,\,\,\,\text{ where }
L_{\Sigma}\phi=\Delta_{\Sigma}\phi+|A|^2\phi
\end{equation}
is the second variational (or Jacobi) operator. Here
$\Delta_{\Sigma}$ is the Laplacian on $\Sigma$ and $A$ is the
second fundamental form.  So $|A|^2=\kappa_1^2+\kappa_2^2$,
 where $\kappa_1,\,\kappa_2$ are the principal curvatures of $\Sigma$
and $H=\kappa_1+\kappa_2$. A minimal surface $\Sigma$ is said to
be stable if
\begin{equation}
\left. \frac{d^2}{dt^2} \right|_{t=0}\Area (\Sigma_{t,\phi})\geq 0
\,\,\,\,\,\,\,\,\,\,\,\text{ for all } \phi\in
C^{\infty}_0(\Sigma)\, .
\end{equation}
A graph (i.e., the set $\{(x_1,x_2,u(x_1,x_2))\,|\,
(x_1,x_2)\in \Omega\}$) of a real
valued function $u$ on a domain $\Omega$ in $\RR^2$ is minimal iff the
function satisfies the minimal surface
equation
\begin{equation}  \label{e:mineq}
\text{div} \left( \frac{d u}{\sqrt{1+|d u|^2}}\right)=0\, ,
\end{equation}
where $du$ is the $\RR^2$ gradient of the function $u$ and
$\text{div}$ is the divergence in $\RR^2$.
 One can
show that a minimal graph is stable and, more generally, so is a
multi--valued minimal graph (see below for the precise
definition).

We will next derive the weak form of the minimal surface equation,
i.e., the so--called first variation formula, which is the basic
tool for working with ``weak solutions'' (typically, stationary
varifolds). Let $X$ be a vector field on $\RR^3$.  We can write
the divergence $\dv_{\Sigma} \, X$ of $X$ on $\Sigma$ as
\begin{equation}  \label{e:o1.3.1}
        \dv_{\Sigma} \, X=\dv_{\Sigma}\, X^T
        + X \cdot H \, ,
\end{equation}
where $X^T$ and $X^N$ are the tangential and normal projections of
$X$. From this and Stokes' theorem, we see that $\Sigma$ is
minimal if and only if for all vector fields $X$ with compact
support and vanishing on the boundary of $\Sigma$,
\begin{equation}   \label{e:o1.3.2}
        \int_{\Sigma}\dv_{\Sigma}\, X=0\, .
\end{equation}
This equation is known as the first variation formula.  It has the
benefit that \eqr{e:o1.3.2} makes sense as long as we can define
the divergence on $\Sigma$.   As a consequence of \eqr{e:o1.3.2},
 we will show the following proposition:

\begin{proposition}
\label{p:haco} $\Sigma^k\subset \RR^n$ is minimal if and only if
the restrictions of the coordinate functions of $\RR^n$ to
$\Sigma$ are harmonic functions.
\end{proposition}

\begin{proof}
Let $\eta$ be a smooth function on $\Sigma$ with compact support
and $\eta |\partial \Sigma=0$, then
\begin{equation}
        \int_{\Sigma} \langle \nabla_{\Sigma} \eta , \,
        \nabla_{\Sigma} x_i \rangle = \int_{\Sigma} \langle  \nabla_{\Sigma}\eta \, ,
        e_i \rangle =  \int_{\Sigma} \dv_{\Sigma} (\eta\, e_i)\, .
\end{equation}
From this, the claim follows easily.
\end{proof}

Recall that if $\Xi\subset \RR^n$ is a compact subset, then the
smallest convex set containing $\Xi$ (the convex hull, $\Conv
(\Xi)$) is the intersection of all half--spaces containing $\Xi$.
The maximum principle forces a  minimal submanifold to lie in the
convex hull of its boundary  (this is the  ``convex hull
property''):

\begin{proposition}
If $\Sigma^k\subset \RR^n$ is a compact minimal submanifold,
 then
$\Sigma\subset \Conv (\partial \Sigma)$.
\end{proposition}

\begin{proof}
A half--space $H \subset \RR^n$ can be written as
\begin{equation}
H = \{ x \in \RR^n \, | \, \langle x , e \rangle \leq a \} \, ,
\end{equation}
 for a vector $e \in \SS^{n-1}$ and constant $a \in \RR$.
  By Proposition \ref{p:haco},   the function
$u(x)=\langle e, x\rangle$ is harmonic on $\Sigma$ and hence
attains its maximum on $\partial \Sigma$ by the maximum principle.
\end{proof}

The argument in the proof of the convex hull property can be
rephrased as saying that as we translate a hyperplane towards a
minimal surface, the first point of contact must be on the
boundary.  When $\Sigma$ is a hypersurface, this is a special case
of the strong maximum principle for minimal surfaces:

\begin{lemma}             \label{l:smp}
Let $\Omega \subset \RR^{n-1}$ be an open connected neighborhood
of the origin. If $u_1$, $u_2:\Omega\to \RR$ are solutions of the
minimal surface equation with $u_1\leq u_2$ and $u_1(0)=u_2(0)$,
then $u_1\equiv u_2$.
\end{lemma}

See \cite{CM1} for a proof of Lemma \ref{l:smp} and further
discussion.
\section{The Bochner formula}
\label{s:3}

On a Riemannian manifold $M$ a very useful formula of S. Bochner
asserts that for any function $u$ on $M$
\begin{equation}    \label{e:bochner}
\frac{1}{2}\Delta |\nabla u|^2 =|\Hess_u|^2  +\langle \nabla
\Delta u,\nabla u\rangle +\Ric_M(\nabla u,\nabla u)\, .
\end{equation}
Two special cases of this formula are particularly useful.  When $u$ is a
distance function, that is, when $|\nabla u|=1$, then the above formula
reduces to
\begin{equation}    \label{e:riccati}
0 = |U|^2 +\text{Tr}(U')+\Ric_M\, ,
\end{equation}
where $U$ is the Hessian of $u$ and the Ricci curvature and the
derivative $U'$ is taken in the direction of the unit vector
$\nabla u$. This is the so--called Ricatti equation.  The other
useful special case of the Bochner formula is when $u$ is a
harmonic function.  In this case, the Bochner formula reduces to
\begin{equation}
\frac{1}{2}\Delta |\nabla u|^2 =|\Hess_u|^2 +\Ric_M(\nabla
u,\nabla u)\, .
\end{equation}
So when $M$ has non--negative Ricci curvature, this formula
implies that the energy density of $u$ is subharmonic.

The  Laplacian and Hessian comparison theorems relate the distance
function on $M$ to a space of constant curvature.  These
comparisons are essentially integrated forms of the Ricatti
equation \eqr{e:riccati}.    For simplicity, we will not state the
most general forms of these theorems, but   rather only state the
comparisons with $\RR^n$.

The Laplacian comparison theorem  compares $\Delta r$, where $r$
is the distance to a point, on $M$ with $\Delta |x| = (n-1)/|x|$
on Euclidean space:

\begin{theorem}     \label{t:laplacian}
If $M$ has non--negative Ricci curvature and $r$ is the distance
function to a fixed point $p$, then
\begin{equation} \label{e:laplacian}
    \Delta r  \leq \frac{n-1}{   r}
    \, .
\end{equation}
Moreover, \eqr{e:laplacian} holds weakly even where $r$ is not
smooth.
\end{theorem}

\begin{proof}
We will prove \eqr{e:laplacian} assuming that $r$ is smooth so
that $|\nabla r| = 1$ (see \cite{Ca} for the extension to the
general case). Let $\gamma$ be a geodesic from $p$ parametrized by
arclength and set $L(t) = \Delta r \circ \gamma (t)$. Note that
$L'(t) = \langle \nabla \Delta r , \nabla r \rangle$ by the chain
rule so that \eqr{e:riccati} gives
\begin{equation}    \label{e:rica}
 L' = \text{Tr}(U') \leq - |U|^2 \leq - L^2   / (n-1) \, .
\end{equation}
Here the second inequality used the Cauchy-Schwarz inequality
\begin{equation}    \label{e:csmat}
    \left( \sum_{i=1}^{n-1} \lambda_i
        \right)^2 \leq (n-1) \, \sum_{i=1}^{n-1} \lambda_i^2
\end{equation}
for the eigenvalues $\lambda_i$ of the matrix $U$ (there are at
most $(n-1)$ non--zero eigenvalues since $U(\nabla r , \nabla r) =
0$). We can rewrite   \eqr{e:rica} as
\begin{equation}    \label{e:rica2}
 (1/L)'  \geq 1 / (n-1) \, .
\end{equation}
(Notice that we get equality in \eqr{e:rica2} for $L(t) =
(n-1)/t$.)  Since any manifold is ``almost Euclidean'' for $r$
small, it is easy to see that
\begin{equation}  \label{e:rsmall}
    \lim_{r
\to 0} r\, \Delta r = (n-1) \, .
\end{equation}
 Integrating the
differential equality \eqr{e:rica2} and substituting the
``boundary condition'' \eqr{e:rsmall} gives
\begin{equation}        \label{e:doner}
    \Delta r \circ \gamma(t) = L(t) \leq \frac{n-1}{   t} \, .
\end{equation}
Since $\gamma$ was arbitrary, the theorem follows.
\end{proof}

We note two immediate consequences of Theorem \ref{t:laplacian}:
\begin{itemize}
\item
Since $\nabla r$ is the unit normal to the geodesic spheres,
    the mean
curvature of these spheres is at most $(n-1)/r$.
\item
The square of the distance function satisfies $\Delta r^2 = 2
|\nabla r|^2 + 2 r \, \Delta r \leq 2n$.
\end{itemize}

\vskip2mm The Hessian comparison theorem is somewhat more
restrictive since it requires bounds on the sectional curvatures
of $M$; of course, the conclusion is correspondingly stronger.
The following theorem is a useful special case of the Hessian
comparison theorem:

\begin{theorem}     \label{t:hessian}
If $M$ is simply connected with non--positive sectional curvature
and $r$ is the distance function to a fixed point $p$, then
\begin{equation} \label{e:hessian}
    \Hess_r (X,X) \geq \frac{\left| X - \langle X , \nabla r \rangle \nabla r
    \right|^2}{
    r}
    \, ,
\end{equation}
for any vector $X$.
\end{theorem}

\vskip2mm An important application of   \eqr{e:bochner} (and
similar formulas) is to prove vanishing theorems relating a
pointwise curvature condition to  global properties of $M$. The
prototype is the Bochner theorem (see also \cite{C1}, \cite{C2}
for an extension of this famous theorem of Bochner that had been
conjectured by M. Gromov):

\begin{theorem}     \label{t:bochner}
\cite{Bc} If $M^n$ is closed with $\Ric_M \geq 0$, then each
harmonic $1$--form is parallel.  In particular, the space of
harmonic $1$--forms is at most $n$-dimensional.
\end{theorem}

\begin{proof}
(Sketch) A harmonic $1$--form $\alpha$ can be written locally as
$du$ where $u$ is a (locally defined) harmonic function.  In
particular, \eqr{e:bochner} implies that
\begin{equation}
    \Delta |\alpha|^2 = \Delta |\nabla u|^2 \geq 2 \, |\Hess_u|^2
    =2  | \nabla \alpha|^2 \, .
\end{equation}
Since $M$ is closed (in particular, $\partial M = \emptyset$),
Stokes' theorem gives
\begin{equation}
   0 = \int \Delta |\alpha|^2   \geq 2 \int | \nabla \alpha|^2 \,
   .
\end{equation}
Therefore, $| \nabla \alpha|^2$ vanishes identically.
\end{proof}

Therefore, by the Hodge theorem,  the first betti number of a
closed manifold $M$ with non--negative Ricci curvature is at most
$n$ with equality  only if the universal cover of $M$ is $\RR^n$.
There have been many geometric applications of this method, where
analytic methods (like the Hodge theorem) use topology to produce
solutions of a pde and then
  a curvature condition (like the Bochner formula)
places restrictions on these solutions.

\vskip2mm Finally, we note that  \eqr{e:bochner} can be used
  to prove an eigenvalue comparison theorem when $M$ has
positive Ricci curvature.  Namely, A. Lichnerowicz showed that if
$\Ric_M \geq \Ric_{\SS^n}$, then the first (non--zero) eigenvalue
$\lambda_1 (M) \geq \lambda_1 (\SS^n)$:

\begin{theorem}     \label{t:lichn}
\cite{Lc} If $M^n$ is closed with $\Ric_M \geq (n-1)$, then
$\lambda_1 (M) \geq n$.
\end{theorem}

\begin{proof}
Let $u$ be a (non--constant) eigenfunction on $M$ with $\Delta u =
- \lambda u$. We will show that $\lambda \geq n$.  After
normalizing so
 $\int u^2 = 1$,
Stokes' theorem gives
\begin{equation}    \label{e:l2}
    \int |\nabla u|^2 = - \int u \, \Delta u = \lambda \int u^2 =
    \lambda \, .
\end{equation}
 Substituting the equation
for $u$ into the Bochner formula \eqr{e:bochner} gives
\begin{equation}    \label{e:bochnerl}
\frac{1}{2}\Delta |\nabla u|^2 \geq |\Hess_u|^2  + (n-1- \lambda)
\, | \nabla u |^2 \geq \frac{\lambda^2}{n} \,   u^2 + (n-1-
\lambda) \, | \nabla u |^2 \, ,
\end{equation}
where the last inequality used   the Cauchy-Schwarz inequality as
in \eqr{e:csmat}.  Integrating \eqr{e:bochnerl}  over $M$ gives
\begin{equation}    \label{e:bochnerll}
0 = \frac{1}{2} \int \Delta |\nabla u|^2 \geq \frac{\lambda^2}{n}
\int   u^2 + (n-1- \lambda) \, \int | \nabla u |^2 =\lambda \,
\frac{n-1}{n} \, \left(n  - \lambda  \right)
 \, .
\end{equation}
\end{proof}

\begin{remark}
These comparison theorems are sharp in the sense that equality is
achieved on the model spaces.  The converse of this is also true,
i.e., equality is achieved only for the model spaces, and is known
as ``rigidity.''  For example, if $M^n$ is closed with $\Ric_M
\geq (n-1)$ and $\lambda_1 (M) = n$, then M. Obata, \cite{Ob},
proved that $M = \SS^n$.  It is then natural to ask how stable is
this rigidity -- i.e., what happens if equality is almost
achieved? These questions, known as ``almost rigidity,'' were
answered by Colding and J. Cheeger--Colding; see \cite{C2},
\cite{ChC4} and references therein.  Moreover, almost rigidity
theorems have played a key role in regularity theory; see
\cite{ChC1}, \cite{ChC2}, \cite{ChC3}.
\end{remark}

\section{Monotonicity and the mean value inequality}
\label{s:4}

Monotonicity formulas and mean value inequalities play a
fundamental role in many areas of geometric analysis.  In this
section, we focus on the specific cases of minimal surfaces and
manifolds with non--negative Ricci curvature.

Before we state and prove the monotonicity formula of volume for
minimal submanifolds, we will need to recall the coarea formula.
This formula asserts (see, for instance, \cite{Fe} for a proof)
that if $\Sigma$ is a manifold and $h:\Sigma\to \RR$ is a proper
(i.e., $h^{-1}((-\infty,t])$ is compact for all $t\in \RR$)
Lipschitz function on $\Sigma$, then for all locally integrable
functions $f$ on $\Sigma$ and $t\in \RR$
\begin{equation}                \label{e:thecoareaf}
        \int_{\{h\leq t\}} f\,|\nabla h|
                =\int_{-\infty}^t \int_{h=\tau} f\,d\tau\, .
\end{equation}

\begin{proposition}
\label{p:pmonot} Suppose that $\Sigma^k\subset \RR^n$ is a minimal
submanifold and $x_0\in \RR^n$; then for all $0<s<t$
\begin{equation}        \label{e:monot0}
        t^{-k}\, \Vol (B_t (x_0) \cap \Sigma) -s^{-k}\, \Vol (B_s\cap \Sigma)
        =\int_{(B_t(x_0) \setminus B_s(x_0) )\cap \Sigma} \frac{|(x-x_0)^N|^2}{|x-x_0|^{k+2}}\, .
\end{equation}
\end{proposition}

\begin{proof}
Within this proof, we set $B_t = B_t (x_0)$.  Since $\Sigma$ is
minimal,
\begin{equation}
        \Delta_{\Sigma} |x-x_0|^2=2 \, \dv_{\Sigma} (x-x_0)=2k \, .
\end{equation}
By Stokes' theorem integrating this gives
\begin{equation}
2\, k\, \Vol (B_s\cap \Sigma)=\int_{B_s\cap \Sigma}
\Delta_{\Sigma} |x-x_0|^2 =2\int_{\partial B_s\cap \Sigma}
|(x-x_0)^T|\, .
\end{equation}
Using this and the coarea formula (i.e., \eqr{e:thecoareaf}), an
easy calculation gives
\begin{align}
        \frac{d}{ds}\left( s^{-k}\, \Vol (B_s\cap \Sigma)\right)
        &=-k\, s^{-k-1}\, \Vol (B_s\cap \Sigma)
        +s^{-k}\,\int_{\partial B_s\cap \Sigma}
                \frac{|x-x_0|}{|(x-x_0)^T|}\notag\\
        &=s^{-k-1}\int_{\partial B_s\cap \Sigma}
        \left( \frac{|x-x_0|^2}{|(x-x_0)^T|}-|(x-x_0)^T| \right)\\
        &=s^{-k-1}\int_{\partial B_s\cap \Sigma}
                \frac{|(x-x_0)^N|^2}{|(x-x_0)^T|}\, .
        \notag
\end{align}
Integrating and applying the coarea formula once more gives the
claim.
\end{proof}

Notice that $(x-x_0)^N$ vanishes precisely when $\Sigma$ is
conical about $x_0$, i.e., when $\Sigma$ is invariant under
dilations about $x_0$. As a corollary, we get the following:

\begin{corollary}             \label{c:cmon}
Suppose that $\Sigma^k\subset \RR^n$ is a minimal submanifold and
$x_0\in \RR^n$; then the function
\begin{equation}                \label{e:thetadef}
        \Theta_{x_0}(s)
        =\frac{\Vol (B_s (x_0) \cap \Sigma)}{\Vol (B_s \subset \RR^k) }\,
\end{equation}
is a nondecreasing function of $s$.  Moreover, $\Theta_{x_0}(s)$
is constant in $s$ if and only if $\Sigma$ is conical about $x_0$.
\end{corollary}

Of course, if $x_0$ is a smooth point of $\Sigma$, then
$\lim_{s\to 0} \Theta_{x_0}(s) = 1$; the Allard regularity theorem
gives the converse of this.

The monotonicity of area is an very useful tool in the
regularity theory for minimal surfaces --- at least when there is
some {\it a priori} area bound.  For instance, this monotonicity
and a compactness argument allow one to reduce many regularity
questions to questions about minimal cones (this was a key
observation of W. Fleming in his work on the Bernstein problem;
see Section \ref{s:5}).  Similar monotonicity formulas have played
key roles in other geometric problems, including  harmonic maps,
Yang--Mills connections, J--holomorphic curves, and regularity
of limit spaces with a lower Ricci curvature bound.

\vskip2mm  Arguing as in Proposition \ref{p:pmonot}, we get a
weighted monotonicity:

\begin{proposition}
\label{p:meanvalue} If\hspace{2pt}  $\Sigma^k\subset \RR^n$ is a
minimal submanifold, $x_0\in \RR^n$, and $f$ is a function on
$\Sigma$, then
\begin{equation}
        t^{-k}\int_{B_t(x_0) \cap \Sigma} f-s^{-k}\int_{B_s(x_0) \cap \Sigma} f
\end{equation}
\begin{equation}
        =\int_{(B_t(x_0)\setminus B_s(x_0))\cap \Sigma} f\,
                \frac{|(x-x_0)^N|^2}{|x-x_0|^{k+2}}
        +\frac{1}{2} \int_{s}^t \tau^{-k-1}
        \int_{B_{\tau}(x_0)\cap \Sigma} (\tau^2-|x-x_0|^2)\,
                \Delta_{\Sigma} f\,d\tau
        \, . \notag
\end{equation}
\end{proposition}

We get immediately the following mean value inequality for the
special case of non--negative subharmonic functions:

\begin{corollary}
Suppose that $\Sigma^k\subset \RR^n$ is a minimal submanifold,
$x_0\in \RR^n$, and $f$ is a non--negative subharmonic function on
$\Sigma$;
 then
\begin{equation}                \label{e:defavg}
        s^{-k}  \int_{B_s(x_0) \cap \Sigma} f
\end{equation}
is a nondecreasing function of $s$. In particular, if $x_0\in
\Sigma$, then for all $s>0$
\begin{equation}
        f(x_0)\leq
        \frac{\int_{B_s (x_0) \cap \Sigma} f}
                {\Vol (B_s \subset \RR^k)} \, .
\end{equation}
\end{corollary}

\vskip2mm Another interesting (and crucial) appearance of
monotonicity is the volume comparison theorem of Bishop--Gromov
for manifolds with non--negative Ricci curvature, \cite{GLPa}; see
also \cite{Pe1} for a generalization to Ricci flow. In the case of
non--negative Ricci curvature, the monotonicity goes the opposite
direction:

\begin{theorem}     \label{t:bg}
If a $k$--dimensional manifold $M$ has non--negative Ricci
curvature, then
\begin{equation}    \label{e:bg}
\frac{\Vol (B_s (x_0) \subset M)}{\Vol (B_s \subset \RR^k) }
\end{equation}
is a non--increasing function of $s$.
\end{theorem}

\begin{proof}(Sketch)
 By the Laplacian comparison theorem, $\Delta_M r^2 \leq 2k$
where $r$ is the distance function to $x_0$. Integrating this by
parts gives
\begin{equation}
    2k \, \Vol (B_s (x_0)) \geq \int_{B_s (x_0)} \Delta r^2 = 2 s
    \, \int_{\partial B_s (x_0)} |\nabla r|  = 2 s \, \frac{d}{ds} \,
    \Vol (B_s (x_0)) \, ,
\end{equation}
where the last equality used the co--area formula (since $|\nabla
r| =1$ almost everywhere).  This differential inequality gives
\eqr{e:bg}.
\end{proof}

\begin{remark}
Equation \eqr{e:bg} immediately implies a volume doubling property
for manifolds with non--negative Ricci curvature:
\begin{equation}
    \Vol (B_{2s}(x_0) \subset M) \leq 2^n \, \Vol
    (B_{s}(x_0) \subset M) \, .
\end{equation}
 This property is very
useful for covering arguments, cf. \cite{G}.
\end{remark}

\vskip2mm We conclude this section with a well--known intrinsic
mean value inequality which is often useful but difficult to find
in the literature (it is often stated only for subharmonic
functions).

\begin{proposition}     \label{p:mv}
There exists $C= C(k)$ so that if $M$ is $k$--dimensional, $\Ric_M
\geq - (k-1) \,  s^{-2}$, and $u\geq 0$ satisfies $\Delta_M u \geq
-   s^{-2} \, u$, then
\begin{equation}
  u^2 (x) \leq \frac{C}{\Vol (B_s (x))} \, \int_{B_s(x)} u^2   \, .
\end{equation}
\end{proposition}

\begin{proof}
After rescaling the metric by $s$, it suffices to prove the case
$s=1$.  Let $N = M \times [-1,1]$ have the product metric, so that
$\Ric_N \geq - (k-1)$.  Define an auxiliary function $w$ on $N$ by
\begin{equation}
    w(x,t) = u(x) \, \e^{ t} \, .
\end{equation}
An easy calculation shows that
\begin{equation}
 \Delta_N w = \e^t \, \Delta_M u + \e^t \, u \geq 0 \, ,
\end{equation}
so that $w$ is subharmonic.  The mean value inequality for
subharmonic functions (see theorem $6.2$ on page $77$ of
\cite{ScYa1}) then gives
\begin{align}
  w^2 (x,0) &\leq \frac{C}{\Vol (B_1 (x,0) \subset N)} \, \int_{B_1 (x,0) \subset N} w^2
  \notag \\
  &\leq 2 \e^2 \, \frac{C}{\Vol (B_{1/2} (x) \subset M)} \, \int_{B_1 (x) \subset M} u^2 \, .
\end{align}
The proposition follows from this after we use   the
Bishop--Gromov volume comparison (cf. Theorem \ref{t:bg})  to
bound   $\Vol (B_{1} (x))/ \Vol (B_{1/2} (x))$.
\end{proof}

\section{The theorems of Bernstein and Bers}        \label{s:5}

 A classical theorem of S. Bernstein from 1916 says that
entire (i.e., defined over all of $\RR^2$) minimal graphs are
planes. This remarkable theorem of Bernstein was one of the first
illustrations of the fact that the solutions to a nonlinear PDE,
like the minimal surface equation, can behave quite differently
from the solutions to a linear equation.  Rather surprisingly,
this result very much depended on the dimension.  The combined
efforts of S. Bernstein \cite{Be}, E. De Giorgi \cite{DG}, F. J.
Almgren, Jr. \cite{Am1}, and J. Simons \cite{Sim} finally gave:

\begin{theorem}  If
$u:\RR^{n-1}\to \RR$ is an entire solution to the minimal surface
equation and $n\leq 8$, then $u$ is an affine function.
\end{theorem}

However, in 1969 E. Bombieri, De Giorgi, and E. Giusti \cite{BDGG}
constructed entire non--affine solutions to the minimal surface
equation on $\RR^8$ and an area--minimizing singular cone in
$\RR^8$. In fact, they showed that for $m \geq 4$ the cones
\begin{equation}        \label{e:cones1}
        C_m = \{ (x_1 , \dots , x_{2m} ) \mid  x_1^2 + \cdots + x_m^2 =
x_{m+1}^2 + \cdots + x_{2m}^2  \} \subset \RR^{2m} \,
\end{equation}
are area--minimizing (and obviously singular at the origin).

One way to prove the Bernstein theorem is to prove a curvature
estimate for minimal graphs.  The basic example is the estimate of
E. Heinz for surfaces:

\begin{theorem}   \cite{He}           \label{t:heinz}
If $D_{r_0} \subset \RR^2$ and
 $u:D_{r_0}\to \RR$ satisfies
the minimal surface equation, then for $\Sigma={\text{Graph}}_u$
and $0<\sigma\leq r_0$
\begin{equation}        \label{e:heinz}
        \sigma^2\,\sup_{D_{r_0-\sigma}}|A|^2\leq C\, .
\end{equation}
\end{theorem}

The original Bernstein Theorem  follows from Theorem \ref{t:heinz}
by taking $r_0 \to \infty$.   By the same reasoning, the examples
of \cite{BDGG} show that \eqr{e:heinz} cannot hold for all
dimensions.   However,  curvature estimates for graphs over
$B_{r_0} \subset \RR^{n-1}$   were proven in \cite{ScSiYa}  for
$n\leq 6$ (the remaining cases, i.e.,  $n=7$ and $8$, were proven
in \cite{Si2}).

\vskip2mm
 In contrast to the entire case, exterior solutions of the minimal graph
 equation, i.e., solutions on $\RR^2 \setminus B_1$, are much more
 plentiful. In this case,  Theorem
\ref{t:heinz}  only gives quadratic curvature decay $|A|^2 \leq C
\, |x|^{-2}$.  In particular, it is not even clear that $|\nabla
u|$ is bounded since $|\nabla |\nabla u|| \leq C \, |x|^{-1}$ is
not integrable along rays.  However, L. Bers proved that $\nabla
u$ actually has an asymptotic limit:

\begin{theorem} \cite{Ber} \label{t:bers}
If $u$ is a $C^2$ solution to the minimal surface equation on
$\RR^2 \setminus B_1$, then $\nabla u$ has a limit at infinity
(i.e., there is an asymptotic tangent plane).
\end{theorem}

 Bers' theorem was extended to higher
dimensions by L. Simon:

\begin{theorem} \cite{Si1} \label{t:berssi}
If $u$ is a $C^2$ solution to the minimal surface equation on
$\RR^n \setminus B_1$, then either \begin{itemize} \item $|\nabla
u|$ is bounded and $\nabla u$ has a limit at infinity. \item All
tangent cones at infinity are of the form $\Sigma \times \RR$
where $\Sigma$ is singular. \end{itemize}
\end{theorem}

Bernstein's theorem has had many other interesting
generalizations, including, e.g., curvature estimates of R. Schoen
for stable surfaces  and Schoen--Simon--Yau for stable
hypersurfaces with bounded density.  In the early
nineteen--eighties Schoen and Simon extended the theorem of
Bernstein to complete simply connected embedded minimal surfaces
in $\RR^3$ with quadratic area growth. A surface $\Sigma$ is said
to have quadratic area growth if for all $r>0$, the intersection
of the surface with the ball in $\RR^3$ of radius $r$ and center
at the origin is bounded by $C\, r^2$ for a fixed constant $C$
independent of $r$.  In corollary $1.18$ in \cite{CM4}, this was
generalized to quadratic area growth for {\it{intrinsic}} balls.

\section{Mean curvature flow}   \label{s:6}

Just as the Laplace equation has the heat equation as a parabolic
analog,  mean curvature flow is the parabolic analog of the
minimal surface equation.
 A one--parameter family of smooth
hypersurfaces $\{ M_t \} \subset \RR^{n+1}$ {\it flows by mean
curvature} if
\begin{equation}
    z_t = {\bf{H}} (z) = \Delta_{M_t} z \, ,
\end{equation}
where $z$ are coordinates on $\RR^{n+1}$ and ${\bf{H}} = - H \nn$
is the mean curvature vector.

\begin{example}     \label{e:ss}
Let $M_t$ be the family of concentric shrinking $n$--spheres of
radius $\sqrt{R^2 -2nt}$.  It is easy to see that $M_t$ flows by
mean curvature and is smooth up to $t=R^2/(2n)$ when it shrinks to
a point.
\end{example}

 Suppose now that each $M_t$ is
  the graph of a function $u(\cdot , t)$.  So, if $z =
(x,y)$ with $x \in \RR^n$, then $M_t$ is given by $y = u (x,t)$
which satisfies
\begin{equation}    \label{e:evq}
    u_t =  (1 + |du|^2)^{1/2} \,
    \dv \left( \frac{du}{(1 + |du|^2)^{1/2} } \right)
    \, ,
\end{equation}
where $du$ is the $\RR^n$ gradient of the function $u$ and $\dv$
is divergence in $\RR^n$.

The monotonicity formula and mean value inequality of Section
\ref{s:4} have analogs in mean curvature flow.  The monotonicity
formula, proven by G. Huisken (and extended to more general weak
solutions by T. Ilmanen and B. White), is:

\begin{theorem} \cite{H} \label{t:mch}
If a smooth one--parameter family of hypersurfaces $M_t$ flows by
mean curvature in $\RR^{n+1} \times [-T,0]$, then
\begin{equation}
    \frac{d}{dt} \, \int_{M_t} \left(-4\pi t\right)^{-n/2} \,
    \e^{ \frac{|x|^2}{4t} } \, =  - \int_{M_t}
    \left| {\bf{H}} - \frac{x^N}{2t} \right|^2 \, \left(-4\pi t\right)^{-n/2} \,
    \e^{ \frac{|x|^2}{4t} }
     \, .
\end{equation}
In particular, the ``density ratio'' $ \int_{M_t} \left(-4\pi
t\right)^{-n/2} \,
    \e^{ \frac{|x|^2}{4t} }$ is non--increasing.
\end{theorem}

The restrictions of the coordinate functions to $M_t$ satisfy the
heat equation
\begin{equation}    \label{e:heqn}
    \partial_t x_i =  \Delta_M x_i \, .
\end{equation}
 From this, we see that the restriction of $|x|^2$ satisfies
$(\partial_t - \Delta_M ) \, |x|^2 = -2n$.  As in the stationary
case (i.e., for minimal surfaces),  this is the key to the proof
of Theorem \ref{t:mch}.

The mean value inequality in this case applies to non--negative
solutions of the heat equation on $M_t$; we refer to \cite{E1} for
more detail on this as well as discussion of the local
monotonicity formula for mean curvature flow proven in \cite{E2}.

\vskip2mm The parabolic maximum principle has been very
useful in mean curvature flow (somewhat similarly to the convex
hull property for minimal surfaces).  Two immediate, but useful,
consequences are:
\begin{enumerate}
\item
Disjoint surfaces stay disjoint.
\item
An embedded surface stays embedded (as long as it evolves
smoothly).
\end{enumerate}
 The reason for (1) and (2) is  quite simple.  Suppose that two
 initially disjoint surfaces touch at a
 first time $t$ at a point $x$.  Clearly, they will be tangent at
 $(x,t)$ so nearby we see two graphs, one above the other.  Hence,
  at $x$ the mean curvature of the upper graph
 is larger (or equal to).  Therefore,
 the upper graph crossed over the lower at a slightly early time,  contradicting that $t$
 is first time of contact.

Combining (1) with the shrinking spheres of Example \ref{e:ss}, we
see that any compact hypersurface flowing by mean curvature has a
finite extinction time.

\section{Ricci flow}        \label{s:7}

The Ricci flow is the parabolic analog of the Einstein equation
$\Ric_M = {\text{ Constant}} \, g$, where $g$ is the metric.
Namely, let $M^3$ be a fixed smooth manifold and let $g(t)$ be a
one--parameter family of metrics on $M$ evolving by the Ricci
flow, so
\begin{equation}  \label{e:eqRic}
 \partial_t g=-2\,\Ric_{M_t}\, .
\end{equation}
Short--time existence for the Ricci flow was established by
Hamilton:

\begin{theorem} \cite{Ha2}  \label{t:shorttime}
Given any smooth compact Riemannian  manifold $(M,g_0)$, there
exists a unique smooth solution $g(t)$ to \eqr{e:eqRic} with
initial condition $g(0) = g_0$ on some time interval
$[0,\epsilon)$.
\end{theorem}

Long--time existence is quite a bit more subtle, see \cite{Ha2}
and \cite{Pe1}.

\vskip2mm There are many formal similarities between the Ricci
flow and the mean curvature flow, including similarities between
the evolution equations for various geometric quantities.  One
interesting distinction is the evolution equation for the scalar
curvature $R= R(t)$ under the Ricci flow (see, for instance, page
16 of \cite{Ha3})
\begin{equation}    \label{e:prescalar}
    \partial_t R = \Delta R + 2 |\Ric|^2   \geq \Delta R +
    \frac{2}{n} \, R^2 \, ,
\end{equation}
where the inequality used the Cauchy--Schwarz inequality ($M$ is
$n$--dimensional).   This differential inequality has an
interesting consequence: After flowing for any positive amount of
time, there is a lower bound for the scalar curvature (the mean
curvature flow has no such analog). To make this precise,  a
straightforward maximum principle argument gives that at time $t >
0$
\begin{equation}    \label{e:scalar}
   R(t) \geq  \frac{1}{1/[ \min R(0)] -2t/n} =  -\frac{n}{2 (t+C)} \, .
\end{equation}
In the derivation of \eqr{e:scalar} we implicitly assumed that
$\min R(0)<0$.  If this was not the case, then \eqr{e:scalar}
trivially holds with $C=0$, since, by \eqr{e:prescalar}, $\min R
(t)$ is always non--decreasing.

\section{Gradient estimates}        \label{s:8}

Gradient estimates  have played a key role in geometry and pde
since at least the early work of Bernstein.  These are probably
the most fundamental a priori estimates for elliptic and parabolic
equations, leading to Harnack inequalities, Liouville theorems,
and compactness theorems for both linear and nonlinear pde.

 A typical example for
linear equations is the well--known gradient estimate of S.Y.
Cheng and S.T. Yau for harmonic functions:

\begin{theorem}  \cite{CgYa}   \label{t:cy}
If $\Delta u = 0 $   on  $B_r(0)$ with non--negative Ricci
curvature, then
\begin{equation}    \label{e:cy}
    |\nabla u|(0) \leq C \, r^{-1} \,
    \|u \|_{\infty} \, ,
\end{equation}
 where
$\|u\|_{\infty}$ is the sup norm of the function $u$ on $B_r(0)$.
\end{theorem}

To give something of the flavor, we will use the maximum principle
to prove  \ref{t:cy} on the Euclidean unit ball   $B_1 (0) \subset
\RR^n$.

\begin{proof} (of Theorem \ref{e:cy} for $B_1 (0) \subset \RR^n$.)
Define the cutoff function $\eta (x) = 1 - |x|^2$, so that
$|\nabla \eta| \leq 2$ and $\Delta \eta = -2n$.  We compute that
\begin{align}    \label{e:cy1}
    \Delta (\eta^2 \, |\nabla u|^2) &\geq -2n \, |\nabla u|^2 - 16
    \, \eta \, |\nabla u| \, |\Hess_u| + 2 \, \eta^2 \, |\Hess_u
    |^2 \notag \\ &\geq - (2n+32) \, |\nabla u|^2 \, ,
\end{align}
where the last inequality used the absorbing inequality $16ab \leq
2a^2 + 32 b^2$.  In particular, the function $w = (n+16) \, u^2 +
\eta^2 \, |\nabla u|^2$ is subharmonic on $B_1(0)$ (i.e., $\Delta
w \geq 0$).  By the maximum principle, the maximum of $w$ occurs
on the boundary so that
\begin{equation}    \label{e:cy2}
     |\nabla u|^2 (0) \leq w(0) \leq \max_{\partial B_1(0)} w =
     (n+16) \, \max_{\partial B_1(0)} u^2  \, .
\end{equation}
\end{proof}

In fact, Cheng and Yau prove a stronger estimate: If in addition
$u$ is  positive on $B_{r}(0)$, then
\begin{equation}     \label{e:pge}
    |\nabla \log u| (0) \leq C \, r^{-1} \, .
\end{equation}
 An important
consequence is the Harnack inequality for positive harmonic
functions

\noindent {\bf{Elliptic Harnack inequality}}:
\begin{equation}    \label{e:harn}
    \sup_{B_{r/2}(0)} u \leq C' \, \inf_{B_{r/2}(0)} u \, .
\end{equation}

\begin{proof}
Suppose the sup and inf are achieved at  $p , q \in \partial
B_{r/2}(0)$.  Fix a curve $\gamma_{p,q} \subset B_{r/2}(0)$ from
$p$ to $q$ of length at most $r$ (e.g., connect each point to $0$
by a ray). Integrating the bound $\sup_{B_{r/2}(0)}|\nabla \log u|
\leq 2 C \, r^{-1}$ over $\gamma_{p,q}$ gives
\begin{equation}
    \log \frac{u(p)}{u(q)} \leq \int_{\gamma_{p,q}} |\nabla \log u| \leq
    2 \, C  \, .
\end{equation}
\end{proof}

This gradient estimate also gives the global Liouville theorem of
Yau, \cite{Ya3}:

\noindent {\bf{Liouville theorem}}: If $u$ is a positive harmonic
function on a complete manifold with non--negative Ricci
curvature, then $u$ is constant.

\begin{proof}
We can take $r \to \infty$ in \eqr{e:pge} to get that $|\nabla u|
= 0$.
\end{proof}

\vskip2mm The parabolic analog of Theorem \ref{t:cy} is the
gradient estimate for the heat equation of P. Li and Yau (we will
state the  version for $M$ complete):

\begin{theorem}   \cite{LiYa}   \label{t:ly}
If $u$ is a positive solution of $\partial_t u = \Delta u$ for
$0\leq t$ on a complete $M$ with non--negative Ricci curvature,
then
\begin{equation}    \label{e:lyh}
    2t \left( |\nabla \log u|^2 - \partial_t \log u \right)
\leq n
     \, ,
\end{equation}
or, equivalently, $2t \, \Delta \log u \geq - n$.
\end{theorem}

\begin{proof}
(Sketch) Set
\begin{equation}
    w = - t \, \Delta \log u= t \left( |\nabla \log u|^2 -
\partial_t \log u \right) \, .
\end{equation}
  The key calculation is (cf. lemma $1$ on page $155$ of
\cite{ScYa1})
\begin{equation}    \label{e:keyly}
    t \, ( \Delta - \partial_t ) \, w \geq \frac{2}{n} \, w^2 -
    w - 2t \, \nabla w \cdot \nabla \log u \, .
\end{equation}
Suppose that $w$ achieves its maximum on $M \times [0,t]$ at
$(x,t)$ (for example, when $M$ is compact; otherwise we use a
cutoff).  The parabolic maximum principle then gives  $\nabla w
(x,t) = 0$ and $(\Delta -
\partial_t ) \, w (x,t) \leq 0$.  Substituting this into
\eqr{e:keyly} gives
\begin{equation}    \label{e:keyly2}
    0 \geq \frac{2}{n} \, w^2 (x,t) -
    w  (x,t) \, ,
\end{equation}
so that $w(x,t) \leq n/2$ as desired.
\end{proof}

Integrating this along curves as in the elliptic case gives for
$t_1 < t_2$ that

\noindent {\bf{Parabolic Harnack inequality}}:
\begin{equation}    \label{e:pharn}
     u(x_1 , t_1) \leq u(x_2 , t_2) \, \left( \frac{t_2}{t_1} \right)^{\frac{n}{2}} \,
        \e^{ \frac{\dist^2 (x_1 , x_2)}{4(t_2 - t_1)} } \, .
\end{equation}
In \cite{Ha1}, R. Hamilton gave an extension of \eqr{e:lyh} to a
full matrix estimate whose trace was \eqr{e:lyh}.  For example, if
$u$ is positive solution of the heat equation on $\RR^n \times
[0,T]$, then \cite{Ha1} implies that
 \begin{equation}  \label{e:lha}
    2 t \, \Hess_{\log u} +   \delta_{ij}  \geq 0 \, .
\end{equation}
Taking the trace of \eqr{e:lha} gives  $ |\nabla \log u|^2 -
\partial_t \log u \leq n/(2t)$.

\subsection{Gradient estimates for nonlinear equations}
 For the (nonlinear) minimal surface equation, the situation is somewhat different.
In this case, i.e., when the graph of $u$ is minimal on $B_r(0)$,
then Bombieri, De Giorgi, and M. Miranda proved in
 \cite{BDM} that
\begin{equation}
    \log |du|(0) \leq   C \, ( 1 + r^{-1} \,
    \|u \|_{\infty} ) \, ,
\end{equation}
where $du$ is the $\RR^n$ gradient of the function $u$  (the case
of surfaces was done by R. Finn in \cite{Fi}).  By an earlier
example of Finn, this exponential dependence cannot be improved.
In \cite{K}, N. Korevaar gave a maximum principle proof of a
weaker form of \cite{BDM}; this weaker form had $\|u\|^2_{\infty}$
in place of $\|u\|_{\infty}$.

In \cite{CM12}, we proved a sharp gradient estimate for graphs
flowing by mean curvature:

\begin{theorem}    \cite{CM12}  \label{t:gbp}
There exists $C= C(n)$ so if the graph of $u: B_{\sqrt{2n+1} r}
 \times [0,r^2] \to \RR$ flows by mean
curvature, then
\begin{equation}    \label{e:gbp}
    \log |du|(0,r^2/[4n]) \leq   C \, ( 1 + r^{-1} \,
\|u(\cdot , 0) \|_{\infty} )^2 \, .
\end{equation}
\end{theorem}

The quadratic dependence on  $\|u(\cdot , 0) \|_{\infty}$ in
\eqr{e:gbp} should be compared with the linear dependence which
holds when the graph of $u$ is minimal (i.e., $u_t = 0$). Somewhat
surprisingly,  examples in \cite{CM12} show that this quadratic
dependence on $\| u (\cdot , 0) \|_{\infty}$ is sharp.

The first gradient estimate for mean curvature flow was proven by
Ecker and Huisken who adapted Korevaar's argument to mean
curvature flow in theorem 2.3 of \cite{EH2} to get
\begin{equation}    \label{e:eh2}
    \log |du|(0,r^2/[4n]) \leq  1/2 \, \log
    \left( 1 + \|du (\cdot,0)\|^2_{\infty} \right)
      + C \, ( 1 + r^{-1} \, \|u(\cdot , 0) \|_{\infty} )^2 \, .
\end{equation}
Note that, unlike \eqr{e:gbp}, the gradient bound \eqr{e:eh2}
depends also on the initial bound for the gradient.

\subsection{Generalizations}        \label{ss:72}

The Harnack inequality actually holds for much more general
spaces. For instance, L. Saloff--Coste and A. Grigor'yan (see
\cite{SC} and \cite{Gr}) have shown that the following two
properties suffice
\begin{description}
\item[Volume doubling]
There exists $C_D$ so that
\begin{equation}
    \Vol (B_{2r}(x)) \leq C_D \, \Vol (B_r(x)) \, ,
\end{equation}
 for all $r>0$ and points $x$.
\item[Neumann Poincar\'e inequality]
There exists $C_N$ so that if $\int_{B_r(x)} f = 0$, then
\begin{equation}
    \int_{B_r(x)} f^2 \leq C_N \, r^2 \int_{B_r(x)} |\nabla f|^2 \, ,
\end{equation}
 for all $r>0$ and points $x$.
\end{description}
These properties, however, do not imply the gradient estimate.
Note that manifolds with non--negative Ricci curvature  satisfy
both conditions (the Poincar\'e inequality essentially follows
from \cite{Bu}, cf. also \cite{Je}).

The De Giorgi, Nash, Moser theory (see chapter $8$ in \cite{GiTr}
or section $4.4$ in \cite{HnLn}) gives a Harnack inequality as
long as we have a volume doubling and a Sobolev inequality.  The
difference between a Sobolev and Poincar\'e  inequality is that a
Sobolev controls an $L^p$--norm of $f$ for some $p > 2$.
Surprisingly, in \cite{HzKo}, P. Hajlasz and P. Koskela showed
that the volume doubling and Neumann Poincar\'e inequality
together imply a Sobolev inequality, thereby recovering the above
result of Saloff--Coste and Grigor'yan.  See \cite{ChC2},
\cite{ChC3}, \cite{Hj} for more such ``low regularity'' analysis
(including analysis on singular spaces).

\section{Simons type inequalities}  \label{s:9}

In this section, we recall a very useful differential inequality
for the Laplacian of the norm squared of the second fundamental
form of a minimal hypersurface $\Sigma$ in $\RR^n$ and illustrate
its role in a priori estimates. This inequality, originally due to
J. Simons (see \cite{CM1} for a proof and further discussion), is:

\begin{lemma}  \cite{Sim}   \label{l:simonsine}
If $\Sigma^{n-1}\subset \RR^n$ is a minimal hypersurface, then
\begin{equation}    \label{e:simtype}
\Delta_{\Sigma} \, |A|^2  = - 2\, |A|^4 + 2 |\nabla_{\Sigma} A|^2
\geq - 2 \, |A|^4 \, .
\end{equation}
\end{lemma}

An inequality of the type \eqr{e:simtype} on its own does not lead
to pointwise bounds on $|A|^2$ because of the nonlinearity.
However, it does lead to estimates if a   ``scale--invariant
energy'' is small.  For example, H. Choi and Schoen used
\eqr{e:simtype}  to prove:

\begin{theorem} \cite{CiSc}     \label{t:cisc}
If $0 \in \Sigma \subset B_r (0)$ with $\partial \Sigma \subset
\partial B_r (0)$ is a minimal surface with sufficiently small
total curvature $\int |A|^2$, then $|A|^2(0) \leq r^{-2}$.
\end{theorem}

Analogs of \eqr{e:simtype} occur in a surprising number of
geometric problems.
 For example, when $u:M^m \to N^n$ is a
harmonic map, the energy density $|du|^2$ satisfies this type of
inequality,
 leading to an {\it a priori} estimate when $u$ has small
 scale--invariant energy $r^{2-m} \, \int_{B_r} |du|^2$
(see \cite{SaUh}, \cite{Sc2}).  Similar inequalities hold for the
curvature of a Yang--Mills connection or the curvature of an
Einstein manifold.  When $M$ is an Einstein manifold, its
curvature tensor $R$ satisfies (see \cite{Ha2} or equation (2.6)
in \cite{An1}; \cite{Ha2} also establishes a parabolic analog for
the Ricci flow)
\begin{equation}    \label{e:simtypeHa}
    \Delta_{M} \, |R|   \geq - C \,  |R|^2  \, .
\end{equation}
We next use \eqr{e:simtypeHa} to prove an estimate for Einstein
manifolds (cf. lemma $2.1$ in \cite{An1}).  For simplicity, we
restrict to the case $\Ric_M = 0$.

\begin{theorem}     \label{t:einstein}
There exist $\epsilon =\epsilon (n)>0$, such that if $M^n$ is an
$n$--dimensional Ricci--flat (Einstein) manifold and for some
$x\in M$ either
\begin{equation}
\Vol( B_r(x))\geq [1-\epsilon] \, \Vol (B_r\subset \RR^n)
\end{equation}
or
\begin{equation}    \label{e:and2}
\int_{B_r(x)}|R|^{n/2}< \epsilon \, r^{-n} \, \Vol( B_r(x))  \, ,
\end{equation}
then $|R| (x) \leq C \, r^{-2}$.
\end{theorem}

\begin{proof}
We will prove that \eqr{e:and2} gives the pointwise curvature
bound; the other case is similar.  Set $F(z)=(r  - 4 \dist_M
(x,z)) \,|R|^{1/2} (z)$, so that
\begin{equation}
F(x) = r \, |R|^{1/2} (x)   \geq 0\, , \text{ and } F\left.
\right|_{B_{r}(x)\setminus B_{r/4}(x)}\leq 0\, .
\end{equation}
Therefore, it suffices to prove that $F \leq C$ for some fixed
constant $C$.

We will assume that $\max_{B_r(x)} F > 32$ and deduce a
contradiction if $\epsilon > 0$ is sufficiently small.
 Let $y$ be a point where the maximum of $F$ is achieved   and
set $s=  |R|^{-1/2}(y)$.  Since $F(y) > 32$, we have $32 s < |r  -
4 \dist_M (x,y)|$ so that for $z \in B_s(y)$
\begin{equation}
1/2 \leq  \frac{|r  - 4 \dist_M (x,y)|}{|r -4 \dist_M (x,z)|} \leq
2  \, .
\end{equation}
 Since $F(z) \leq F(y)$, it follows that
$B_s(y)$ satisfies
\begin{equation}    \label{e:bup}
    \sup_{B_s (y)} |R|^{1/2} \leq 2 |R|^{1/2} (y) = 2/ s \, ,
\end{equation}
so that \eqr{e:simtypeHa} gives on $B_s (y)$ that
\begin{equation}    \label{e:sand}
    \Delta_{M} \, |R|   \geq - C \,  |R|^2 \geq - 4\, C \,  s^{-2} \,
          |R|    \, .
\end{equation}
 Furthermore, the Bishop--Gromov volume comparison , i.e., Theorem
 \ref{t:bg},
gives
\begin{equation}
    \frac{\Vol( B_s(y)) }{ s^n } \geq \frac{ \Vol( B_{r/2}(y)) }{ (r/2)^n }
    \geq 2^{-n} \, \frac{\Vol( B_r(x))
}{ r^n }\, .
\end{equation}
It follows from this and \eqr{e:and2} that $B_s(y)$ also satisfies
\begin{equation}    \label{e:and3}
\int_{B_s(y)}|R|^{n/2}< 2^n \, \epsilon \, \, \frac{\Vol( B_s(y))
}{ s^n }\, .
\end{equation}
Using \eqr{e:bup}, \eqr{e:sand}, and \eqr{e:and3},  the mean value
inequality (Proposition \ref{p:mv}) gives
\begin{equation}
     s^{-n} = |R|^{n/2} (y) \leq \frac{C}{\Vol(
     B_s(y)) } \,  \int_{B_s(y)} |R|^{n/2}    < C \, 2^n \, \epsilon \, s^{-n} \, .
\end{equation}
This gives a contradiction for $\epsilon$ sufficiently small.
\end{proof}

\begin{remark}
We could alternatively have proven Theorem \ref{t:einstein} by
integral methods, i.e.,  using Moser iteration.  However, the
above proof by scaling is both shorter and more elementary.
\end{remark}

Finally, we mention that \eqr{e:simtype} has a parabolic version
as well (see proposition $2.15$ in \cite{E1}): If $M_t$ flows by
mean curvature, then
\begin{equation}    \label{e:simtypeH}
        \left( \frac{\partial}{\partial t} - \Delta_{M_t} \right) \, |A|^2
        = 2\, |A|^4 - 2 |\nabla_{M_t} A|^2  \,        .
\end{equation}
As in the elliptic case, this Simons' type inequality is a crucial
ingredient for establishing curvature estimates.

\section{Minimal annuli with small total curvature are graphs}  \label{s:10}

It is easy to see that a minimal disk with small total curvature
must be a graph  away from its boundary:

\begin{quotation}
 If $\int_{B_R \cap \Sigma} |A|^2 < \epsilon$, then Theorem
\ref{t:cisc} gives $|A|^2 < C \, \epsilon / R^2$ on $B_{R/2} \cap
\Sigma$; integrating this (since $|\nabla \nn| \leq |A|$) implies
that each component of $B_{R/2} \cap \Sigma$ is a graph if
$\epsilon > 0$ is small enough.
\end{quotation}

\vskip1mm \noindent However, the corresponding question for
minimal annuli is more subtle.   We shall discuss this and some
related problems in this section.

In \cite{CM10}, we gave three proofs that a minimal annulus with
small total curvature is a graph.  The first used a singular
integral formula which had previously been useful for estimating
nodal and singular sets; see Proposition \ref{p:annulus} below and
compare \cite{Do}.  The second, and easiest, applies more
generally to surfaces with quasi--conformal Gauss maps; see
Proposition \ref{p:annulus3}.  The third, which is outlined in
Lemma \ref{l:supe}, was the one which could be extended to
``annuli with slits'' --- i.e., embedded double--valued minimal
graphs.

\vskip2mm In this section, $\Sigma \subset \RR^3$ is a compact
connected oriented immersed surface.  If $a \in \SS^2$,
$a^{\perp}$ denotes $\{ x \in \RR^3 \, | \, \langle x , a \rangle
= 0 \}$. For  $a , b \in \SS^2$, $\Angle (a,b)$ is the angle
between  $a^{\perp} , b^{\perp}$; i.e.,

\begin{equation}
    \Angle (a,b) = \dist_{\SS^2} (a, \{ b , -b \}) \, .
\end{equation}

Let $f$ be  harmonic  on $\Sigma^2$ with
 critical points $\{ y_i \}$ with multiplicities $\{m_i\}$.
Suppose that none of the $y_i$'s lie on $\partial \Sigma$. The
Bochner formula on $\Sigma\setminus \{y_i\}$ gives
\begin{equation}  \label{e:bochner2}
\Delta_{\Sigma} \,\log |\nabla_{\Sigma} f|^2
=2\frac{|\Hess_f|^2}{|\nabla_{\Sigma} f|^2}+2\,K
-\frac{|\nabla_{\Sigma} |\nabla_{\Sigma} f|^2|^2}{|\nabla_{\Sigma}
f|^4} =2\,K \, .
\end{equation}
Here we used that since $\Delta_{\Sigma}f=0$ and $\Sigma$ is
$2$--dimensional, then
\begin{equation}
    2\,|\Hess_f|^2\,|\nabla_{\Sigma} f|^2 =|\nabla_{\Sigma}
|\nabla_{\Sigma} f|^2|^2 \, .
\end{equation}
 Hence, by Stokes' theorem
\begin{align}    \label{e:formula}
\int_{\partial \Sigma}\frac{d\,\log |\nabla_{\Sigma} f|^2}{dn}
&=\int_{\Sigma\setminus \{y_i\}}\Delta_{\Sigma} \,\log
|\nabla_{\Sigma} f|^2+4\,\pi\sum_i m_i \notag \\
&=2\int_{\Sigma}K +4\,\pi\sum_i m_i\, .
\end{align}

\begin{proposition}  \cite{CM10} \label{p:annulus}
If $\Sigma$ is connected and
 minimal with boundaries
$\sigma_1$ and $\sigma_2$, $\int_{\sigma_1\cup \sigma_2} \, |A| <
\pi / 8$, and $\int_{\Sigma}K \geq -\pi$, then $\Sigma$ is
graphical.
\end{proposition}

\begin{proof}
Fix $q_i \in \sigma_i$. Since $|\nabla_{\Sigma} \dist_{\SS^2} (\nn
(q_i) , \nn ( \cdot  ) ) | \leq |A|$, the assumption
 on $\partial \Sigma$ gives
\begin{equation} \label{e:nsma}
\sum_{i} \sup_{z_i \in \sigma_i} \dist_{\SS^2} (\nn (q_i) , \nn
(z_i)) \leq \sum_{i} \int_{\sigma_i} |A| < \pi / 8 \, .
\end{equation}
Choose $b \in \SS^2$ with $\Angle ( \nn(q_i) , b )  \leq \pi /4$
for $i=1,2$. We will show that
 $\Sigma$ is
graphical over the plane $b^{\perp}$. By the triangle inequality
and \eqr{e:nsma}, for $i=1,2$,
\begin{equation} \label{e:nsma2}
\sup_{z_i \in \sigma_i} \Angle (b , \nn (z_i)) \leq \pi / 4 +
\sup_{z_i \in \sigma_i} \dist_{\SS^2} (\nn (q_i) , \nn (z_i))
                                            <  3 \, \pi / 8 \, .
\end{equation}
 Rotate coordinates so that $b=(0,0,1)$ and
$b^{\perp}$ is the $x_1$--$x_2$--plane.
  Fix $\theta$ and set
\begin{equation}
        f= x_1 \cos \theta + x_2 \sin \theta \, .
\end{equation}
Given $x \in \Sigma$,
\begin{equation} \label{e:ala}
|\nabla_{\Sigma} f|^2 (x) = 1 - \langle (\cos \theta , \sin \theta
, 0 ) , \nn (x) \rangle^2 \geq \langle b , \nn(x) \rangle^2
  \, .
\end{equation}
On $\partial \Sigma = \sigma_1 \cup \sigma_2$, \eqr{e:nsma2}  and
\eqr{e:ala} imply that
\begin{equation} \label{e:that}
\inf_{\partial \Sigma} |\nabla_{\Sigma} f| \geq \inf_{\partial
\Sigma} \, | \langle b , \nn(x) \rangle |
> \cos (3 \, \pi / 8) > 1/3 \, .
\end{equation}
Since $|\nabla_{\Sigma} |\nabla_{\Sigma} f|| \leq |\Hess_f| \leq
|A|$, \eqr{e:that} gives on $\partial \Sigma$
\begin{equation}
| \nabla_{\Sigma} \log |\nabla_{\Sigma} f|^2| = \frac{ 2 \,
|\nabla_{\Sigma} |\nabla_{\Sigma} f|| }{ |\nabla_{\Sigma} f|} \leq
6 \, |A| \, . \label{e:this}
\end{equation}
 Integrating \eqr{e:this}, we   get
\begin{equation}    \label{e:formula1}
\int_{\partial \Sigma} \left| \frac{d\,\log |\nabla_{\Sigma}
f|^2}{dn} \right| \leq 6 \, \int_{\partial \Sigma} |A| < 3 \, \pi
/ 4 \, .
\end{equation}
Since $\Sigma$ is minimal, $\Delta_{\Sigma} f = 0$. Substituting
\eqr{e:formula1} into \eqr{e:formula},
\begin{equation}
4\,\pi\sum_i m_i = \int_{\partial \Sigma}\frac{d\,\log
|\nabla_{\Sigma} f|^2}{dn} - 2\int_{\Sigma}K <  3 \, \pi / 4 + 2
\pi < 4 \, \pi \, ;
\end{equation}
hence, $f$ has no critical points.  Since this is true for any
$\theta$,  $\Sigma$ is graphical over $b^{\perp}$.
\end{proof}

In fact, Proposition \ref{p:annulus} holds for $\Sigma$ whose
Gauss map  $\nn$ is quasi--conformal:

\begin{proposition}  \cite{CM10} \label{p:annulus3}
If $\Sigma\subset \RR^3$ is connected, $\int_{\Sigma} \, | K |
\leq \pi$,
 $|A|^2 \leq C \, |K|$, and $\partial \Sigma$ has
 components
$\{ \sigma_i \}_{1 \leq i \leq n}$ with
\begin{equation}
 \sum_{i=1}^n \inf_{a\in \SS^2} \sup_{x\in \sigma_i} \, \{ \dist_{\SS^2} (\nn
(x),a )\} <\epsilon < \pi / 8 \, ,
\end{equation}
then $\nn (\Sigma) \subset \cB_{2 \, \epsilon}(a)$ for some $a \in
\SS^2$ and $\Sigma$ is the graph of $u$ over $a^{\perp}$ with
$|\nabla u| \leq 4 \, \epsilon$.
\end{proposition}

\subsection{Holomorphic functions on annuli} \label{s:sefa}

We next   estimate the oscillation of a   holomorphic function
with small gradient on annuli.  Since the Gauss map of a  minimal
annulus is conformal and its derivative is the second fundamental
form $A$, this serves as a good model (in fact, the proof can
easily be adapted to that case).

Suppose that $f$ is holomorphic on an annulus $D_{R} \setminus
D_{\delta}$ and satisfies $|\nabla f| \leq \epsilon / (2\pi \,
|z|)$.  Integrating this around each circle, we see that the
oscillation of $f$ on each circle is at most $\epsilon$.  However,
if we simply integrate this bound radially to compare $f$ on the
different circles, we get the log of the ratio of the radii.  The
next lemma improves on this to get a bound independent of this
ratio; the key is to keep track of cancellation rather than
integrating the bound on $|\nabla f|$.

\begin{lemma} \label{l:supe}
If  $f:D_{R} \setminus D_{\delta} \to \CC$
 is holomorphic and
$\int_{\partial D_{R} \cup \partial D_{\delta}} |\nabla f| \leq
\epsilon$, then
\begin{equation} \label{e:supe}
\min_{c \in \CC} \, \max_{z}   |f(z) - c| \leq \epsilon \, .
\end{equation}
\end{lemma}

\begin{proof}
 For $\delta \leq s \leq R$, define the circular average of $f$ by
\begin{equation} \label{e:iofs}
I(s) = (2 \, \pi \, s)^{-1} \, \int_{\partial D_s} f = (2 \,
\pi)^{-1} \, \int_0^{2 \, \pi} f(s \, \e^{i \, \theta}) \, d
\theta \, ,
\end{equation}
and $c = I(\delta)$.  Differentiating \eqr{e:iofs}, we have
\begin{align} \label{e:iofsp}
2 \, \pi \, s \, I'(s) &= \int_{\partial D_s} \frac{\partial
f}{\partial r} =  - i \, s^{-1} \, \int_{\partial D_s}
 \frac{\partial f}{\partial \theta}
= -i \,  \int_0^{2 \, \pi} \frac{\partial f}{\partial \theta}  \,
d \theta
\notag \\
&= -i \,  [ f(s \, \e^{i \, 2 \, \pi}) - f(s \, \e^0)] = 0  \, ,
\end{align}
where we used that $\frac{\partial f}{\partial r} = -i \, r^{-1}
\, \frac{\partial f}{\partial \theta}$ since $f$ is holomorphic.
In particular, $I(R) = c$.  Since $I(s)$ is the average of $f$
over $\partial D_s$, there exist $y_1 , y_2 \in \partial D_{R}$
with
\begin{equation} \label{e:avg1}
c = \rp (f(y_1)) + i \,  \ip (f(y_2))  \, .
\end{equation}
Combining \eqr{e:avg1} with $\int_{\partial D_{R} \cup \partial
D_{\delta}} |\nabla f| \leq  \epsilon$,
\begin{equation} \label{e:avg2}
\max_{y \in \partial D_{R}} \left|  \rp (f(y)-c) \right| \leq
\epsilon  / 2
 {\text{ and }} \max_{y \in \partial D_{R}}
\left|  \ip (f(y)-c) \right| \leq \epsilon / 2  \, ,
\end{equation}
so that $|f-c| \leq  \epsilon$ on $\partial (D_{R} \setminus
D_{\delta})$.   The maximum principle then gives \eqr{e:supe}.
\end{proof}

Note that Lemma \ref{l:supe} does not hold for  harmonic functions
(in particular, \eqr{e:iofsp} does not hold); e.g., take $\epsilon
\, \log r / (4 \, \pi)$ and $R> \e^{8\, \pi} \, \delta$.

\subsection{Bers' Theorem revisited}

We have seen that an annulus with small total curvature must be
graphical, even as the outer radius goes to infinity.  However,
Bers' Theorem indicates that much more should be true: the unit
normal actually goes to a limit at infinity.  This follows from a
sharper decay estimate for the curvature.  To keep things simple,
we will give the argument for a holomorphic function with small
gradient as in the previous lemma.

\begin{lemma} \label{l:supeB}
If  $f: D_R \setminus D_{1} \to \CC$
 is holomorphic and
$|\nabla f(z)| \leq 1 /  |z|$, then
\begin{equation} \label{e:supeB}
     \int_{
D_{R^{1/2} \e^{t}} \setminus D_{R^{1/2} \e^{-t}} } |\nabla f|^2
\leq 2 \pi \, \e^{2t} / R \, .
\end{equation}
\end{lemma}

\begin{proof}
By  \eqr{e:iofsp}, we can subtract a constant so that
$\int_{\partial D_r} f = 0$ for all $r$.  Set
\begin{equation}        \label{e:Et}
    E(t) = \int_{ D_{R^{1/2} \e^{t}} \setminus D_{R^{1/2}
    \e^{-t}} } |\nabla f|^2 \, .
\end{equation}
 Since $f$ is holomorphic, the Cauchy--Riemann equations give
\begin{equation}
    d \bar{f} \wedge d f = i \, |\nabla f|^2 \,
    dx\wedge dy \, .
\end{equation}
Applying Stokes' theorem, we get
\begin{equation}        \label{e:qs}
        E(t) = -  i  \, \int_{ D_{R^{1/2} \e^{t}} \setminus  D_{R^{1/2}
\e^{-t}}
 } d \bar{f} \wedge d f  =
        -   i   \,  \int_{\partial \left( D_{R^{1/2} \e^{t}} \setminus  D_{R^{1/2}
\e^{-t}} \right) } \bar{f} \, df  \, .
\end{equation}
 Differentiating \eqr{e:Et}, the chain rule and
$|\nabla f| \leq 1/ |z|$ give
\begin{equation}        \label{e:qds}
        E'(t) = R^{1/2} \e^{-t} \, \int_{\partial   D_{R^{1/2} \e^{-t}}} |\nabla f|^2  +
    R^{1/2} \e^{t} \, \int_{\partial   D_{R^{1/2} \e^{t}}} |\nabla f|^2
    \leq 4 \, \pi   \, .
\end{equation}
The Cauchy--Schwarz and Wirtinger inequalities (recall that
$\int_{\partial D_r} f = 0$)    give
\begin{align}   \label{e:b2}
        \left|  \int_{\partial D_r} \bar{f} \, df \right|
          &=  \left|  \int_{\partial D_r} \bar{f} \, f_{\theta} \, d\theta \right|
          \leq \left(  \int_{\partial D_r} |\bar{f}|^2 \, d\theta
          \right)^{1/2}
  \left(  \int_{\partial D_r} |f_{\theta}|^2 \, d\theta
          \right)^{1/2} \notag \\
          &\leq \int_{\partial D_r} |f_{\theta}|^2 \, d\theta =
          r/2  \, \int_{\partial D_r} |\nabla f |^2  \, ,
\end{align}
where the last equality used  the Cauchy--Riemann equations to
relate $f_{\theta}$ and $f_r$.  By \eqr{e:qs}--\eqr{e:b2}, we get
$2 \, E(t) \leq E'(t)$ and $E(1/2 \, \log R) \leq 2 \pi$.
Integrating this differential inequality yields \eqr{e:supeB}.
\end{proof}

Combining Lemma \ref{l:supeB} (with $t=\log 2$) and the mean value
inequality, we get
\begin{equation} \label{e:supeB2}
       \max_{\partial D_{ R^{1/2} } } |\nabla f|^2 \leq
       \left( \pi \, R / 4 \right)^{-1} \, \int_{D_{2 R^{1/2}}
       \setminus D_{R^{1/2} /2} } |\nabla f|^2 \leq
32 \, R^{-2} \, .
\end{equation}
If we now take $R\to \infty$ as in Bers' theorem, then we see that
\begin{equation}    \label{e:rtib}
    |\nabla f(z)| \leq C \, |z|^{-2} \, .
\end{equation}
The bound \eqr{e:rtib} is integrable radially (as $|z| \to
\infty$), so we see that $f$ has an asymptotic limit.

Finally, note that   $f(z) = 1/z$ shows that \eqr{e:supeB2} is
sharp (up to the constant).

\part{The role of multi--valued graphs in minimal surfaces} \label{p:2}

There are two local models for embedded minimal disks (by an {\it
  embedded disk} we mean a smooth injective map from the closed unit ball in
$\RR^2$ into $\RR^3$).  One model is the plane (or, more
generally, a minimal graph) and the other is a piece of a
helicoid.

The second model comes from the helicoid which was discovered by
Meusnier in 1776.  Meusnier had been a student of Monge.  He also
discovered that the surface now known as the catenoid is minimal
in the sense of Lagrange, and he was the first to characterize a
minimal surface as a surface with vanishing mean curvature. Unlike
the helicoid, the catenoid is not topologically a plane but rather
a cylinder.

\begin{figure}[htbp]
    \setlength{\captionindent}{20pt}
    \begin{minipage}[t]{0.5\textwidth}
    \centering\includegraphics{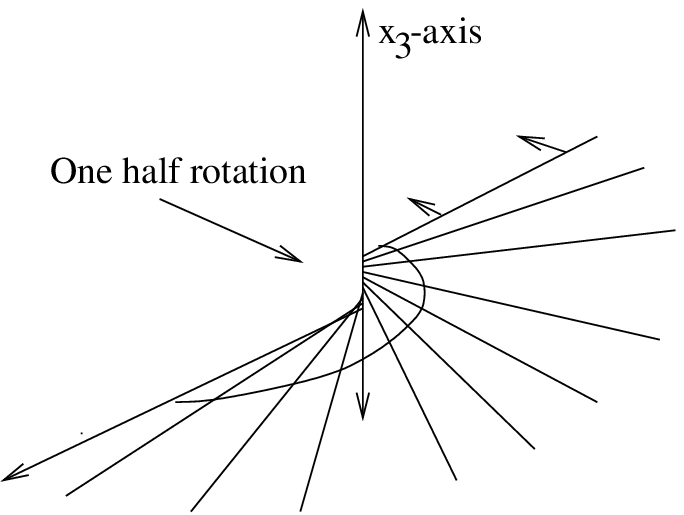}
    \caption{Multi--valued graphs.  The helicoid is obtained
    by gluing together two $\infty$--valued graphs along a line.}\label{f:f1}
    \end{minipage}%
    \begin{minipage}[t]{0.5\textwidth}
    \centering\input{pl8a.pstex_t}
    \caption{The separation $w$ grows/decays in $\rho$ at most sublinearly
for a multi--valued minimal graph; see \eqr{e:slgc}.}\label{f:f2}
    \end{minipage}
\end{figure}

The helicoid is a ``double spiral staircase'' (see \cite{CM17}):

\vskip2mm \noindent {\bf{Example 2}}: (Helicoid; see fig.
\ref{f:f1}).  The helicoid is the minimal surface in $\RR^3$ given
by the parametrization
\begin{equation}  \label{e:helicoid}
(s\cos t,s\sin t,t)\, ,\,\,\,\,\,\text{ where }s,\,t\in \RR\, .
\end{equation}

\section{Basic properties of multi--valued graphs}      \label{s:11}

We will need the notion of a multi--valued graph, each staircase
will be a multi--valued graph. Intuitively, an (embedded)
multi--valued graph is a surface such that over each point of the
annulus, the surface consists of $N$ graphs. To make this notion
precise, let $D_r$ be the disk in the plane centered at the origin
and of radius $r$ and let $\cP$ be the universal cover of the
punctured plane $\CC\setminus \{0\}$ with global polar coordinates
$(\rho, \theta)$ so $\rho>0$ and $\theta\in \RR$.  An {\it
$N$--valued graph} on the annulus $D_s\setminus D_r$ is a single
valued graph of a function $u$ over
\begin{equation}
    \{(\rho,\theta)\,|\,r< \rho\leq s\, ,\, |\theta|\leq
    N\,\pi\} \, .
\end{equation}
 For working purposes, we generally
think of the intuitive picture of a multi--sheeted surface in
$\RR^3$, and we identify the single--valued graph over the
universal cover with its multi--valued image in $\RR^3$.

The multi--valued graphs that we will consider will all be
embedded, which corresponds to a nonvanishing separation between
the sheets (or the floors).  Here the {\it separation} is the
function (see fig. \ref{f:f2})
\begin{equation}
w(\rho,\theta)=u(\rho,\theta+2\pi)-u(\rho,\theta)\, .
\end{equation}
If $\Sigma$ is the helicoid, then $\Sigma\setminus
\{x_3-\text{axis}\}=\Sigma_1\cup \Sigma_2$,
 where $\Sigma_1$, $\Sigma_2$ are $\infty$--valued graphs on
$\CC\setminus \{0\}$. $\Sigma_1$ is the graph of the function
$u_1(\rho,\theta)=\theta$ and $\Sigma_2$ is the graph of the
function $u_2(\rho,\theta)=\theta+\pi$.  ($\Sigma_1$ is the subset
where $s>0$ in \eqr{e:helicoid} and $\Sigma_2$ the subset where
$s<0$.)  In either case the separation $w=2\,\pi$.  A {\it
multi--valued minimal graph} is a multi--valued graph of a
function $u$ satisfying the minimal surface equation.

Note that for an embedded multi--valued graph, the sign of $w$
determines whether the multi--valued graph spirals in a
left--handed or right--handed manner, in other words, whether
upwards motion corresponds to turning in a clockwise direction or
in a counterclockwise direction.

\subsection{The sublinear growth of the separation}

As we have seen, the separation is constant for the multi--valued
graphs coming from each half of the helicoid.  This can be viewed
as a type of Liouville Theorem reflecting the conformal properties
of an infinite--valued graph.  In Proposition II.2.12 of
\cite{CM3}, we proved a corresponding gradient estimate:

\begin{proposition} \cite{CM3}  \label{p:slg}
 Given $\alpha > 0$, there exists $N_g$ so if
 $u$ satisfies the minimal surface equation
on $\{ \e^{-N_g} \, R  \leq \rho \leq \e^{N_g} \, R \, , \, - N_g
\leq \theta \leq 2 \pi + N_g \}$, $|\nabla u| \leq 1$, and
separation $w \ne 0$, then
\begin{equation}    \label{e:slg}
     |\Hess_u| (R,0) +   |\nabla \log |w ||(R,0)  \leq \alpha  / R  \, .
\end{equation}
\end{proposition}

One important consequence of \eqr{e:slg} is that, for $\alpha <
1$, the separation grows sublinearly:

\begin{corollary}     \label{c:slg}
 Given $\alpha > 0$, there exists $N_g$ so if
 $u$ satisfies the minimal surface equation
on $\{ \e^{-N_g} \, r_1  \leq \rho \leq \e^{N_g} \, r_2 \, , \, -
N_g \leq \theta \leq 2 \pi + N_g \}$, $|\nabla u| \leq 1$, and
separation $w \ne 0$, then
\begin{equation}    \label{e:slgc}
    |w | (r_2,0)  \leq  |w| (r_1 , 0) \left( \frac{r_2}{r_1} \right)^{\alpha}   \, .
\end{equation}
\end{corollary}

\begin{proof}
 Integrate $|\nabla \log |w ||  \leq \alpha / \rho$ along the ray
$\theta = 0$ to get
\begin{equation}    \label{e:slg2}
   \log \frac{ w (r_2,0) }{w(r_1,0)} \leq \alpha \int_{r_1}^{r_2}
   \rho^{-1} \, d \rho = \log \left( \frac{r_2}{r_1}
   \right)^{\alpha} \, .
\end{equation}
\end{proof}

Since $u (\cdot , \cdot)$ and its $2\pi$--rotation $u(\cdot ,
\cdot + 2\pi)$ are both solutions of the minimal surface equation,
the difference $w$ is almost a solution of the linearized equation
(which is the Jacobi equation in this case). Since the graphs have
bounded gradient, this equation is not too far from the Laplace
equation.
 To give some indication of why \eqr{e:slg} holds,
we will give an elementary proof when $u$ and $w$ are harmonic.

\begin{proof}
(of Proposition \ref{p:slg} when $u$  is harmonic.) After
rescaling, we can assume that $R=1$.  By making the conformal
change of coordinates $(\rho,\theta)\to (\log \rho,\theta)$ we get
a positive harmonic function
\begin{equation}
    \tilde w(x,y)=w(\e^x,y)
\end{equation}
 defined on the square $[-N_g , N_g] \times [-N_g , N_g]$.
Since the chain rule gives
\begin{equation}
    \nabla \log w (1,0) =
        \nabla \log \tilde w(0,0) \, ,
\end{equation}
 applying the Euclidean gradient estimate to $\tilde w$ yields
\begin{equation}    \label{e:ege}
    |\nabla \log  w (1,0)|  = |\nabla \log \tilde w (0,0)| \leq C/N_g \, .
\end{equation}
This gives the sublinear gradient estimate for $w$ if  $N_g$ is
sufficiently large.  The bound on $\Hess_u$ follows similarly.
\end{proof}

Proposition \ref{p:slg} allows us to assume (after rotating so
$\nabla u(1,0) = 0$) that
\begin{equation} \label{e:wantit}
    |\nabla u| + \rho \, |\Hess_u| +  4 \, \rho \, |\nabla  w |/|w| +
    \rho^2 \, |\Hess_w  |/|w|
      \leq \epsilon < 1/(2\pi) \, .
\end{equation}
The bound on $|\Hess_w  |$ follows from the other bounds and
standard elliptic theory.

\subsection{Curvature decay}

In corollary $1.14$ of \cite{CM8}, we proved faster than quadratic
curvature decay for double--valued minimal graphs whose separation
grows sublinearly:

\begin{proposition}     \label{c:decay}
There exists $C$ so if $u$ satisfies the minimal surface equation
and \eqr{e:wantit} on $\{ 1  \leq \rho \leq   R \, , \, - \pi \leq
\theta \leq 3 \pi  \}$, then on $\{ 1  \leq \rho \leq   R^{1/2} \,
, \, 0 \leq \theta \leq 2\pi  \}$
\begin{equation}   \label{e:cdecay1}
          \rho \, |\Hess_u| \leq C \, \epsilon \, \rho^{-5/12} \,  ,
\end{equation}
and therefore, after possibly rotating $\RR^3$ (and replacing
$u$), we get
\begin{equation}   \label{e:cdecay2}
         |\nabla u| \leq C \, \epsilon \, \rho^{-5/12} \,  .
\end{equation}
\end{proposition}

Of course, Proposition \ref{c:decay} is a generalization of Lemma
\ref{l:supeB}  which proved curvature decay for annuli.  As in the
annuli case, the proof uses the quasi--conformality of the Gauss
map to deduce a differential inequality.  However, the ``slit''
(i.e., where the double--valued graph does not close up)
contributes new terms which are controlled using the estimate for
the separation.  Notice that the second conclusion \eqr{e:cdecay2}
of the proposition proves a generalization of Bers theorem:
Embedded multi--valued minimal graphs have an asymptotic tangent
plane.

\section{Sharp estimates on the separation for
multi--valued graphs}       \label{s:12}

We will describe in this section two sharp estimates on the
separation of a multi--valued minimal graph; see Propositions
\ref{p:cm7} and \ref{p:upandlo} below. These estimates will,
unlike the earlier estimates \eqr{e:slgc} and \eqr{e:cdecay1},
require a rapidly growing number of sheets (growing in $\rho$).

Suppose for a moment that we are looking at an embedded surface
which is the $\infty$--valued graph of a harmonic function so that
in particular the separation $w$ is a harmonic function which
after reflection we may assume is positive. By making the
conformal change of coordinates $(\rho,\theta)\to (\log
\rho,\theta)$ we get a positive harmonic
\begin{equation}
    \tilde w(x,y)=w(\e^x,y)
\end{equation}
 defined on the half--plane $\{x\geq 0,\, y\in \RR\}$. By the
mean value inequality and the Harnack inequality (since $\tilde w$
is positive)
\begin{equation}\label{e:newone}
\tilde w (y,0)=\frac{1}{2\pi\,y}\int_{\partial D_y(y,0)}\tilde w
\geq \frac{1}{2\pi\,y}\int_{\partial D_y(y,0)\cap D_1}\tilde w
\geq \frac{C}{2\pi\,(y+1)}\, \tilde w(0,0)\, .
\end{equation}
Similarly, by an inversion formula one may show that
\begin{equation}  \label{e:two}
\tilde w (y,0)\leq \frac{C}{2\pi\,(y+1)}\, \tilde w(0,0)\, .
\end{equation}
For the original function $w$, \eqr{e:newone} combined with
\eqr{e:two} gives  for some constant $C$
independent of $w$
\begin{equation}
\frac{1}{C\,\log \rho} \leq
    \frac{w(\rho,0)}{w(1,0)}
\leq C\,\log \rho \, .
\end{equation}

\vskip2mm In the case of embedded multi--valued minimal graphs we
get similarly:

\begin{proposition} \label{p:cm7}
\cite{CM8}. Let $\Sigma$ be an embedded multi--valued minimal
graph of a function $u$ and with a rapidly growing number of
sheets, then for the separation $w$ we have for some constant $C$
\begin{equation}
\frac{1}{C\,\log \rho} \leq
    \frac{w(\rho,0)}{w(1,0)}
\leq C\,\log \rho \, .
\end{equation}
\end{proposition}

Suppose again for a moment that $u$ and hence $w$ is harmonic.
Similarly to \eqr{e:newone} we get that
\begin{equation}   \label{e:one}
\tilde w (0,y)\leq C\,2\,\pi\,(y+1)\, \tilde w(y,y)\, .
\end{equation}
By the Harnack inequality $\tilde w (y,y)\leq C\,\tilde w(y,0)$,
combining this with \eqr{e:one} and \eqr{e:two} we get
\begin{equation}  \label{e:desired}
\tilde w(0,y)\leq C\,(y^2+1)\,\tilde w(0,0)\, .
\end{equation}
For the original function $w$, this gives
\begin{equation}
\frac{1}{C\,(\theta^2+1)} \leq
    \frac{w(\rho,\theta)}{w(\rho,0)}
\leq C\,(\theta^2+1)\, .
\end{equation}

Again in the case of embedded multi--valued minimal graphs we get
similarly:

\begin{proposition} \label{p:upandlo}
Let $\Sigma$ be an embedded multi--valued minimal graph of a
function $u$  and with a rapidly growing number of sheets, then
for the separation $w$ we have for some constant $C$
\begin{equation}
\frac{1}{C\,(\theta^2+1)} \leq
    \frac{w(\rho,\theta)}{w(\rho,0)}
\leq C\,(\theta^2+1)\, .
\end{equation}
\end{proposition}

The lower bound in \eqr{e:one} for the decay of the separation is
sharp. It is achieved for the $\infty$--valued graph of the
harmonic function (graphs of multi--valued harmonic functions are
good models for multi--valued minimal graphs)
\begin{equation}  \label{e:exofu}
u(\rho,\theta)=\arctan  \frac{\theta}{\log \rho} \, .
\end{equation}
Note that the graph of $u$ is embedded and lies in a slab in
$\RR^3$, i.e., $|u|\leq \pi/2$, and hence in particular is not
proper.  On the top it spirals into the plane $\{x_3=\pi/2\}$ and
on the bottom into $\{x_3=-\pi/2\}$, yet it never reaches either
of these planes.

\begin{question}
It would be interesting to construct an infinite--valued exterior
solution of the minimal graph equation with the same properties as
$\arctan \frac{\theta}{\log \rho}$; i.e., one which spirals
infinitely in a slab (see \cite{CM18} for a local example).
\end{question}

\section{Double--valued minimal graphs}     \label{s:13}

We will now describe how to bound the oscillation of the Gauss map
of a double--valued minimal graph.  This bound was proven in
\cite{CM10}.  Rather than give   the precise statement here, we
will instead illustrate a few of the key ideas by considering the
analogous situation for a holomorphic function $f$.  Here we are
of course thinking of $f$ as being the stereographic projection of
the Gauss map. The sublinear growth of the separation, i.e.,
\eqr{e:wantit} and its integrated form, correspond to
\begin{align}    \label{e:wantit2a}
    |f| + \rho \, |\nabla f|
      &\leq \epsilon < 1/(2\pi) \, , \\
|f(\rho , 2\pi) - f(\rho , 0)|
      &\leq  2 \pi \,  \epsilon \, \left( \frac{\delta}{\rho} \right)^{1 -\epsilon } \, . \label{e:wantit2b}
\end{align}

The bound on the oscillation of $f$ (Lemma \ref{l:supeq} below)
now follows by modifying the argument for annuli (Lemma
\ref{l:supe}). Since $f$ does not match up at $\theta =0$ and
$2\pi$, we get an   additional term which is estimated using the
sublinear growth \eqr{e:wantit2b}.

 \begin{lemma} \label{l:supeq} If
$f:\{ \delta \leq \rho \leq R \, , \, 0 \leq \theta \leq 2\pi \}
\to \CC$
 is holomorphic and satisfies \eqr{e:wantit2a} and \eqr{e:wantit2b}, then
\begin{equation} \label{e:supeq}
\min_{c \in \CC} \, \max  |f  - c| \leq
\frac{\epsilon}{1-\epsilon} + 2\pi \, \epsilon \, .
\end{equation}
\end{lemma}

\begin{proof}
(Following Lemma \ref{l:supe}.)
 For $\delta \leq s \leq R$, define the circular average
\begin{equation} \label{e:iofsq}
I(s) =  (2 \, \pi)^{-1} \, \int_0^{2 \, \pi} f(s , \theta) \, d
\theta \, .
\end{equation}
Note that integrating \eqr{e:wantit2a} gives
\begin{equation} \label{e:ivsf}
\left| I(s) - f (s,\theta) \right| \leq 2\pi \, \epsilon \, .
\end{equation}
 Differentiating \eqr{e:iofsq} and using
$\frac{\partial f}{\partial \rho} = -i \, \rho^{-1} \,
\frac{\partial f}{\partial \theta}$ since $f$ is holomorphic, we
have
\begin{align} \label{e:iofspq}
2 \, \pi \, s \, I'(s) &= s \, \int_{0}^{2\pi} \frac{\partial
f(s,\theta)}{\partial \rho} \, d \theta = -i \,  \int_0^{2 \, \pi}
\frac{\partial f(s,\theta)}{\partial \theta}  \, d \theta
\notag \\
&= -i \,  [ f(s  , 2 \, \pi) - f(s , 0)]   \, .
\end{align}
Using the bound \eqr{e:wantit2b} along the slit, \eqr{e:iofspq}
gives
\begin{equation}
     s \,| I'(s)| \leq \epsilon \, \delta^{1 -\epsilon } \,
    s^{\epsilon - 1} \, .
\end{equation}
 In particular, integrating this gives for $\delta \leq \rho \leq
 R$ that
 \begin{equation} \label{e:wb3}
  |I(\rho) - I(\delta)| \leq \epsilon \, \delta^{1 -\epsilon } \, \int_{\delta}^{\rho}
 s^{\epsilon - 2} \, ds \leq
\epsilon / (1-\epsilon) \, .
\end{equation}
 The bound \eqr{e:supeq}  follows from  \eqr{e:ivsf} and \eqr{e:wb3}.
\end{proof}

Modifying this argument to apply to double--valued minimal graphs
introduces new  difficulties.  In that case, one works directly on
the graph (where the Gauss map is holomorphic) and uses averages
over ``geodesic sectors'' rather than circles in the plane.   One
  difficulty is a new term in the analog of \eqr{e:iofspq} which
results from differentiating the measures of the level sets.

\section{Approximation by standard pieces}  \label{s:14}

In this section we show that any embedded multi--valued minimal
graph has a sub--graph which is close to the sum of a piece of a
catenoid  and a piece of a helicoid.  This generalizes a similar
representation for minimal graphs over an annulus given in
proposition 1.5 in \cite{CM9}; of course, there was no helicoid
term in that case. These results are new and have not appeared in
the literature elsewhere.

Recall that half of a catenoid, i.e.,
\begin{equation}
    \{(x_1/s)^2+(x_2/s)^2=\cosh^2 (x_3/s)\, , \,\pm x_3>0\} \, ,
\end{equation}
 is a minimal graph of $u(z)=\pm s \cosh^{-1} (|z|/s)$ over
$\CC\setminus D_s$.  Note that $s \cosh^{-1} (|z|/s)$ is
asymptotic to $s \log [2\,|z|/s]$. Recall also that half of the
helicoid is the multi--valued graph of the function $u$ given in
polar coordinates by $u(\rho,\theta)= \theta$. Our approximation
result (see Corollary \ref{c:approxbyst} below) is therefore that
any embedded multi--valued minimal graph has a sub--graph which is
close to the graph of a multi--valued function $v$ given by
($a,\,b,\,c\in\RR$ are constants)
\begin{equation}
v (\rho, \theta ) =
  a + b \, \log (\rho/ r) + c\, \theta/(2\pi) \, .
\end{equation}

As in the two previous sections,   $u$ will be a multi-valued
function and $w$ will be its separation.
 For convenience, we will write $S_{r_1,r_2}^{\theta_1 ,
 \theta_2}$ to denote the ``rectangle''
 \begin{equation}
    \{ (\rho , \theta) \, | \, r_1 \leq \rho \leq r_2 , \,
 \theta_1 \leq \theta \leq \theta_2 \} \, .
 \end{equation}

\vskip2mm
 We begin with a representation formula for the gradient
of an ``almost harmonic'' function. Recall that if $\Delta u = 0$
on an annulus, then the function $f= u_x - i \, u_y$ is
holomorphic. In particular,   $f$ has a  Laurent expansion which
can be recovered
 using the Cauchy integral formula. The next lemma uses a
variation on this for multi--valued ``almost harmonic'' functions:

\begin{lemma}     \label{l:decomp}
Given a function $u$ on $S_{1,{\sqrt{R}}}^{-\pi,3\pi}$, set $f =
u_x - i \, u_y$. If
\begin{align}
    |f(\zeta)| &\leq C \, |\zeta
    |^{-5/12} \, , \\
 |\Delta u(\zeta)| &\leq C \, |\zeta |^{-9/4} \, , \\
 \rho \, |\nabla  w |/|w|
     & \leq \epsilon < 1/(2\pi) \, ,
\end{align}
then
 for $r_1 \geq 1$ and $\zz\in
S_{2r_1,\sqrt{R}/2}^{0,2\pi}$
\begin{equation}    \label{e:dd1}
f (\zz) = (b+i\,c/(2\pi))\,\zz^{-1} + g(\zz)
\end{equation}
where $b,\,c \in \RR$ and
\begin{equation}    \label{e:dd2}
 |g (\zeta)| \leq  C_1 \,  R^{-5/24}
+  C_1 \,  r_1^{-1/4} \, |\zeta|^{-1}
+ C_1 \, \epsilon \, r_1^{-1} \, |w(r_1 ,0)|   \, .
\end{equation}
\end{lemma}

\begin{proof}
We will first use the Cauchy integral formula (this is just
Stokes' theorem applied to the one--form $f(z)/ (z-\zz) \, dz$) on
the domain $S_{1,\sqrt{R}}^{0,2\pi}$ to get a representation
formula on $S_{2\,r_1,\sqrt{R}/2}^{\pi/2,3\pi/2}$. Assume that
$\zz\in S_{2\,r_1,\sqrt{R}/2}^{\pi/2,3\pi/2}$. The Cauchy integral
formula gives
\begin{align}   \label{e:cifs}
        2 \, \pi \, i \, f(\zz) &
=  \int_{S_{{\sqrt{R}},{\sqrt{R}}}^{0,2\pi}}  \frac{ f(z) }{ z-
\zz} \, dz
    - \int_{S_{1,1}^{0,2\pi}} \frac{ f(z) }{ z-
\zz} \, dz  - i \, \int_{S_{1,{\sqrt{R}}}^{0,2\pi} }
        \frac{\Delta u}{z-\zz} \notag  \\ &\quad+
\left( \int_{{\sqrt{R}}}^{1}  \frac{ f(\rho  , 2 \, \pi)}{ \rho -
\zz} \, d\rho  + \int^{{\sqrt{R}}}_{1} \frac{ f(\rho  , 0) }{ \rho
- \zz} \, d\rho  \right) \, .
\end{align}
The first two terms correspond to the usual formula for annuli,
the third term vanishes when $f$ is holomorphic, and the last term
arises since $f$ is multi--valued. (The first three are almost
identical to the corresponding ones in lemma 1.7 in \cite{CM9}.)
 To prove the lemma, we will show that the first term is small and
 the other three are small after subtracting a multiple of
 $1/\zz$.

 First, $|f(z)| \leq C \, |z|^{-5/12}$ and
  $|\zz| \leq \sqrt{R}/2$ give that
\begin{equation}   \label{e:ft}
        \left| \int_{S_{{\sqrt{R}},{\sqrt{R}}}^{0,2\pi}} \frac{ f(z) }{ z- \zz}  \, dz \right|
    \leq 4 \, \pi \, \sup_{S_{{\sqrt{R}},{\sqrt{R}}}^{0,2\pi}} |f|
\leq 4 \pi \, C \, R^{-5/24} \, .
\end{equation}

 For the remaining terms, we will use the identity
\begin{equation}    \label{e:theid}
    \frac{1}{z-\zz} = \frac{-1}{\zz}  + \frac{z}{\zz(z -\zz)} \, .
\end{equation}

Second, using \eqr{e:theid}, $2 \leq |\zz|$, and $|f|\leq C$,
\begin{equation}        \label{e:st1}
         \left| \int_{S_{1,1}^{0,2\pi}} \frac{f(z)}{z- \zz} \, dz
+ \frac{1}{\zz}
  \,  \int_{S_{1,1}^{0,2\pi}} f(z) \, dz  \right| \leq
\frac{2}{ |\zz|^{2}} \,  \int_{S_{1,1}^{0,2\pi}} |z \, f(z)|
\leq   \frac{4 \pi \, C}{ |\zz|^{2}}
 \,  .
\end{equation}

 To bound the third integral, we separate out  the disk of radius $|\zz|/2$ about
 $\zz$ and divide the remainder into two regions using the circle
 of radius $r_1$.  Using \eqr{e:theid}, $2\, r_1 \leq |\zz|$, and $|\Delta
u(z)| \leq C \, |z|^{-9/4}$, we get
\begin{align}   \label{e:nonh}
         \left| \int_{S_{1,{\sqrt{R}}}^{0,2\pi} }
        \frac{\Delta u}{z-\zz} + \frac{1}{\zz}
  \,  \int_{S_{1,r_1}^{0,2 \pi} } \Delta u
 \right|
&\leq 4\, \pi \, C \, |\zz|^{-2} \, \int_{1}^{r_1}
        r^{-1/4} \, dr
        +  2^{9/4} \, C\,  \int_{D_{\frac{|\zz|}{2}}(\zz)}
\frac{|\zz|^{-9/4}}{|z-\zz|} \notag \\
    &\quad \quad
+ 4 \, \pi C \, |\zz|^{-1} \int_{r_1}^{\sqrt{R}} r^{-5/4} \, dr
\notag \\
          &\leq
    12 \, \pi  \, C \,  |\zz|^{-5/4} +
    16 \, \pi  \, C \,  r_1^{-1/4} \, |\zz|^{-1} \, .
\end{align}

Finally, this leaves the slit (i.e., the terms which arise because
$u$ is not well-defined over the annulus), where we divide the
integral into three parts
\begin{align}   \label{e:tt}
 &\int^{{\sqrt{R}}}_{1}  \frac{ f( \rho , 2 \pi) - f(\rho  ,0 ) }{  \rho  - \zz} \, d\rho   =
\int_{r_1}^{{\sqrt{R}}} \frac{ f(\rho  , 2 \pi) - f(\rho  ,0 ) }{ \rho  - \zz} \, d\rho   \\
&\quad \quad \quad - \zz^{-1} \, \int^{r_1}_{1} (f(\rho  , 2 \pi)
- f(\rho  ,0 ) ) \, d\rho  + \int^{r_1}_{1} \frac{ \rho (f(\rho  ,
2 \pi) - f(\rho  ,0 ) ) }{ \zz(\rho -\zz)} \, d\rho \, . \notag
\end{align}
Using  $\rho \, |\nabla  w |/|w|
      \leq \epsilon$, we
  get for $\rho \geq r_1$ that
\begin{equation}
 \label{e:sf0}
        |f(\rho , 2 \pi) - f(\rho ,0)| = |\nabla w (\rho,0)|
\leq \epsilon \, \frac{|w|(\rho,0)}{\rho} \leq
         \epsilon \,  \left( \frac{\rho}{r_1}\right)^{\epsilon}
\, \frac{ |w|(r_1,0)}{\rho}  \, .
\end{equation}
Note that the real part of $\zz$ is negative since $\zz \in
S_{2r_1,\sqrt{R}/2}^{\pi/2 , 3\pi/2}$ and hence $|\rho - \zz| \geq
\rho$.  To bound the first term in \eqr{e:tt}, we use this and
\eqr{e:sf0} to get
\begin{equation}
\left| \int_{r_1}^{{\sqrt{R}}} \frac{f(\rho  , 2 \pi) - f(\rho  ,0 ) }{\rho  - \zz}
     \, d\rho
\right| \leq  \frac{\epsilon}{r_1^{\epsilon}} \,
|w|(r_1 ,0) \,
\int^{{\sqrt{R}}}_{r_1} \rho^{\epsilon - 2}
     \, d\rho   \leq  2 \epsilon \, \frac{|w|(r_1,0)}{ r_1} \, .
\end{equation}
Similarly, to bound the third,  use  $\rho \, |\nabla  w |/|w|
      \leq \epsilon$
 to get
\begin{align}
\left| \int^{r_1}_{1}
\frac{\rho \, (f(\rho  , 2 \pi) - f(\rho  ,0 ) )}
{\zz(\rho  -\zz)}
\, d\rho  \right| &\leq \epsilon \, r_1^{-2} \, \int_{1}^{r_1} |w|(\rho ,0)
\, d\rho   \\ &\leq \epsilon \, r_1^{-2} \, \int_{1}^{r_1} |w|(r_1,0) \,
(r_1/\rho )^{\epsilon} \, d\rho  \leq 2 \epsilon \, \frac{|w|(r_1,0)}{r_1} \, .
\notag
\end{align}

Putting all of this together, we now get the desired
representation formula   \eqr{e:dd1} and remainder estimate
\eqr{e:dd2} for $\zz \in S_{2r_1,\sqrt{R}/2}^{\pi/2 , 3\pi/2}$.
Namely, we set   $g=g_1+g_2+g_3$ and
\begin{align}
2 \pi  i \, (b+i\, \frac{c}{2\pi}) &=  \int_{S_{1,1}^{0,2\pi}} f(z) \, dz +
\int^{r_1}_{1} (f(\rho  , 2 \pi) - f(\rho  ,0 ) ) \, d\rho
+i   \,  \int_{S_{1,r_1}^{0,2\pi} } \Delta u \, ,  \notag \\
2 \pi i
\, g_1(\zz) &= \int_{S_{{\sqrt{R}},{\sqrt{R}}}^{0,2\pi}} \frac{f(z)}{z- \zz} \, dz \, ,
\\
2 \pi i
\, g_2(\zz) &=  -
\int_{r_1}^{{\sqrt{R}}} \frac{f(\rho , 2 \pi) - f(\rho  ,0 ) }{\rho  - \zz} \, d\rho  -
\int^{r_1}_{1} \frac{f(\rho  , 2 \pi) - f(\rho  ,0 )}{\zz(\rho  -\zz)} \, \rho  \,
d\rho  \, , \notag
\\
2 \pi  i \, g_3 (\zz) &= - \int_{S_{1,1}^{0,2\pi}} \frac{f(z)}{z- \zz} \,
dz - \frac{1}{\zz}
  \,  \int_{S_{1,1}^{0,2\pi}} f(z) \, dz
- i \, \int_{S_{1,{\sqrt{R}}}^{0,2\pi} }
        \frac{\Delta u}{z-\zz}
-  \frac{i}{\zz}
  \,  \int_{S_{1,r_1}^{0,2\pi} } \Delta u  \, .  \notag
\end{align}
We can repeat this  by integrating over
$S_{1,{\sqrt{R}}}^{-\pi/2, 3\pi/2}$ and $S_{1,{\sqrt{R}}}^{\pi/2, 5\pi/2}$ to
get similar representations on
$S_{2r_1,\sqrt{R}/2}^{0 , \pi}$ and $S_{2r_1,\sqrt{R}/2}^{\pi , 2\pi}$
(with different values of $b$ and $c$).  By continuity, the same
representation holds on all three regions (since they overlap), giving
the lemma.
\end{proof}

\begin{definition}  \label{d:stapiece}
Fix $\mu >0$.
A standard piece $\Sigma_{v}$ on the scale $r$ is a
graph of $v$ over $S_{r,\mu r}^{-\pi,3\pi}(y)$
given in polar coordinates $(\rho_p , \theta_p)$
centered at $p$  by
\begin{equation}
v (\rho_p, \theta_p ) =
  a + b  \, \log (\rho_p/ r) +   c \, \theta_p / (2\pi) \, .
\end{equation}
\end{definition}

The constant $a$ is the ``plane coefficient'', $b$ is the
``catenoid coefficient'' and $c$ is the ``helicoid coefficient''
(this gives the separation $w$). Note that $v$ is harmonic.

\begin{corollary}  \label{c:approxbyst}
Let $u$ be a solution of the minimal graph equation on
$S_{1/2,R}^{-3\pi,5\pi}$ with $w < 0$ satisfying \eqr{e:wantit}
and $|u| \leq \epsilon \, \rho$. There is a rotation of $\RR^3$ so
given $r_1 > 2$, we get a standard piece
\begin{equation}
    v =a + b
        \, \log \rho + c \, \theta/(2\pi) \, ,
\end{equation}
  for $a,b,c
\in \RR$ with
\begin{equation}  \label{e:capprox}
\sup_{S_{r_1,\mu r_1}^{0,2\pi}} |u (\rho ,\theta)-v (\rho,\theta)|\leq C_2 \, \mu \,
\epsilon \, |w|(r_1,0)  +
C_2 \, \mu \, r_1 \, R^{-5/24}
+  C_2 \, \mu \,  r_1^{-1/4}
\, .
\end{equation}
\end{corollary}

\begin{proof}
Corollary 1.14 in \cite{CM8}  give a  rotation of $\RR^3$ so that
 for $z \in S_{1,\sqrt{R}}^{-2\pi,4\pi}$
\begin{equation}    \label{e:qcm1}
    |\nabla u(z)| + |z| \, |\Hess_u (z)|  \leq
     C \, \epsilon \, |z|^{-5/12}\,   .
\end{equation}
Since, by (1.6) of \cite{CM8}, $|\Delta \,u|\leq |\nabla
u|^2\,|\Hess_u|$, \eqr{e:qcm1} gives on
$S_{1,\sqrt{R}}^{-2\pi,4\pi}$ that
\begin{equation}   \label{e:cdecay1a}
| \Delta u (z)|  \leq C \, \epsilon^3 \, |z|^{-9/4} \,  .
\end{equation}
Note that $(\log \rho)_x - i \, (\log \rho)_y = 1/ (x+ i y)$ and
$\theta_x - i \, \theta_y = i / (x+ i y)$.  Using this and Lemma
\ref{l:decomp},
 we get that
on $S_{r_1,\mu r_1}^{0,2\pi}$
\begin{equation}        \label{e:diffg}
        |\nabla ( u - b \, \log \rho - c \, \theta/(2\pi))|
        \leq C_1 \, \left( R^{-5/24}
+  2 \,  r_1^{-1/4} \, \rho^{-1}
+ 2 \, \epsilon \, \frac{|w(r_1 ,0)|}{r_1} \right) \, .
\end{equation}
Integrating \eqr{e:diffg} gives \eqr{e:capprox}.
\end{proof}

\begin{question}
Corollary \ref{c:approxbyst} shows that embedded multi-valued
minimal graphs are closely approximated by helicoids (plus a
catenoid term) on each scale.  However, a priori, the helicoid
coefficient can change from scale to scale (the point of
\cite{CM10} is that the axis of the  helicoid   does not change so
that all  are vertical).  It would be interesting to estimate how
quickly this helicoid coefficient can change and construct
examples demonstrating this.
\end{question}

\part{Regularity theory}    \label{p:3}

In this part, we survey some of the key  ideas in classical
regularity theory, recent developments on embedded minimal disks,
and some global results for minimal surfaces in $\RR^3$. Sections
\ref{s:16} and \ref{s:17} focus on Reifenberg type conditions,
where a set is assumed to be close to a plane at all points and at
all scales (``close'' is in the Hausdorff or Gromov--Hausdorff
sense and is defined in Section \ref{s:16}).  This condition
automatically gives H\"older regularity (and hence higher
regularity if the set is also a weak solution to a natural
equation).  Section \ref{s:18} surveys the role of monotonicity
and scaling in regularity theory, including $\epsilon$-regularity
theorems (such as Allard's theorem) and tangent cone analysis
(such as Almgren's refinement of Federer's dimension reducing).
Section \ref{s:19} briefly reviews recent results of the authors
for embedded minimal disks, developing a regularity theory in a
setting where the classical methods cannot be applied and in
particular where there is no monotonicity.  The estimates and
ideas discussed in Section \ref{s:19} have applications to the
global theory of minimal surfaces in $\RR^3$.  In Section
\ref{s:20}, we give a quick tour of some recent results in this
classical, but rapidly developing, area.

\section{Hausdorff and Gromov--Hausdorff distances}   \label{s:16}

Recall that the
Hausdorff distance between
two subsets $A$ and $B$ of a Euclidean space
is no greater than $\delta >0$ provided that each is contained in a
$\delta $-neighborhood of the other.

There is a natural generalization of this classical Hausdorff distance
between subsets of Euclidean space to a distance function on all metric
spaces.  This is the Gromov-Hausdorff distance and it gives a
good tool to study metric spaces.

Suppose that $(X,d_X)$ and $(Y,d_Y)$ are two compact
metric spaces.  We say that the
Gromov-Hausdorff distance between them is at most $\epsilon>0$ if there exist
maps $f:X\to Y$ and $g:Y\to X$ such that
\begin{equation}
\forall x_1,\, x_2\in X:\, \, |d_X(x_1,x_2)-d_Y(f(x_1),f(x_2))|<\epsilon\, ,
\end{equation}
\begin{equation}
\forall x\in X:\, \, d_X(x,g\circ f(x))<\epsilon\, ,
\end{equation}
and the two symmetric properties in $Y$ hold.
The Gromov-Hausdorff distance between $X$ and $Y$, denoted by $d_{GH}(X,Y)$,
is then the infimum of all such $\epsilon$.

Using this topology we can then say that a sequence of metric
space converges to another metric space.  For noncompact metric
spaces there is a more useful notion of convergence which is
essentially convergence on compact subsets.  Namely, if
$(X_i,x_i,d_{X_i})$ is a pointed sequence of metric spaces then we
say that $(X_i,x_i,d_{X_i})$ converges to some pointed metric
space $(X,x,d_X)$ in the pointed Gromov-Hausdorff topology if for
all $r_0>0$ fixed the compact metric spaces $B_{r_0}(x_i)$
converges to $B_{r_0}(x)$.

Gromov's compactness theorem is the statement that any pointed
sequence, $(M_i^n, m_i )$, of $n$-dimensional manifolds, with
\begin{equation}   \label{e:o1}
\Ric_{M_i^n} \geq (n-1)\,\Lambda\, ,
\end{equation}
has a subsequence, $(M_j^n, m_j )$, which converges
in the pointed Gromov-Hausdorff topology to some length space
$(M_{\infty},m_{\infty})$.

The proof of this compactness theorem relies only on
the volume comparison theorem.  In fact it only uses the
volume doubling which is implied by the volume comparison.

\section{Reifenberg type conditions}    \label{s:17}

One concept or idea that often plays a central role in regularity theory
is that of being
Reifenberg flat, a notion introduced by E. R. Reifenberg
\cite{Re}.  In the case of Reifenberg it was introduced to measure the
closeness of a subset of Euclidean space to being an affine space.  The
deviation was measured on all scales and a number was assigned which was
the maximal scale invariant Hausdorff distance to a affine plane, see
\cite{To}, \cite{Se}, \cite{DaTo}, \cite{DaKeTo}
for interesting results in this classical direction.

To explain this point about the importance of this condition further
we recall the following definition from \cite{ChC1} which was
inspired by the classical
work of Reifenberg.

Let $(Z,d_Z)$ be a complete metric space.  We will say that $Z$ satisfies the
$(\epsilon,\rho,n)$-$\mathcal{G_R}$ (or {\it Generalized Reifenberg})
condition at $z\in Z$ if for all
$0<\sigma< \rho$ and all $y\in B_{\rho-\sigma}(z)$,
\begin{equation}
d_{GH}(B_{\sigma}(y),B_{\sigma}(0))<\epsilon\, \sigma\, ,
\end{equation}
where $B_{r_0}(0)\subset \RR^n$.

\begin{theorem}
(Appendix of \cite{ChC2}.)
Given $\epsilon>0$ and $n\geq 2$, there exists $\delta>0$ such that if
$(Z,d_Z)$ is a complete metric space and $Z$ satisfies the
$(\delta,r_0,n)$-$\mathcal{G_R}$ condition
at $z\in Z$ then there exists a bi-H\"older
homeomorphism $\Phi: B_{\frac{r_0}{2}}(z)\to B_{\frac{r_0}{2}}(0)$ such that
for all $z_1,\, z_2\in Z$,
\begin{equation}
r_0^{-\epsilon}\, |\Phi (z_1)-\Phi (z_2)|^{1+\epsilon}
\leq d_{Z}(z_1,z_2)\leq
r_0^{\epsilon}\,|\Phi (z_1)-\Phi (z_2)|^{1-\epsilon} \, .
\end{equation}
\end{theorem}

Here is one example where such a Reifenberg type condition
naturally come up; see \cite{C2}, \cite{ChC1} for more on this.

\begin{theorem}
(\cite{C1}, \cite{ChC1}).
Given $\epsilon>0$ and $n\geq 2$, there exist
$\delta =\delta (\epsilon,n)>0$ and
$\rho=\rho (\epsilon,n)>0$, such
that if $M^n$ has $\Ric_{M^n}\geq - (n-1)$ and
$0<r_0\leq \rho$ with
\begin{equation}
d_{GH}(B_{2r_0}(x),B_{2r_0}(0))<\delta\, r_0\, ,
\end{equation}
or
\begin{equation}
\Vol (B_{2r_0}(x))\geq (1-\delta)\, \V_0^n(2r_0)\, ,
\end{equation}
then $M^n$ satisfies the $(\epsilon,r_0,n)$-$\mathcal{G_R}$ condition at
$x$.
\end{theorem}

\section{Monotonicity and regularity theory}        \label{s:18}

As we have already seen, the monotonicity of a ``scale--invariant
energy'' has played a key role in the regularity theory for
geometric variational problems.  In many cases, this monotonicity
is useful for establishing two key tools:
\begin{enumerate}
\item
An $\epsilon$--regularity theorem which guarantees that a weak
solution is actually smooth when the scale--invariant energy is
small.
\item
The existence of tangent cones which are dilation invariant.
\end{enumerate}
In this section, we briefly review the role of monotonicity in
regularity theory, emphasizing these two tools, and give some
examples.

To keep things concrete, it may be useful to mention some examples
of variational problems and their scale--invariant energies, some
of which we have already encountered.  The primary example is for
a minimal $k$--dimensional submanifold; in this case, the
appropriate  scale--invariant energy is the density
\begin{equation}                \label{e:thetadefvv}
        \Theta_{x_0}(s)
        =\frac{\Vol (B_s (x_0) \cap \Sigma)}{\Vol (B_s \subset \RR^k)
        }\, .
\end{equation}
There are many other examples, including
\begin{itemize}
\item
If $\Sigma \subset \RR^n$ is a ($2$--dimensional) minimal surface,
then another possible scale--invariant energy is $\int_{B_r}
|A|^2$.
\item
If $u:\RR^k \to M$ is a harmonic map, then the scale--invariant
energy is  $r^{2-k} \int_{B_r} |\nabla u|^2$, see \cite{ScUh}.
\item
If $M^n$ has non--negative Ricci curvature, then the
scale--invariant energy is the volume density $r^{-n} \, \Vol (B_r
(x_0) \subset M )$ -- here, the monotonicity goes the opposite
direction.
\item
If $A$ is a Yang--Mills connection on $M^k$ with curvature $F_A$,
then the scale--invariant energy is $r^{4-k} \, \int_{B_r}
|F_A|^2$, see \cite{Pr}.
\end{itemize}
Interestingly, there are also parabolic analogs for all of these
examples (cf. Huisken's monotonicity formula for the mean
curvature flow, i.e., Theorem \ref{t:mch} or M. Struwe's
monotonicity for the harmonic map heat flow, theorem $1.10$ in
\cite{St}).

  Using this monotonicity, we can define the density
$\Theta_{x_0}$ at the point $x_0$ to be the limit as $r \to 0$ of
the scale--invariant energy in $B_r (x_0)$. It follows easily from
the monotonicity that the density is semi--continuous.

\begin{remark}
These energies are all scale--invariant in the sense that if a
solution is rescaled by a factor $\lambda$, then the energy of the
new solution on a ball of radius $\lambda \, r$ is the same as the
original solution on a ball of radius $r$.  It is not difficult to
find scale--invariant quantities -- the difficulty lies in finding
one that is monotone and can be estimated.  For example, if $u:
\RR^k \to M$ is a harmonic map, then $\int |\nabla u|^{k}$ is
scale--invariant and monotone; however, it is not at all clear
that there should be an a priori bound on this when $k> 2$.
\end{remark}

\vskip2mm \subsection{$\epsilon$--regularity and the singular set}

This monotonicity is then important in proving an
$\epsilon$--regularity theorem. Recall that an
$\epsilon$--regularity theorem gives that a weak (or generalized)
solution is actually smooth at a point if the scale--invariant
energy is small enough there. The standard example is the Allard
regularity theorem:
\begin{theorem} \cite{Al}   \label{t:allard}
There exists $\delta (k , n) > 0$ such that if $\Sigma \subset
\RR^{n}$ is a $k$--rectifiable stationary varifold (with density
at least one a.e.), $x_0 \in \Sigma$, and
\begin{equation}
    \Theta_{x_0} = \lim_{r\to 0} \frac{ \Vol (B_r (x_0) \cap \Sigma)}
        {\Vol ( B_r \subset \RR^k )}
        < 1 + \delta \, ,
\end{equation}
 then $\Sigma$ is smooth in a neighborhood of $x_0$.
\end{theorem}

 Another well--known example is the $\epsilon$--regularity theorem
 of Schoen and K. Uhlenbeck for harmonic maps:

\begin{theorem} \cite{ScUh}   \label{t:scsu}
There exists $\epsilon (k ,N) > 0$ such that if $u: B_r \subset
\RR^k \to N$ is an energy minimizing map and
\begin{equation}
    r^{2-k} \int_{B_r} |\nabla u|^2 < \epsilon \, ,
\end{equation}
 then $u$ is smooth in a neighborhood of $0$ and $|\nabla u|(0)
 \leq C \, r^{-1}$.
\end{theorem}

See \cite{ChCTi} for $\epsilon$-regularity results for limits of
K\"ahler--Einstein metrics and their applications to regularity of
such limit spaces and see
 F.H. Lin, \cite{Ln}, for recent developments on the regularity of harmonic
 maps which are not minimizing.

\begin{remark}
The $\epsilon$--regularity theorem is closely related to  a priori
estimates when the energy is small (cf. Section \ref{s:9}). In one
direction, combining the regularity with a compactness argument
usually directly gives an a priori estimate.
\end{remark}

\vskip2mm
  The singular set $\cS$ is
defined to be the set where the scale--invariant energy is not
small. The first application of these $\epsilon$--regularity
theorems is some  control on the singular set.   For example, the
semi--continuity of the density immediately gives that $\cS$ is
closed.  In order to bound the size of the singular set (e.g., the
Hausdorff measure), we combine the $\epsilon$--regularity with
simple covering arguments:

\begin{lemma}     \label{l:std}
If $u:B_1 \subset \RR^k \to M$ is an energy minimizing  (harmonic)
map, then the $(k-2)$--dimensional Hausdorff measure of $\cS$ is
zero.
\end{lemma}

\begin{proof}
Given $\delta > 0$ and $x \in \cS$, then the
$\epsilon$--regularity theorem (Theorem \ref{t:scsu}) yields a
ball $B_{r_x}(x)$ so that
\begin{equation}
    0 < \epsilon \leq r_x^{2-k} \, \int_{ B_{r_x} (x) } |\nabla u|^2   \, .
\end{equation}
Using the $5$--times covering lemma, we can find a disjoint
collection of balls $B_{r_i}(x_i)$ so that
\begin{align} \label{e:bad}
    \epsilon &\leq r_i^{2-k} \,
    \int_{ B_{r_i}(x_i) } |\nabla u|^2  \, , \\
    \cS &\subset \cup_{i} B_{5r_i}(x_i) \, . \label{e:good}
\end{align}
Since the balls $B_{r_i}(x_i)$ are disjoint, \eqr{e:bad} implies
that
\begin{equation}    \label{e:ugly}
    \epsilon \, \sum_i r_i^{k-2} \leq \sum_i \int_{ B_{r_i}(x_i) } |\nabla
    u|^2 \leq \int_{\cup_i B_{r_i}(x_i) } |\nabla u|^2 \, .
\end{equation}
Combining this with \eqr{e:good}, we get a uniform bound for the
$(k-2)$--dimensional Hausdorff measure of $\cS$.  In particular,
$\cS$ has Lebesgue measure zero.

Finally, we show that the $(k-2)$--dimensional Hausdorff measure
of $\cS$ is zero (and not just finite).  First, notice that as
$\delta \to 0$, the Lebesgue measure of $\cup_i B_{r_i}(x_i)$ goes
to zero.  Since $|\nabla u|^2$ is an $L^1$ function, the dominated
convergence theorem implies that
\begin{equation}    \label{e:ugly2}
    \lim_{ \delta \to 0 } \int_{\cup_i B_{r_i}(x_i) } |\nabla u|^2 = 0 \, .
\end{equation}
Substituting this back into \eqr{e:ugly}, gives the claim.
\end{proof}

\vskip2mm  This preliminary analysis of the singular set can be
refined by doing a so--called tangent cone analysis.

\subsection{Tangent cone analysis}

Each of these variational problems comes with a natural scaling
which preserves the space of solutions.   For example, if $\Sigma
\subset \RR^n$ is a minimal submanifold, then so is
\begin{equation}
     \Sigma_{y,\lambda}  = \{ y + \lambda^{-1}\, (x-y) \, | \, x \in \Sigma \} \, .
\end{equation}
(To see this, simply note that this scaling multiplies the
principal curvatures by $\lambda$.) Similarly, if $u: \RR^n \to M$
is an energy minimizing map, then so is the map $u_{y,\lambda}$
defined by
\begin{equation}
    u_{y,\lambda} (x) = u (y+ \lambda (x-y) ) \, .
\end{equation}
Suppose now that we fix the point $y$ and take a sequence
$\lambda_j \to 0$.  The monotonicity formula (for either area or
energy) bounds the density of the rescaled solution, allowing us
to extract a convergent subsequence and limit.  This limit, which
is called a {\it tangent cone} at $y$,   achieves equality in the
monotonicity formula and, hence, must be homogeneous (i.e.,
invariant under dilations about $y$).

The usefulness of tangent cone analysis in
 regularity theory is based on two key facts.  For simplicity, we
 illustrate these when $\Sigma \subset \RR^n$ is an area
 minimizing hypersurface.  (See \cite{ChC1} for similar results for
singular limit spaces of manifolds with lower Ricci curvature bounds.)
First, if any tangent cone at $y$ is a
 hyperplane $\RR^{n-1}$, then $\Sigma$ is smooth in a neighborhood
 of $y$.  This follows easily from the Allard regularity theorem since
 the density at $y$ of the tangent cone is the same as the density
 at $y$ of $\Sigma$.  The second key fact, known as ``dimension
 reducing,'' is due to Almgren, \cite{Am2},
  and is a refinement of an argument of
 Federer.  To state this, we first stratify the singular set $\cS$
 of $\Sigma$ into subsets
 \begin{equation}
    \cS_0 \subset \cS_1 \subset \cdots \subset \cS_{n-2} \, ,
 \end{equation}
 where we define $\cS_i$ to be the set of points $y\in \cS$ so
 that  any linear space contained in any tangent cone at $y$ has
 dimension at most $i$.  (Note that $\cS_{n-1} = \emptyset$ by the
 previous fact.)  The dimension reducing argument then gives that
 \begin{equation}       \label{e:dimred}
    {\text{dim}} \, (\cS_i) \leq i \, ,
 \end{equation}
 where dimension   means the Hausdorff dimension.  In particular, the
 solution of the Bernstein problem then gives codimension $7$
 regularity of $\Sigma$, i.e., ${\text{dim}} \, (\cS) \leq n-8$.
 See lecture $2$ in \cite{Si3} for a proof of \eqr{e:dimred}.

 \begin{remark}
Using that the $(k-2)$--dimensional Hausdorff measure of $\cS$ is
zero for an energy minimizing map $u:\RR^k \to M$ by Lemma
\ref{l:std}, it is not hard to see that $\cS_{k-2} = \emptyset$.
In particular, we get that $\cS$ is at most $(k-3)$-dimensional
 \end{remark}

\vskip2mm This approach has been applied fruitfully to many
problems since it requires only a natural scaling, a monotonicity
formula, and a compactness theorem.  A variation of this, giving
tangent flows rather than tangent cones, has also been useful in
parabolic problems (see, e.g., \cite{W} for an application to mean
curvature flow).

\vskip2mm Note that the tangent cones produced in this way may
very well depend on the particular convergent subsequence.  In
some cases, one can prove uniqueness of the tangent cone and this
is often quite useful (see, for instance, section $3.4$ in
\cite{Si3} for one such application).  However, in many settings
tangent cones are not unique; see, for instance, \cite{Pe4}, \cite{ChC1}.

\section{Embedded minimal disks}        \label{s:19}

We next survey  recent results of the authors for embedded minimal
disks.  The main result is a compactness and singular convergence
theorem (Theorem \ref{t:t0.1} below) in a setting where the
classical methods cannot be applied.  In particular,  there is no
useful monotonicity formula or natural a priori bound.  The main
tools for overcoming these difficulties are a ``classification of
singularities'' which describes a neighborhood of points of large
curvature (Theorem \ref{t:blowupwinding0}) and our one-sided
curvature estimate (Theorem \ref{t:t2} below).  We will keep
things brief here, attempting to highlight a few key points.  We
refer to \cite{CM22} for a more detailed survey.

As we will see, a fundamental theorem about  embedded minimal
disks is that such a disk is either a minimal graph or can be
approximated by a piece of a rescaled helicoid depending on
whether the curvature is small or not; see Theorem \ref{t:t0.1}
below.  To avoid tedious dependence of various quantities we state
this, our main result, not for a single embedded minimal disk with
sufficiently large curvature at a given point but instead for a
sequence of such disks where the curvatures are blowing up.
Theorem \ref{t:t0.1} says that a sequence of embedded minimal
disks mimics the following behavior of a sequence of rescaled
helicoids:

\begin{quotation}
Consider the sequence $\Sigma_i = a_i \, \Sigma$ of rescaled
helicoids where $a_i \to 0$. (That is, rescale $\RR^3$ by $a_i$,
so points that used to be distance $d$ apart will in the rescaled
$\RR^3$ be distance $a_i\,d$ apart.)  The curvatures of this
sequence of rescaled helicoids are blowing up along the vertical
axis. The sequence converges (away from the vertical axis) to a
foliation by flat parallel planes. The singular set $\cS$ (the
axis) then consists of removable singularities.
\end{quotation}

Let now $\Sigma_i\subset B_{2R}$ be a sequence of embedded minimal
 disks with
$\partial \Sigma_i\subset \partial B_{2R}$.  Clearly (after
possibly going to a subsequence) either (A) or (B) occur:
\begin{enumerate}
\item[(A)]
 $\sup_{B_{R}\cap\Sigma_i}|A|^2\leq C<\infty$ for some constant $C$.
 \item[(B)]
 $\sup_{B_{R}\cap\Sigma_i}|A|^2\to \infty$.
\end{enumerate}
 In (A) (by a standard argument) the intrinsic ball
$\cB_s(y_i)$ is a graph for all $y_i\in B_{R}\cap \Sigma_i$, where
$s$ depends only on $C$. Thus the main case is (B) which is the
subject of the next theorem.

\vskip6mm Using the notion of multi--valued graphs, this the main
theorem of \cite{CM6}, can now be stated:

\begin{figure}[htbp]
    \setlength{\captionindent}{20pt}
    \begin{minipage}[t]{0.5\textwidth}
    \centering\input{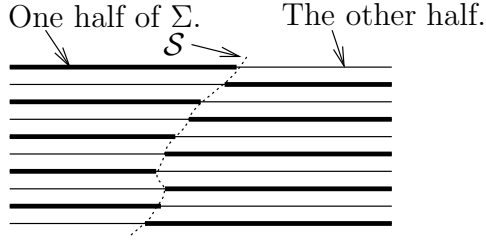}
    \caption{Theorem \ref{t:t0.1} -- the singular set, $\cS$, and
the two multi--valued graphs.}\label{f:f3}
    \end{minipage}
\end{figure}

\begin{theorem} \label{t:t0.1}
(Theorem 0.1 in \cite{CM6}).  See fig. \ref{f:f3}. Let $\Sigma_i
\subset B_{R_i}=B_{R_i}(0)\subset \RR^3$ be a sequence of embedded
minimal disks with $\partial \Sigma_i\subset \partial B_{R_i}$
where $R_i\to \infty$. If $\sup_{B_1\cap \Sigma_i}|A|^2\to
\infty$, then there exists a subsequence, $\Sigma_j$, and a
Lipschitz curve $\cS:\RR\to \RR^3$ such that after a rotation of
$\RR^3$:
\begin{enumerate}
\item
 $x_3(\cS(t))=t$.  (That is, $\cS$ is a
graph over the $x_3$--axis.)
\item  Each $\Sigma_j$
consists of exactly two multi--valued graphs away from $\cS$
(which spiral together). \item  For each $1>\alpha>0$,
$\Sigma_j\setminus \cS$ converges in the $C^{\alpha}$--topology to
the foliation, $\cF=\{x_3=t\}_t$, of $\RR^3$. \item
$\sup_{B_{r}(\cS (t))\cap \Sigma_j}|A|^2\to\infty$ for all $r>0$,
$t\in \RR$.  (The curvatures blow up along $\cS$.)
\end{enumerate}
\end{theorem}

In (2), (3) that $\Sigma_j\setminus \cS$ are multi--valued graphs
and converges to $\cF$ means that for each compact subset
$K\subset \RR^3\setminus \cS$ and $j$ sufficiently large $K\cap
\Sigma_j$ consists of multi--valued graphs over (part of)
$\{x_3=0\}$ and $K\cap \Sigma_j\to K\cap \cF$ in the sense of
graphs.

Theorem \ref{t:t0.1} (as many of the other results discussed
below) is modelled by the helicoid and its rescalings. Take a
sequence $\Sigma_i = a_i \, \Sigma$ of rescaled helicoids where
$a_i \to 0$. The curvatures of this sequence are blowing up along
the vertical axis. The sequence converges (away from the vertical
axis) to a foliation by flat parallel planes. The singular set
$\cS$ (the axis) then consists of removable singularities.

We now come to our key results for embedded minimal disks. These
are some of the main ingredients in the proof of Theorem
\ref{t:t0.1}.
 The
first says that if the curvature of such a disk $\Sigma$ is large
at some point $x\in \Sigma$, then nearby $x$
 a multi--valued graph forms (in $\Sigma$) and this extends
(in $\Sigma$) almost all the way to the boundary.  Precisely this
is:

\begin{theorem} \label{t:blowupwinding0}
(Theorem $0.2$ in \cite{CM4}).  See fig. \ref{f:f9} and fig.
\ref{f:f10}. Given $N\in \ZZ_+$, $\epsilon > 0$, there exist
$C_1,\,C_2>0$ so: Let $0\in \Sigma \subset B_{R}\subset \RR^3$ be
an embedded minimal disk, $\partial \Sigma\subset \partial B_{R}$.
If $\max_{B_{r_0} \cap \Sigma}|A|^2\geq 4\,C_1^2\,r_0^{-2}$ for
some $R>r_0>0$, then there exists (after a rotation) an
$N$--valued graph $\Sigma_g \subset \Sigma$ over $D_{R/C_2}
\setminus D_{2r_0}$ with gradient $\leq \epsilon$ and $\Sigma_g
\subset \{ x_3^2 \leq \epsilon^2 \, (x_1^2 + x_2^2) \}$.
\end{theorem}

\begin{figure}[htbp]
    \setlength{\captionindent}{20pt}
    \begin{minipage}[t]{0.5\textwidth}
\centering\input{shn3.pstex_t}
    \caption{Part 1 of the proof of Theorem \ref{t:blowupwinding0};
 finding a small
    multi--valued graph in a disk near a point of large curvature.}\label{f:f9}
\end{minipage}\begin{minipage}[t]{0.5\textwidth}
    \centering\input{shn2.pstex_t}
    \caption{Part 2 of the proof of Theorem \ref{t:blowupwinding0};
    extending a small multi--valued
graph
    in a disk.}\label{f:f10}
    \end{minipage}
\end{figure}

As a consequence of Theorem \ref{t:blowupwinding0}, one easily
gets that if $|A|^2$ is blowing up near $0$ for a sequence of
embedded minimal disks $\Sigma_i$, then there is a sequence of
$2$--valued graphs $\Sigma_{i,d}\subset \Sigma_i$, where the
$2$--valued graphs start off on a smaller and smaller scale
(namely, $r_0$ in Theorem \ref{t:blowupwinding0} can be taken to
be smaller as the curvature gets larger). Consequently, by the
sublinear separation growth, such $2$--valued graphs collapse and,
hence, a subsequence converges to a smooth minimal graph through
$0$.  To be precise, given any fixed $\rho > 0$, \eqr{e:slg}
bounds the separation $w$ at $(\rho , 0)$ by
\begin{equation}    \label{e:closeup}
    |w(\rho , 0)| \leq \left( \frac{\rho}{r_0} \right)^{\alpha} \,
    |w(r_0,0)| \leq 2 \pi \, \epsilon \, \rho^{\alpha} \,
    r_0^{1-\alpha} \, ,
\end{equation}
and this goes to $0$ as $r_0 \to 0$ since $\alpha < 1$.  The bound
$|w(r_0,0)| \leq 2 \pi \, \epsilon \, r_0$ in \eqr{e:closeup} came
from integrating the gradient bound on the graph around the circle
of radius $r_0$.  (Here $0$ is a removable singularity for the
limit.) Moreover, if the sequence of such disks is as in Theorem
\ref{t:t0.1}, i.e., if $R_i\to \infty$, then the minimal graph in
the limit is entire and hence, by Bernstein's theorem (theorem
$1.16$ in \cite{CM1}), is a plane.

\begin{figure}[htbp]
    \setlength{\captionindent}{4pt}
    \begin{minipage}[t]{0.5\textwidth}
    \centering\input{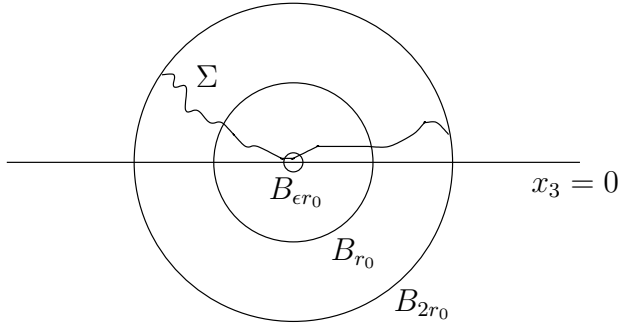}  
    \caption{Theorem \ref{t:t2} -- the one--sided curvature estimate for an
embedded minimal disk $\Sigma$ in a half--space with $\partial
\Sigma\subset \partial B_{2r_0}$:  The components of $B_{r_0}\cap
\Sigma$ intersecting $B_{\epsilon r_0}$ are graphs.}\label{f:f11}
    \end{minipage}
\end{figure}

\vskip2mm The second key result is the  curvature estimate for
embedded minimal disks in a half--space.  This theorem says
roughly that if an embedded minimal disk lies in a half--space
above a plane and comes close to the plane, then it is a graph
over the plane. Precisely, this is the following theorem:

\begin{theorem}  \label{t:t2}
(Theorem 0.2 in \cite{CM6}). See fig. \ref{f:f11}. There exists
$\epsilon>0$, such that if $\Sigma \subset B_{2r_0} \cap \{x_3>0\}
\subset \RR^3$ is an embedded minimal disk with $\partial
\Sigma\subset \partial B_{2 r_0}$, then for all components
$\Sigma'$ of $B_{r_0} \cap \Sigma$ which intersect $B_{\epsilon
r_0}$
\begin{equation}  \label{e:graph}
\sup_{\Sigma'} |A_{\Sigma}|^2 \leq r_0^{-2} \, .
\end{equation}
\end{theorem}

Using the minimal surface equation and that $\Sigma'$ has points
close to a plane, it is not hard to see that, for $\epsilon>0$
sufficiently small, \eqr{e:graph} is equivalent to the statement
that $\Sigma'$ is a graph over the plane $\{x_3=0\}$.

An embedded minimal surface $\Sigma$ which is as in Theorem
\ref{t:t2} is said to satisfy the $(\epsilon , r_0)$--{\it
effective one--sided Reifenberg condition}; cf. appendix A of
\cite{CM6} and the appendix of \cite{ChC1}. We will often refer to
Theorem \ref{t:t2} as {\it the one--sided curvature estimate}.

\begin{figure}[htbp]
    \setlength{\captionindent}{20pt}
    \begin{minipage}[t]{0.5\textwidth}
    \centering\input{pl2a.pstex_t}
    \caption{The catenoid given by revolving $x_1= \cosh x_3$
around the $x_3$--axis.}  \label{f:f12}
    \end{minipage}\begin{minipage}[t]{0.5\textwidth}
    \centering\input{unot7.pstex_t}
    \caption{Rescaling the catenoid shows that simply connected
(and embedded) is
    needed in the one--sided curvature estimate.}  \label{f:f13}
    \end{minipage}%

\end{figure}

Note that the assumption in Theorem \ref{t:t2} that $\Sigma$ is
simply connected is crucial as can be seen from the example of a
rescaled catenoid. The catenoid is the minimal surface in $\RR^3$
given by $(\cosh s\, \cos t,\cosh s\, \sin t,s)$ where
$s,t\in\RR$; see fig. \ref{f:f12}. Under rescalings this converges
(with multiplicity two) to the flat plane; see fig. \ref{f:f13}.
Likewise, by considering the universal cover of the catenoid, one
sees that embedded, and not just immersed, is needed in Theorem
\ref{t:t2}.

As an almost immediate consequence of Theorem \ref{t:t2} and a
simple barrier argument we get that if in a ball two embedded
minimal disks come close to each other near the center of the ball
then each of the disks are graphs.  Precisely, this is the
following:

\begin{figure}[htbp]
    \setlength{\captionindent}{20pt}
    \begin{minipage}[t]{0.5\textwidth}
    \centering\input{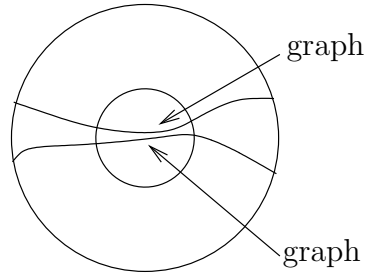}
    \caption{Corollary \ref{c:barrier}:  Two sufficiently close
    components of an embedded minimal disk must each be a graph.}\label{f:f14}
    \end{minipage}
\end{figure}

\begin{corollary}         \label{c:barrier}
(Corollary 0.4 in \cite{CM6}). See fig. \ref{f:f14}. There exist
$c > 1$, $\epsilon >0$ so: Let $\Sigma_1,\, \Sigma_2 \subset
B_{cr_0} \subset \RR^3$ be disjoint embedded minimal surfaces with
$\partial \Sigma_i \subset \partial B_{cr_0}$ and $B_{\epsilon \,
r_0 } \cap \Sigma_i \ne \emptyset$. If $\Sigma_1 $ is a disk,
 then for all components $\Sigma_1'$ of
$B_{r_0} \cap \Sigma_1$ which intersect $B_{\epsilon \, r_0}$
\begin{equation}        \label{e:onece}
    \sup_{\Sigma_1'}   |A|^2
        \leq  r_0^{-2}  \, .
\end{equation}
\end{corollary}

This estimate has also been useful in the global theory of minimal
surfaces, cf. \cite{CM9}, \cite{CM22}, and \cite{MeRo}.  It would
be very interesting to find an intrinsic version of it (i.e., for
intrinsic balls on one side of a plane):

\begin{question}
It would be very useful to prove an intrinsic version of the
one--sided curvature estimate.  Namely, does Theorem \ref{t:t2}
hold when $\Sigma$ is an intrinsic ball (in an embedded minimal
disk)?  If true, then this would likely have important
consequences for proving properness of embedded minimal disks.
\end{question}

One of the topics that we have suppressed is what we call
``properness'' of the limit.  Basically, this is proving that we
get a foliation in the limit or, equivalently, that   the points
of $\cS$ cannot be isolated. This is the one place where the
assumption $R_i \to \infty$ is used in an essential way; see
\cite{CM18} for a nonproper limit when $R_i$ does not go to
$\infty$.

\begin{question}
Suppose that a sequence of embedded minimal planar domains
$\Sigma_i \subset B_1 \subset \RR^3$ with $\partial \Sigma_i
\subset
\partial B_1$  converges away from $0$ to a minimal
lamination $\cL'$ of $B_1 \setminus \{ 0 \}$.   Does $\cL'$ extend
to a smooth lamination of $B_1$?  In other words, is $0$ a
removable singularity?  An example constructed in \cite{CM18}
shows that this need not be the case when the $\Sigma_i$'s are
disks.  It would be interesting to find non--disk examples (cf.
\cite{Ka}, \cite{Tr}).
\end{question}

\section{Global theory of minimal surfaces in $\RR^3$}      \label{s:20}

Recent years have seen breakthroughs on many long--standing
problems in the global theory of minimal surfaces in $\RR^3$. This
is an enormous subject and, rather than give a comprehensive
treatment, we will mention a few important results which fit well
with the theme of this survey. Throughout this section, $\Sigma$
will be a complete properly embedded minimal surface in $\RR^3$
(recall that properness here means that the intersection of
$\Sigma$ with any compact subset of $\RR^3$ is compact).

We say that $\Sigma$ has finite topology if it is homeomorphic to
a closed Riemann surface with a finite number of punctures; the
genus of $\Sigma$ is then the genus of this Riemann surface and
the number of punctures is the number of ends.  It follows that a
neighborhood of each puncture corresponds to a properly embedded
annular end of $\Sigma$.  Perhaps surprisingly at first, the more
restrictive case is when $\Sigma$ has more than one end. The
reason for this is that a barrier argument gives a stable minimal
surface between any pair of ends.  This stable surface is then
asymptotic to a plane (or catenoid), essentially forcing each end
to live in a half--space. Using this restriction, P. Collin
proved:

\begin{theorem}  \cite{Co}      \label{t:co}
Each end of a complete properly embedded minimal surface with
finite topology and at least two ends is asymptotic to a plane or
catenoid.
\end{theorem}

In particular, outside some compact set, $\Sigma$ is given by a
finite collection of disjoint graphs over a common plane (and has
finite total curvature).  See \cite{CM21} for a proof of Theorem
\ref{t:co} using the one--sided curvature estimate.

\vskip2mm When $\Sigma$ has only one end (e.g., for the helicoid),
it need not have finite total curvature so the situation is more
delicate.  However, the regularity results of the previous section
can be applied.  For example, if $\Sigma$ is a (non--planar)
embedded minimal disk, then we get a multi--valued graph structure
 away from a ``one--dimensional singular set.''  Using
Theorems \ref{t:t0.1} and \ref{t:t2}, W. Meeks and H. Rosenberg
proved the uniqueness of the helicoid:

\begin{theorem}     \cite{MeRo}     \label{t:mero}
The plane and helicoid are the only complete properly embedded
simply--connected minimal surfaces in $\RR^3$.
\end{theorem}

This uniqueness should have many applications. Recall that if we
take a sequence of rescalings of the helicoid, then the singular
set $\cS$ for the convergence is the vertical axis perpendicular
to the leaves of the foliation. In \cite{Me}, W. Meeks used this
fact together with the uniqueness of the helicoid to prove that
the singular set $\cS$ in Theorem \ref{t:t0.1} is always a
straight line perpendicular to the foliation.  Recently, W. Meeks
and M. Weber have constructed a {\it local} example (i.e., a
sequence of embedded minimal disks in a unit ball) where $\cS$ is
a circle.

\vskip2mm Properness is an important ingredient in many of these
results and it is not known to what extent this assumption can be
relaxed. In \cite{Na}, N. Nadirashvili constructed complete
non--proper minimal {\it immersions} (in fact, contained in a
ball).  It would be interesting to know whether this is possible
for embeddings:

\begin{question}
Suppose that $\Sigma \subset \RR^3$ is a complete embedded minimal
surface with finite topology.  Does $\Sigma$ have to be proper?
\end{question}

\vskip2mm We have not even touched on the case where $\Sigma$ has
infinite topology (e.g., when $\Sigma$ is one of the Riemann
examples). This is an area of much current research, see
\cite{CM5},  the work of Meeks, J. Perez and A. Ros,
\cite{MePRs1}, \cite{MePRs2}, and \cite{MePRs3}, and references
therein.

\vskip2mm We close this section with a local analog of the
two--ended case.  Namely, in \cite{CM9}, we proved
 that any embedded minimal annulus in a ball
(with boundary in the boundary of the ball and) with a small neck
can be decomposed by a simple closed geodesic into two graphical
sub--annuli.  Moreover, we gave a sharp bound for the length of
this closed geodesic in terms of the separation (or height)
between the graphical sub--annuli.  This  serves to illustrate our
``pair of pants'' decomposition from \cite{CM5} in the special
case where the embedded minimal planar domain is an annulus  (we
will not touch on this further here). The catenoid
$\{x_1^2+x_2^2=\cosh^2 x_3\}$ is the prime example of an embedded
minimal annulus.

The precise statement of this decomposition for annuli is:

\begin{theorem} \cite{CM9} \label{t:nitsche}
There exist $\epsilon>0$, $C_1 ,\, C_2,\, C_3>1$ so:   If
$\Sigma\subset B_{R}\subset \RR^3$ is an embedded minimal annulus
with $\partial \Sigma\subset \partial B_{R}$ and
$\pi_1(B_{\epsilon R}\cap \Sigma)\ne 0$, then there exists a
simple closed geodesic $\gamma \subset \Sigma$ of length $\ell$ so
that:

\begin{itemize}
\item
 The curve $\gamma$ splits the connected component of $B_{R/C_1}\cap
\Sigma$ containing it into
 two annuli
$\Sigma^{+} , \Sigma^{-}$ each with $\int |A|^2 \leq 5 \, \pi$.
 \item
 Furthermore, $\Sigma^{\pm}
\setminus \cT_{C_2 \, \ell}(\gamma)$ are graphs with gradient
$\leq 1$.
\item
Finally, $\ell \log (R/\ell) \leq C_3\,h$ where the separation $h$
is given by \begin{equation}
 h=\min_{x_{\pm} \in\partial B_{
R/C_1 }\cap \Sigma^{\pm}}|x_+ -x_- | \, .
\end{equation}
\end{itemize}
\end{theorem}

Here  $\cT_{s}(S) \subset \Sigma$ denotes the intrinsic
$s$--tubular neighborhood of a subset $S\subset \Sigma$.

\part{Constructing minimal surfaces and applications}   \label{p:4}

Thus far, we have mainly dealt with regularity and a priori
estimates but have ignored questions of existence.  In this part we
surveys some of the most useful existence results for minimal
surfaces and gives an application to Ricci flow.  Section
\ref{s:21} recalls the classical Weierstrass representation,
including a few modern applications, and the Kapouleas
desingularization method.  Section \ref{s:22} deals with producing
area minimizing surfaces (whether for fixed boundary, fixed
homotopy class, etc.) and questions of embeddedness.  The next
section discusses unstable (hence {\it not} minimizing) surfaces
and the corresponding questions for geodesics, concentrating on
whether the Morse index can be bounded uniformly.  Section
\ref{s:24} recalls the min--max construction for producing
unstable minimal surfaces and, in particular, doing so while
controlling the topology and guaranteeing embeddedness.  Finally,
Section \ref{s:25} discusses a recent application of min--max
surfaces to bound the extinction time for Ricci flow, answering a
question of Perelman.

\section{The Weierstrass representation}        \label{s:21}

The classical Weierstrass representation (see \cite{HoK} or
\cite{Os}) takes holomorphic data (a Riemann surface, a
meromorphic function, and a holomorphic one--form) and associates
a minimal surface in $\RR^3$.
 To be precise, given
a Riemann surface $\Omega$, a meromorphic function $g$ on
$\Omega$, and a holomorphic one--form $\phi$ on $\Omega$, then we
get  a (branched) conformal minimal immersion $F: \Omega \to
\RR^3$ by
\begin{equation}    \label{e:ws1}
    F(z) = {\text{Re }} \int_{\zeta \in \gamma_{z_0,z}}
\left( \frac{1}{2} \, (g^{-1} (\zeta) - g (\zeta) )
    , \frac{i}{2} \, (g^{-1} (\zeta)
    +g (\zeta) ) , 1 \right) \, \phi (\zeta) \, .
\end{equation}
Here $z_0 \in \Omega$ is a fixed base point and the integration is
along a path $\gamma_{z_0,z}$ from $z_0$ to $z$. The choice of
$z_0$ changes $F$ by adding a constant. In general, the map $F$
may depend on the choice of path (and hence may not be
well--defined);  this is known as ``the period problem'' (see M.
Weber and M. Wolf, \cite{WeWo}, for the latest developments).
However, when $g$ has no zeros or poles and $\Omega$ is simply
connected, then
 $F(z)$ does not depend on the choice of path
$\gamma_{z_0,z}$.

Two standard constructions of minimal surfaces from Weierstrass
data are
 \begin{align}
 &g (z) = z, \, \phi (z)  = dz/z ,
\, \Omega = \CC \setminus \{ 0 \} {\text{ giving a catenoid}}
\, , \\
  &g (z) = \e^{iz} , \,
\phi (z) = dz , \, \Omega = \CC  {\text{ giving a helicoid}} \, .
\label{e:hel}
\end{align}

The Weierstrass representation is particularly useful for
constructing immersed minimal surfaces.  For example, in
\cite{Na},  Nadirashvili used it to construct a complete immersed
minimal surface in the unit ball in $\RR^3$ (see also \cite{JXa}
for the case of a slab). In particular, Nadirashvili's surface is
not proper, i.e., the intersections with compact sets are not
necessarily compact.

 Typically, it is rather difficult to prove that the resulting
immersion is an embedding (i.e., is $1$--$1$), although there are
some interesting cases where this can be done.    The first modern
example was \cite{HoMe} where D. Hoffman and Meeks proved that the
surface constructed by Costa was embedded; this was the first new
complete finite topology properly embedded minimal surface
discovered since the classical catenoid, helicoid, and plane. This
led to the discovery of many more such surfaces
 (see \cite{HoK} and \cite{Ro} for more discussion).

In \cite{CM18}, we used the Weierstrass representation to
construct a sequence of embedded minimal disks $\Sigma_i \subset
B_1 = B_1 (0)\subset \RR^3$ with $\partial \Sigma_i \subset
\partial B_1$ where
 the curvatures blow up only at $0$  and $\Sigma_i \setminus \{ {\text{$x_3$--axis}} \}$
consists of two multi--valued graphs for each $i$. Furthermore,
$\Sigma_i  \setminus \{ x_3 = 0 \}$ converges to two
  embedded minimal disks $\Sigma^- \subset \{ x_3 < 0 \}$
and  $\Sigma^+ \subset \{ x_3 > 0 \}$ each of which spirals
 into $\{ x_3 =
0 \}$ and thus is not proper.  (This should be contrasted with
Theorem \ref{t:t0.1} where the {\it complete} limits are planes
and hence proper.)

\vskip2mm N. Kapouleas has developed another method to construct
  complete embedded minimal surfaces with finite total
curvature.  For instance, in \cite{Ka}, he shows that (most)
collections of coaxial   catenoids and planes can be
desingularized to get complete embedded minimal surfaces with
finite total curvature. The Costa surface above had genus one and
three ends (that is to
 say, it is homeomorphic to a torus with three punctures).  In the
 Kapouleas construction, one could start with a plane and catenoid
 intersecting in a circle and then desingularize this circle using
 suitably scaled and bent Scherk surfaces
 to get a finite genus embedded surface with three
 ends.  (This desingularization process adds handles, i.e.,
 increases the genus.)  In this manner, Kapouleas gets an enormous
 number of new examples; see also the gluing construction of S.D. Yang,
 \cite{Y}, which uses catenoid necks to glue together nearby
 minimal surfaces.

\section{Area--minimizing surfaces}     \label{s:22}

Perhaps the most natural way to construct minimal surfaces is to
look for ones which minimize area, e.g., with fixed boundary, or
in a homotopy class, etc.  This has the advantage that often it is
possible to show that the resulting surface is embedded.  We
mention a few results along these lines.

The first embeddedness result, due to   Meeks and Yau, shows that
if the boundary curve is embedded and lies on the boundary of a
smooth mean convex set (and it is null--homotopic in this  set),
then it bounds an embedded least area disk.

\begin{theorem}  \cite{MeYa1}             \label{t:my1}
Let $M^3$ be a compact Riemannian three--manifold whose boundary
is mean convex and let $\gamma$ be a simple closed curve in
$\partial M$ which is null--homotopic in $M$; then $\gamma$ is
bounded by a least area disk and any such least area disk is
properly embedded.
\end{theorem}

Note that some restriction on the boundary curve $\gamma$ is
certainly necessary. For instance, if the boundary curve was
knotted (e.g., the trefoil), then it could not be spanned by any
embedded disk (minimal or otherwise).  Prior to the work of Meeks
and Yau, embeddedness was known for extremal boundary curves in
$\RR^3$ with small total curvature by the work of R. Gulliver and
J. Spruck \cite{GuSp}; see chapter $4$ in \cite{CM1} for other
results and further discussion.

 If we instead fix a homotopy class of maps, then the two fundamental existence results
 are due to Sacks--Uhlenbeck and Schoen--Yau (with embeddedness proven
 by Meeks--Yau and Freedman--Hass--Scott,
 respectively):

\begin{theorem}     \label{p:existence1} \cite{SaUh}, \cite{MeYa2}
Given $M^3$, there exist conformal  (stable) minimal immersions
$u_1 , \dots , u_m  : \SS^2 \to M$ which generate $\pi_2 (M)$ as a
$\ZZ[\pi_1 (M)]$ module.  Furthermore,
\begin{itemize}
\item
If $u: \SS^2 \to M$ and $[u]_{\pi_2} \ne 0$, then $\Area (u) \geq
\min_i \Area (u_i)$.
 \item Each $u_i$ is either an
embedding or a $2$--$1$ map onto an embedded $2$--sided $\RP^2$.
\end{itemize}
\end{theorem}

\begin{theorem}     \label{p:existence2} \cite{ScYa2}, \cite{FHS}
If $\Sigma^2$ is a closed surface with genus $g>0$ and $i_0 :
\Sigma \to M^3$ is an embedding which induces an injective map on
$\pi_1$, then there is a least area embedding with the same action
on $\pi_1$.
\end{theorem}

In \cite{MeSiYa}, Meeks, Simon, and Yau find an embedded sphere
minimizing area in an isotopy class in a closed $3$--manifold.

\vskip2mm We end this section by mentioning two applications of
Theorem \ref{p:existence2}.  First, in \cite{CM20}, we showed that
any topological $3$--manifold $M$ had an open set of metrics so
that, for each such metric, there was a sequence of embedded
minimal tori whose area went to infinity.  In \cite{De}, B. Dean
showed that this was true for every genus $g \geq 1$.  This leaves
an obvious interesting question:

\begin{question}
Given a topological $3$--manifold $M$, does there exist an open
set of metrics which have embedded minimal spheres with
arbitrarily large area?
\end{question}

It would be interesting  to answer this question even when the
minimal spheres are stable (the examples constructed in
\cite{CM20} and \cite{De} were all locally minimizing and hence
also stable)

\section{Index bounds for geodesics and minimal surfaces}   \label{s:23}

The minimal surfaces discussed in the previous section were all stable
and in fact locally area minimizing.  This is a very special and strong
property of a minimal surface.  In general, like, for instance the
helicoid and the catenoid, most minimal surfaces are not stable but
have non--zero index.  In this section we will discuss the Morse index
of simple
closed geodesics on surfaces and of embedded minimal surfaces in
$3$--manifolds.  First let us discuss the situation of simple closed
geodesics in surfaces.

Let $M^2$ be a closed orientable surface with curvature $K$
and $\gamma\subset M$
a closed geodesic.  The {\it Morse index} of
$\gamma$ is the index of the
critical point $\gamma$ for the length functional, i.e., the number of
negative
eigenvalues (counted with multiplicity) of the second
derivative of length (throughout curves will always be in $H^1$).
Since the second derivative of length at $\gamma$
in the direction of a normal variation $u\,\nn$ is
$-\int_{\gamma}u\,L_{\gamma}\,u$ where $L_{\gamma} \,u= u'' + K\,u$,
the Morse index is the number of
negative eigenvalues of $L_{\gamma}$.
(By convention,  an
eigenfunction $\phi$ with eigenvalue
$\lambda$ of $L_{\gamma}$ is a solution of
$L_{\gamma}\,\phi+\lambda\, \phi=0$.)  Note that if $\lambda=0$,
then $\phi$ (or $\phi\,\nn$) is a (normal) Jacobi field.
$\gamma$ is {\it stable}
if the index is zero.
The {\it index} of a noncompact geodesic is the dimension
of a maximal vector space of compactly supported variations for which the
second derivative of length is negative definite.  We also say that
such a geodesic
is {\it stable} if the index is $0$.

As the following result shows then it turns out that in general
there are no Morse index bounds for simple closed geodesics on
surfaces.

\begin{theorem}  \label{t:example}
\cite{CH1}.
On any $M^2$,
there exists a metric with a geodesic lamination
with infinitely many unstable leaves.  Moreover, there is such a metric
with simple closed geodesics of arbitrary high Morse index.
\end{theorem}

A codimension one {\it lamination}
on a surface $M^2$ is a collection $\cL$ of
smooth disjoint curves (called leaves)
such that
$\cup_{\ell \in \cL} \ell$ is closed.
Moreover, for each $x\in M$ there exists an
open neighborhood $U$ of $x$ and a $C^0$ coordinate chart, $(U,\Phi)$, with
$\Phi (U)\subset \RR^2$
so that in these coordinates the leaves in $\cL$
pass through $\Phi (U)$ in slices of the
  form $(\RR\times \{ t\})\cap \Phi(U)$.
A {\it foliation} is a lamination for which
the union of the leaves is all of $M$
and a {\it geodesic lamination} is a lamination whose leaves are geodesics.

Similarly, to the geodesic case, for an immersed minimal surface
$\Sigma$ in a $3$--manifold $M$, we set $L_{\Sigma}\,\phi
=\Delta_{\Sigma}\,\phi+|A|^2\,\phi+\Ric_M(\nn,\nn)\,\phi$.
(Note that by the second variational formula (see, for instance,
section 1.7 of \cite{CM1}), then
\begin{equation}
\frac{\partial^2}{\partial r^2}_{r=0}\Area (\Sigma_r)
=-\int_{\Sigma}\phi\,L_{\Sigma}\,\phi\, ,
\end{equation}
where $\Sigma_r=\{x+r\,\phi (x)\,\nn_{\Sigma} (x)\,|\,x\in
\Sigma\}$.)  Recall also that by definition the index of a minimal
surface $\Sigma$ is the number of negative eigenvalues (counted with
multiplicity) of $L_{\Sigma}$.  (A function $\eta$ is an eigenfunction
of $L_{\Sigma}$ with eigenvalue $\lambda$ if
$L_{\Sigma}\,\eta+\lambda\,\eta=0$.)   Thus in
particular, since $\Sigma$ is assumed to be closed, the index is always finite.

Theorem \ref{t:example} was proven by first constructing a
metric on the disk with convex boundary having no Morse index bounds
and then completing the metric to a metric on the given
$M^2$.  By taking the product of this metric on the disk
with a circle we get, on a solid torus,
a metric with convex boundary and without Morse index bounds for
embedded minimal tori, and with a minimal lamination with infinitely
many unstable leaves.  By completing this metric we get:

\begin{theorem}  \label{t:example2}
\cite{CH2}
On any $M^3$,
there exists a metric with a minimal lamination
with infinitely many unstable leaves.  Moreover, there is such a metric
with embedded minimal tori of arbitrary high Morse index.
\end{theorem}

By construction the embedded minimal tori in Theorem \ref{t:example2} and the
leaves of the lamination can be taken to be totally geodesic.

We will equip the space of metrics on a given manifold with the
$C^{\infty}$-topology.  A subset of the set of metrics on the manifold is said
to be {\it residual} if it is a countable intersection of open dense
subsets.
A metric on a surface is {\it bumpy} if each  closed geodesic
is a nondegenerate critical point, i.e.,
$L_{\gamma} u = 0$ implies $u\equiv 0$.  It follows from
results of Abraham and Anosov that bumpy metrics are generic;
that is the set of bumpy metrics contain a
residual set.

To check that any given metric is bumpy is virtually impossible; however
it seems that the metric in Theorem \ref{t:example} can be chosen
to be bumpy.  Thus it seems
unlikely that a bumpy metric is enough to ensure a bound for the Morse
index of simple closed geodesics on $M^2$.  What is needed is a
nondegeneracy condition for noncompact simple geodesics, rather
than one for closed geodesics; cf. \cite{CH2}, \cite{CH3}.

In \cite{HaNoRu} examples were given of metrics on any $M^3$ that
have embedded minimal spheres without bounds and in \cite{CD2} the
following was shown:  For any $3$-manifold $M^3$ and any
nonnegative integer ${\bf{g}}$, there are examples of metrics on
$M$ each of which has a sequence of embedded minimal surfaces of
genus ${\bf{g}}$ and without Morse index bounds. On any spherical
space form $\SS^3/\Gamma$  \cite{CD2} constructed such a metric
with positive scalar curvature. More generally \cite{CD2}
constructed such a metric with $\text{Scal} >0$ (and such
surfaces) on any $3$-manifold which carries a metric with
$\text{Scal}>0$. In all but one of the examples in \cite{CD2} the
Hausdorff limit is a singular minimal lamination.  The
singularities being in each case exactly two points lying on a
closed leaf (the leaf is a strictly stable sphere).

\cite{CD2} used in part ideas
of Hass-Norbury-Rubinstein
\cite{HaNoRu}.  As in \cite{HaNoRu}, but
unlike the examples in \cite{CH1},
the surfaces in \cite{CD2} have no uniform curvature bounds.
In fact, it follows easily (see appendix  B of \cite{CM4}) that
 if
$\Sigma_i\subset M^3$ is a sequence of embedded  minimal surfaces with
uniformly
bounded curvatures, then a subsequence  converges to a smooth lamination.
Moreover, with the right notion of being generic,
the following seems likely (by \cite{CH1} bumpy is
not the right generic notion):

\begin{question}
Let $M^3$ be a closed $3$-manifold with a
generic metric and $\Sigma_i\subset M$ a sequence of embedded minimal
surfaces of a given genus.  If any limit
of the $\Sigma_i$'s is a
{\underline{smooth}} (minimal) lamination, then the sequence $\Sigma_i$ has a
uniform Morse index bound.
\end{question}

A codimension one {\it lamination} of $M^3$ is a collection $\cL$ of
smooth disjoint connected surfaces (called leaves)
such that  $\cup_{\Lambda \in \cL} \Lambda$ is closed.
Moreover, for each $x\in M$ there exists an
open neighborhood $U$ of $x$ and a local coordinate chart, $(U,\Phi)$, with
$\Phi (U)\subset \RR^3$
such that in these coordinates the leaves in $\cL$
pass through the chart in slices of the
form $(\RR^2\times \{ t\})\cap \Phi(U)$.

A lamination is said to be minimal if the leaves are (smooth)
minimal surfaces.  If the union of the leaves  is all of $M$, then it
is a foliation.

There are two results that support this question.
The first concerns the corresponding question
in one dimension less (that is for geodesics on surfaces); see
\cite{CH2}, \cite{CH3}.
The second concerns the question for
$3$-manifolds with
positive scalar curvature.  However, there are
examples where the limit is not smooth; see \cite{CD2}.

Finally, we refer to \cite{CH2} and \cite{CH3} for further
discussion of Morse index bounds for geodesics including some
positive results about when one has such bounds.

\section{The min--max construction of minimal surfaces} \label{s:24}

Variational arguments can also be used to construct higher index
(i.e., non--minimizing) minimal surfaces using   the topology of
the space of surfaces. There are two basic approaches:
\begin{itemize}
\item
Applying Morse theory to the energy functional on the space of
maps from a fixed surface $\Sigma$ to $M$.
\item
Doing a min--max argument over  families of (topologically
non--trivial) sweep--outs of $M$.
\end{itemize}
The first approach has the advantage that the topological type of
the minimal surface is easily fixed; however, the second approach
has been  more successful at producing embedded minimal surfaces.
We will highlight a few key  results below but refer to \cite{CD1}
for a thorough treatment.

Unfortunately, one cannot directly apply Morse theory to the
energy functional on the space of maps from a fixed surface
because of a lack of compactness (the Palais--Smale Condition C
does not hold).  To get around this difficulty, \cite{SaUh}
introduce a family of perturbed energy functionals which do
satisfy Condition C and then obtain   minimal surfaces as limits
of critical points for the perturbed problems:

\begin{theorem} \cite{SaUh}     \label{t:nonasph}
If   $\pi_k (M) \ne 0$ for some $k>1$, then there exists a
branched immersed minimal $2$--sphere in $M$ (for any metric).
\end{theorem}

This was sharpened somewhat by \cite{MiMo} (showing that the index
of the minimal sphere was at most $k-2$), who used it to prove a
generalization of the sphere theorem.  See \cite{Jo} and \cite{St}
for approaches which avoid using the perturbed functionals and
\cite{Fr} for a generalization to a free boundary problem.

\begin{figure}[htbp]
    \setlength{\captionindent}{4pt}
    \begin{minipage}[t]{0.5\textwidth}
    \centering\input{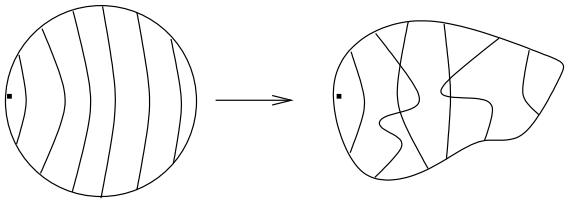}
    \caption{A $1$--parameter family of curves on a $2$--sphere which
induces a map $F:\SS^2 \to
    \SS^2$ of degree $1$.}\label{f:fsweep}
    \end{minipage}
\end{figure}

\vskip2mm The basic idea of constructing minimal surfaces via
min--max arguments and sweep--outs goes back to Birkhoff, who
developed it to construct simple closed geodesics on spheres.  In
particular, when $M$ is a topological $2$--sphere, we can find a
$1$--parameter family of curves starting and ending at point
curves so that the induced map $F:\SS^2 \to \SS^2$ (see fig.
\ref{f:fsweep}) has nonzero degree.    The min--max argument
produces a nontrivial closed geodesic of length less than or equal
to the longest curve in the initial one--parameter family.  A
curve shortening argument gives that the geodesic obtained in this
way is simple.

In \cite{Pi}, J. Pitts applied a similar argument and   geometric
measure theory to get that every closed Riemannian three manifold
has an embedded minimal surface (his argument was for dimensions
up to seven), but he did not estimate the genus of the resulting
surface.  Finally, F. Smith (under the direction of
L. Simon) proved (see \cite{CD1}):

\begin{theorem} \cite{Sm}   \label{t:smith}
Every metric   on  a topological $3$--sphere $M$ admits an
embedded minimal $2$--sphere.
\end{theorem}

The main new contribution of Smith was to control the
topological type of the resulting minimal surface while keeping it
embedded; see also Pitts and Rubinstein, \cite{PiRu}, for some generalizations.

\section{An application of min--max surfaces to Ricci flow} \label{s:25}

We will in this section give an application of the min--max construction
of minimal surfaces to the Ricci flow.  This application
is that on $\SS^3$ starting
at any given metric the Ricci flow becomes extinct in finite time.  The
treatment here follows \cite{CM19} and was inspired by a question of
Perelman; see \cite{CM19}, \cite{Pe3}.  (See also paragraph 4.4 of
\cite{Pe2} for the precise definition of extinction time in the case
that surgery occurs.)

Throughout this section we let $M^3$ be the $3$--sphere and let
$g(t)$ be a one--parameter family of metrics on $M$ evolving by
the Ricci flow, so
\begin{equation}  \label{e:eqRic2}
 \partial_t g=-2\,\Ric_{M_t}\, .
\end{equation}

Since $\pi_3(M) = \ZZ$, it follows from suspension, as in lemma
$3$ of \cite{MiMo}, that the space of maps from $\SS^2$ to $M$ is
not simply connected.

Fix a continuous map $\beta: [0,1] \to C^0\cap L_1^2
(\SS^2 , M)$ where $\beta (0)$ and $\beta (1)$ are constant maps
so that $\beta$ is in  the nontrivial homotopy class $[\beta]$. We
define the width $W=W(g,[\beta])$ by
\begin{equation}    \label{e:w3}
   W(g) = \min_{\gamma \in [ \beta]}
\, \max_{s \in [0,1]} \Energy (\gamma(s)) \,
   .
\end{equation}

One could equivalently define the width using the area rather than
the energy, but  the energy is somewhat easier to work with.  As
for the Plateau problem, this equivalence follows using the
uniformization theorem and the inequality $\Area (u) \leq \Energy
(u)$  (with equality when $u$ is a branched conformal map); cf.
lemma $4.12$ in \cite{CM1}.

\vskip2mm The next theorem gives an upper bound for the derivative
of $W(g(t))$ under the Ricci flow which forces the solution $g(t)$
to become extinct in finite time.

\begin{theorem}     \label{t:upper}
(\cite{CM19} and cf. \cite{Pe3}).
Let $M^3$ be the $3$--sphere equipped with a Riemannian metric $g=g(0)$. Under
the Ricci flow, the width $W(g(t))$ satisfies
\begin{equation}   \label{e:di1a}
\frac{d}{dt} \, W(g(t))  \leq -4 \pi + \frac{3}{4 (t+C)} \,
W(g(t))   \, ,
\end{equation}
in the sense of the limsup of forward difference quotients. Hence,
 $g(t)$  must become extinct in finite time.
\end{theorem}

 Suppose that $\Sigma\subset M$ is a closed immersed
 surface (not necessarily minimal), then using \eqr{e:eqRic2}
an easy calculation
gives (cf. page 38--41 of \cite{Ha3})
\begin{equation}        \label{e:diffAn}
\frac{d}{dt}_{t=0}\Area_{g(t)}( \Sigma ) =-\int_{\Sigma} [R -
\Ric_M (\nn,\nn)] \, .
\end{equation}
If $\Sigma$ is also minimal, then
\begin{align}        \label{e:diffA}
\frac{d}{dt}_{t=0}\Area_{g(t)}( \Sigma )
&=-2\int_{\Sigma}\text{K}_{\Sigma} -\int_{\Sigma}[|A|^2+\Ric_M
(\nn,\nn)] \\
&=-\int_{\Sigma}K_{\Sigma}-\frac{1}{2}\int_{\Sigma}[ |A|^2 + R] \,
.   \notag
\end{align}
Here $\text{K}_{\Sigma}$ is the (intrinsic) curvature of
$\Sigma$, $\nn$ is a unit
normal for $\Sigma$ (our $\Sigma$'s below will be $\SS^2$'s and
hence have a well--defined unit normal), $A$ is the second
fundamental form of $\Sigma$ so that $|A|^2$ is the sum of the
squares of the principal
curvatures, $\Ric_M $ is the Ricci curvature of
$M$, and $R$ is the scalar curvature of $M$. (The curvature
is normalized so that on the unit $\SS^3$ the Ricci curvature is
$2$ and the scalar curvature is $6$.) To get \eqr{e:diffA}, we
used that by the Gauss equations and minimality of $\Sigma$
\begin{equation}
\text{K}_{\Sigma}=\text{K}_M-\frac{1}{2}|A|^2\, ,
\end{equation}
where $\text{K}_M$ is the sectional curvature of $M$ on the two--plane
tangent to $\Sigma$.

Our first lemma gives an upper bound for the rate of change of
area of minimal $2$--spheres.

\begin{lemma}     \label{l:upper}
If $\Sigma\subset M^3$ is a branched minimal immersion of the
$2$--sphere, then
\begin{equation}
\frac{d}{dt}_{t=0}\Area_{g(t)}(\Sigma)  \leq -4 \pi -
\frac{\Area_{g(0)}(\Sigma)}{2} \, \min_{M} R(0) \, .
\end{equation}
\end{lemma}

\begin{proof}
Let $\{ p_i \}$ be the set of branch points of $\Sigma$ and $b_i >
0$ the order of branching at $p_i$.
 By \eqr{e:diffA}
\begin{equation}
 \frac{d}{dt}_{t=0} \Area_{g(t)}( \Sigma )
\leq -\int_{\Sigma}K_{\Sigma}-\frac{1}{2}\int_{\Sigma}R = -4\pi
-2\pi
 \sum b_i -\frac{1}{2}\int_{\Sigma}R  \, ,
\end{equation}
where the  equality used the Gauss--Bonnet theorem  with branch
points.
\end{proof}

We will need to recall a  result of J.
Jost, theorem $4.2.1$ of
\cite{Jo}, which gives the existence of minimal spheres which
realize the width $W(g)$. (The bound for the index is not stated
explicitly in \cite{Jo} but follows immediately as in
\cite{MiMo}.)

\begin{proposition}     \label{p:existence}
Given a metric $g$ on $M$ and a nontrivial $[\beta] \in \pi_1
(C^0\cap L^2_1 (\SS^2 , M))$, there exists  a sequence of
sweep--outs $\gamma^j: [0,1] \to C^0\cap L_1^2 (\SS^2 , M)$ with
$\gamma^j \in [\beta]$ so that
\begin{equation}
    W(g) = \lim_{j \to \infty} \, \max_{s \in [0,1]} \, \Energy
        (\gamma^j_s)   \, .
\end{equation}
Furthermore, there exist $s_j \in [0,1]$ and branched conformal
minimal immersions $u_0 , \dots , u_m : \SS^2 \to M$ with index at
most one so that, as $j \to \infty$, the maps $\gamma^j_{s_j}$
converge to $u_0$ weakly in $L_1^2$ and uniformly on compact
subsets of $\SS^2 \setminus \{ x_1 , \dots , x_k \}$, and
\begin{equation}
    W(g) =   \sum_{i=0}^m \Energy (u_i) = \lim_{j \to \infty} \Energy
        (\gamma^j_{s_j}) \,  .
\end{equation}
 Finally,
for each $i > 0$, there exists a point $x_{k_i}$ and a sequence of
conformal dilations $D_{i,j} : \SS^2 \to \SS^2$  about $x_{k_i}$
so that the maps $\gamma^j_{s_j} \circ D_{i,j}$ converge to $u_i$.
\end{proposition}

 We will also need a  standard additional property for
 the min--max sequence of
  sweep--outs $\gamma^j$ of Proposition \ref{p:existence}
which can be achieved by modifying
  the sequence   as in section $4$ of \cite{CD1}
(cf. proposition $4.1$ on page 85 in \cite{CD1}).
  Loosely speaking this is the property that any subsequence
$\gamma_{s_k}^k$ with energy converging to $W(g)$ converges (after
possibly going to a further subsequence) to the union of branched
immersed minimal $2$--spheres each with index at most one.
Precisely this is that we can choose $\gamma^j$ so that:
 Given $\epsilon >
 0$, there exist $J$ and $\delta > 0$ (both depending on $g$ and $\gamma^j$)
so that if $j > J$ and
\begin{equation}    \label{e:eclose}
      \Energy (\gamma^j_s) > W(g) - \delta \, ,
\end{equation}
then there is a collection of branched minimal $2$--spheres $\{
\Sigma_i \}$ each of index at most one and with
\begin{equation}    \label{e:vclose}
    \dist  \, (\gamma^j_s ,
    \cup_i \Sigma_i ) < \epsilon \, .
\end{equation}
Here, the distance means  varifold distance (see, for instance,
section $4$ of \cite{CD1}). Below we will use that, as an immediate
consequence of \eqr{e:vclose}, if $F$ is a quadratic form on $M$
and $\Gamma$ denotes $\gamma^j_s$, then
\begin{equation}     \label{e:hclose}
   \left| \int_{  \Gamma } [\text{Tr} (F) - F(\nn_{\Gamma} ,
   \nn_{\Gamma})] - \sum_i \int_{\Sigma_i}
   [\text{Tr} (F) - F(\nn_{\Sigma_i} ,
   \nn_{\Sigma_i})]  \right| <   C \, \epsilon \, \| F \|_{C^1}
\, \Area (\Gamma) \, .
\end{equation}

\vskip2mm
\begin{proof}
(of Theorem \ref{t:upper})   Fix a time $\tau$.  Below $\tilde{C}$
denotes a constant depending only on $\tau$ but will be allowed to
change from line to line.  Let $\gamma^j (\tau)$ be the sequence
of sweep--outs for the metric $g(\tau)$ given by
 Proposition \ref{p:existence}. We will use the sweep--out at
 time $\tau$ as a comparison to get an upper bound for the width
 at times $t > \tau$.   The key for this is the following claim
(the first inequality
 in \eqr{e:acomp1} below):
Given $\epsilon > 0$, there exist $J$ and $\bar{h} > 0$ so that if
$j> J$ and $0 < h < \bar{h}$, then
\begin{align}        \label{e:acomp1}
    \Area_{g(\tau + h)}( \gamma^j_s (\tau ) )
    &- \max_{s} \, \Energy_{g(\tau)}( \gamma^j_s (\tau ))  \notag
    \\
    &\leq
     [-4 \pi + \tilde{C} \, \epsilon
- \frac{ \max_{s} \, \Energy_{g(\tau)}( \gamma^j_s (\tau ))}{2} \, \min_{M}
R(\tau) ] \, h + \tilde{C} \, h^2   \notag \\
&\leq
     [-4 \pi + \tilde{C} \, \epsilon +
 \frac{3}{4 (\tau+C)} \, \max_{s} \,
\Energy_{g(\tau)}( \gamma^j_s (\tau )) ] \, h + \tilde{C} \, h^2 \, ,
\end{align}
where the second inequality used the lower bound \eqr{e:scalar}
for $R(\tau)$.
 To see why
\eqr{e:acomp1} implies \eqr{e:di1a}, we use the definition of the
width to get
\begin{equation}        \label{e:defw}
    W(g (\tau + h) ) \leq \max_{s \in [0,1]}
    \Area_{g(\tau + h)}( \gamma^j_s (\tau ) ) \, ,
\end{equation}
and then  take the limit as $j\to \infty$ (so that $\max_{s} \,
\Energy_{g(\tau)}( \gamma^j_s (\tau ))\to W(g(\tau))$) in
\eqr{e:acomp1} to get
\begin{equation}        \label{e:defwq}
    \frac{W(g (\tau + h) ) - W(g (\tau ))}{h}
\leq    -4 \pi + \tilde{C} \, \epsilon +
 \frac{3}{4 (\tau+C)} \, W(g(\tau))     + \tilde{C} \, h
     \, .
\end{equation}
Taking $\epsilon \to 0$ in \eqr{e:defwq} gives \eqr{e:di1a}.
  It remains to prove
\eqr{e:acomp1}. First,  let $\delta > 0$ and $J$, depending on
$\epsilon$ (and on $\tau$), be given by \eqr{e:eclose}--\eqr{e:hclose}.
If $j > J$ and $\Energy_{g(\tau)} (\gamma^j_s (\tau )  ) >  W(g) -
\delta$, then let $\cup_i \Sigma_{s,i}^j (\tau) $ be the
collection of minimal spheres in \eqr{e:hclose}. Combining
\eqr{e:diffAn}, \eqr{e:hclose} with $F = \Ric_M$, and Lemma
\ref{l:upper} gives
\begin{align}        \label{e:diffAn2}
    \frac{d}{dt}_{t=\tau}\Area_{g(t)}( \gamma^j_s (\tau )) &
    \leq \frac{d}{dt}_{t=\tau}\Area_{g(t)}( \cup_i \Sigma_{s,i}^j (\tau) )
    +  \tilde{C} \, \epsilon \, \| \Ric_M \|_{C^1}
\, \Area_{g(t)}( \gamma^j_s (\tau ))\notag \\
        &\leq -4 \pi  - \frac{\max_{s} \,
\Energy_{g(\tau)}( \gamma^j_s (\tau ))}{2} \, \min_{M}
R(\tau) + \tilde{C} \, \epsilon \, .
\end{align}
Since the metrics $g(t)$ vary smoothly and every sweep--out
$\gamma^j$ has uniformly bounded energy, it is easy to see that
$\Energy_{g(\tau + h)} (\gamma^j_s (\tau )  )$ is a smooth
function of $h$ with a uniform $C^2$ bound independent of both $j$
and $s$ near $h=0$ (cf. \eqr{e:diffAn}).  In particular,
\eqr{e:diffAn2} and Taylor expansion gives $\bar{h} > 0$
 (independent of $j$) so that \eqr{e:acomp1} holds
for  $s$ with $\Energy_{g(\tau)} (\gamma^j_s (\tau )  ) >  W(g) -
\delta$.  In
 the remaining case, we have  $\Energy (\gamma^j_s (\tau )) \leq
W(g) - \delta$ so the continuity of $g(t)$ implies that
\eqr{e:acomp1} automatically holds after possibly shrinking
$\bar{h}> 0$.

Finally, we claim that \eqr{e:di1a} implies finite
extinction time.   Namely, rewriting \eqr{e:di1a} as $\frac{d}{dt} \left(
W(g(t)) \, (t+C)^{-3/4} \right) \leq - 4\pi \, (t+C)^{-3/4}$ and
integrating gives
\begin{equation}  \label{e:last}
 (T+C)^{-3/4} \, W(g(T)) \leq C^{-3/4} \, W(g(0))
- 16 \, \pi \, \left[ (T+C)^{1/4} - C^{1/4} \right]   \, .
\end{equation}
Since $W \geq 0$ by definition and the right hand side of \eqr{e:last} would
become
negative for $T$ sufficiently large, the theorem follows.
\end{proof}

\part{Growth  of harmonic functions}    \label{p:5}

We next discuss some global results for harmonic functions and a
few applications of function theory. We will focus on the function
theory of manifolds with non--negative Ricci curvature, with the
exception of Section \ref{s:29} where we discuss two estimates
related to nodal sets of eigenfunctions.

Recall that the classical Liouville theorem states that any
bounded (or even just positive) harmonic function is constant on
Euclidean space. In fact, the Euclidean gradient estimate shows
that a nonconstant harmonic function must grow at least linearly.
Since the partial derivatives of a Euclidean harmonic function are
again harmonic, iterating this gives that on Euclidean space any
harmonic function of polynomial growth is a harmonic polynomial.
 In particular, the dimensions
of these spaces are finite on Euclidean space.

The picture gets quite a bit more complicated when we look at more
general manifolds.  For example, one can prescribe asymptotic
values  at infinity on hyperbolic space (cf. \cite{An3}), so that
even the space of bounded harmonic functions is infinite
dimensional in this case.

However,  \cite{CM14} proved  that each space of harmonic
functions of polynomial growth is  finite dimensional for
manifolds with non--negative Ricci curvature. (This had been
conjectured by S.T. Yau; see \cite{Ya1}, \cite{Ya2}.  The case of
surfaces was settled in \cite{LiTa2}.). An interesting feature of
\cite{CM14} was that only two properties were used: a volume
doubling and a Neumann Poincar\'e inequality, cf. \ref{ss:72}.

\vskip2mm Given an open manifold $M$ and $d> 0$, we define the
spaces of harmonic functions of polynomial growth of order at most
$d$, $\cH_d (M)$, using the distance function from a fixed point
$p$:

\begin{definition}   \rm
A function $u$ is in  $\cH_d (M)$ if $u$ is harmonic on $M$ and
\begin{equation}
    |u(x)| \leq C (1 +\dist_M (p,x)^d) \, ,
\end{equation}
 for some constant $C$ and point $p \in M$.
\end{definition}

\section{Harmonic functions and spherical harmonics}        \label{s:26}

 It  is worthwhile to recall the Euclidean case using polar
coordinates $(\rho , \theta)$, where $\theta \in \SS^{k-1}$. In
these coordinates, the Laplacian is
\begin{equation}
\Delta_{\RR^k}     = \rho^{-2} \Delta_{\theta}
        + (k-1)\,\rho^{-1} \frac{\partial}{\partial \rho}
                + \frac{\partial^2}{\partial \rho^2}  \, .
\end{equation}
In particular, the restriction of a homogeneous harmonic
polynomial of degree $d$ to $\SS^{k-1}$ gives an eigenfunction
with eigenvalue $d^2 + (k-2)d$. It is then not hard to see that
understanding $\cH_d (\RR^k)$ is a spectral problem on the compact
manifold $\SS^{k-1}$.

  A similar ``cone construction'' holds more generally.
 Given a  manifold $N^{k-1}$,  the {\it cone} over $N$ is
 the manifold $C(N) = N \times [0,\infty)$ with the metric
\begin{equation}
    ds_{C(N)}^2 = dr^2 + r^2 \, ds^2_N \, .
\end{equation}
 Usually we identify $N \times \{ 0 \}$ and refer to this point as
 the vertex.   A direct computation shows that the Laplacians of $N$ and $C(N)$
are related by the following simple formula at $(x,r) \in N \times
(0,\infty)$:
\begin{equation}                \label{e:sepv}
        \Delta_{C(N)} u   = r^{-2} \Delta_N
                u
        + (k-1)\,r^{-1} \frac{\partial}{\partial r} u
                + \frac{\partial^2}{\partial r^2} u \, .
\end{equation}
Using \eqr{e:sepv}, we can now reinterpret the spaces $\cH_d
(C(N))$:

\begin{lemma}   \label{l:equiv}
  If $\Delta_N g = - \lambda \, g$   on $N^{k-1}$, then
  $r^p \, g \in \cH_p (C(N))$ where
  \begin{equation}
    p^2 +       (k-2) p = \lambda \, .
\end{equation}
\end{lemma}

As a consequence of  Lemma \ref{l:equiv},  the spectral properties
of $N$ are equivalent to properties of harmonic functions of
polynomial growth on   $C(N)$.

\vskip2mm When   $N^{k-1}$  is   a submanifold  of
$\SS^{n-1}\subset \RR^n$, the   cone  over $N$ can be
isometrically embedded in $\RR^n$ as
\begin{equation}
    C(N) = \{ x \in \RR^n \mid x / |x|  \in N \} \, .
\end{equation}
     Note that
$C(N)$ is then invariant under dilations about the origin. We get
the following simple lemma whose proof is left for the reader:

\begin{lemma}
Suppose that $N^{k-1} \subset  \SS^{n-1}$.  The following are
equivalent: \begin{itemize} \item $N$ is minimal.  \item The
Euclidean mean curvature of $N \subset \RR^n$ is normal to
$\SS^{n-1} \subset \RR^n$. \item The coordinate functions are
eigenfunctions on $N$ with eigenvalue $k-1$. \item The cone $C(N)$
is minimal. \end{itemize}
\end{lemma}

\section{Manifolds with non--negative Ricci curvature}  \label{s:27}

In \cite{Ya3},  Yau extended the classical Liouville theorem to
complete manifolds $M$ with non--negative Ricci curvature:

  A positive (or bounded) harmonic function on $M$ must be
constant.

\noindent
In fact, by the gradient estimate of Cheng and Yau
(Theorem \ref{t:cy}), any harmonic function with sublinear growth
must be constant:

\begin{corollary}
\cite{CgYa} If $M$ is complete with $\Ric_M \geq 0$ and $d< 1$,
then $\cH_d (M) = \{ {\text{Constant functions}} \}$.
\end{corollary}

Since $\RR^n$ has non--negative Ricci curvature and the coordinate
functions are harmonic, this is obviously sharp. Therefore, when
$d\geq 1$ a different approach is needed.  Namely, instead of
showing a Liouville theorem, the point is to control the size of
the space of solutions. Over the years, there were many
interesting partial results (including two proofs when $M$ is a
surface with non--negative sectional curvature, \cite{LiTa2} and
\cite{DF}).  For example, in \cite{LiTa1}, Li and L.F. Tam
obtained the borderline case $d=1$, showing that
\begin{equation}    \label{e:linearg}
    \dim (\cH_1 (M)) \leq n+1 \, ,
\end{equation}
 for an $n$-dimensional manifold with $\Ric_M \geq 0$. This is
similar in   spirit to the classical comparison theorems since
$\dim (\cH_1 (\RR^n)) = n+1$ (the $n$ coordinate functions plus
the constant functions).  This corresponding rigidity theorem was
proven in \cite{ChCM} (see \cite{Li} for the special case where
$M$ is K\"ahler):

\begin{theorem}     \label{t:chcm}
\cite{ChCM}
 If $M$ is complete with $\Ric_M \geq 0$, then every
tangent cone at infinity $M_{\infty}$ splits isometrically as
\begin{equation}
    M_{\infty} = N \times \RR^{ \dim (\cH_1 (M)) - 1 } \, .
\end{equation}
  Hence, if $\dim (\cH_1 (M)) = n+1$, then \cite{C1} implies
that $M = \RR^n$.
\end{theorem}

Finally, in \cite{CM14},  the spaces of polynomial growth harmonic
functions were shown to be finite dimensional:

\begin{theorem}  \label{t:cm14}
\cite{CM14} If $M$ is complete with $\Ric_M \geq 0$, then $\cH_d
(M)$ is finite dimensional for each $d$.
\end{theorem}

 The proof of Theorem \ref{t:cm14}
consists of  two independent steps (the first does not use
harmonicity):
\begin{itemize} \item
Given a $2k$-dimensional subspace $H \subset \cH_d (M)$ and $h \in
(0,1]$, there exists a $k$-dimensional subspace $K \subset H$ and
$R> 0$ so that
\begin{equation}        \label{e:step1}
    \sup_{  v \in K \setminus \{ 0\} } \, \, \frac{\int_{ B_{(1+h)^2 R} }
    v^2}{\int_{ B_{R} }
    v^2} \leq C_1 \, (1+h)^{8d} \, .
\end{equation}
  \item
We  bound the dimension of a subspace $K$ of harmonic functions
satisfying \eqr{e:step1} in terms of $h$ and $d$.
\end{itemize}

To give some feel for the argument, we will sketch a proof of the
second step.

\begin{proof}(Sketch of second step)
 For simplicity, suppose that $R=1$ and $h=1$.  Fix a scale $r \in (0,1)$
 to be chosen small.
  We will use two properties of manifolds with $\Ric_M \geq 0$:

   First,
  we can find $N \leq
 C_n
 \, r^{-n}$  balls $B_r(x_i)$ with
 \begin{equation}       \label{e:vd}
    \chi_{B_1} \leq \sum_i \chi_{ B_{r} (x_i) } \leq C_n \,
    \chi_{B_2} \, ,
 \end{equation}
where $\chi_E$ is the characteristic function of a set $E$. (To do
this,  choose a maximal disjoint collection of balls of radius
$r/2$ and then use the volume comparison to get the second
inequality in \eqr{e:vd} and bound   $N$.)

Second, there is a uniform Neumann Poincar\'e inequality: If
$\int_{B_s(x)} f = 0$, then
\begin{equation}    \label{e:np}
    \int_{B_s(x)} f^2 \leq C_N \, s^2  \int_{B_s(x)} |\nabla f|^2
    \, .
\end{equation}

To bound the dimension of $K$, we will  construct a linear map
$\cM : K \to \RR^N$  and show that $\cM$ is injective for $r>0$
sufficiently small.  We define $\cM$ by
\begin{equation}
    \cM (v) = \left( \int_{B_{r}(x_1)} v \, ,  \cdots , \int_{B_{r}(x_N)}
    v  \right) \, .
\end{equation}
We will deduce a contradiction if $v \in K \setminus \{ 0 \}$ is
in the kernel of $\cM$. In particular, \eqr{e:np} gives that for
each $i$
\begin{equation}    \label{e:np2}
    \int_{B_r(x_i)} v^2 \leq C_N \, r^2  \int_{B_r(x_i)} |\nabla v|^2
    \, .
\end{equation}
Combining this with \eqr{e:vd} gives
\begin{equation}    \label{e:np3}
    \int_{B_1} v^2 \leq \sum_{i=1}^N \int_{B_r(x_i)} v^2 \leq C_N \, r^2 \sum_{i=1}^N \int_{B_r(x_i)} |\nabla v|^2
    \leq C_n \, C_N \, r^2 \,  \int_{B_2} |\nabla v|^2
    \, .
\end{equation}
We now (for the only time) use that $v$ is harmonic.  Namely, the
Caccioppoli inequality (or reverse Poincar\'e inequality) for
harmonic functions gives
\begin{equation}    \label{e:rp}
    \int_{B_2} |\nabla v|^2 \leq \int_{B_4} v^2 \, .
\end{equation}
Combining \eqr{e:np3} and \eqr{e:rp}, we get
\begin{equation}    \label{e:np4}
    \int_{B_1} v^2 \leq  C_n \, C_N \, r^2 \,  \int_{B_4} v^2
    \, .
\end{equation}
This contradicts \eqr{e:step1} if $r$ is sufficiently small,
completing the proof.
\end{proof}

On Euclidean space $\RR^n$, the spaces $\cH_d$ are given by
harmonic polynomials of degree at most $d$.  In particular, it is
not hard to see that
\begin{equation}    \label{e:ruppern}
    \dim (\cH_d (\RR^n)) \approx
    C \, d^{n-1} \, .
\end{equation}
 Using the correspondence between
harmonic polynomials and eigenfunctions on $\SS^{n-1}$ (see Lemma
\ref{l:equiv}), this
 is   closely related to Weyl's asymptotic formula on $\SS^{n-1}$.
 In \cite{CM15}, the authors proved a
similar sharp polynomial bound for manifolds with non--negative
Ricci curvature:

\begin{theorem}     \label{t:cm15}
\cite{CM15}  If $M^n$ is complete with $\Ric_M \geq 0$ and $d\geq
1$, then
\begin{equation}    \label{e:cm15}
    \dim ( \cH_d (M)) \leq C \, d^{n-1} \, .
\end{equation}
\end{theorem}

Taking $M = \RR^n$,   \eqr{e:ruppern}  illustrates that the
exponent $n-1$ is sharp in \eqr{e:cm15}.  However, as in Weyl's
asymptotic formula, the constant in front of $d^{n-1}$ can be
related to the volume. Namely,   we actually showed the stronger
statement
\begin{equation}    \label{e:cm15a}
    \dim ( \cH_d (M)) \leq C_n \, \V_M \, d^{n-1} + o (d^{n-1}) \,
    ,
\end{equation}
where
\begin{itemize} \item $C_n$ depends only on the
dimension $n$. \item
  $\V_M$ is the ``asymptotic volume ratio'' $\lim_{r\to \infty} \,
  \Vol (B_r)/ r^n$. \item $o(d^{n-1})$ is a function of $d$ with
  $\lim_{d\to \infty} \, o(d^{n-1})/d^{n-1} = 0$.
  \end{itemize}
As noted above, Theorem \ref{t:cm15} also gives lower bounds for
eigenvalues on a manifold $N^{n-1}$ with $\Ric_N \geq (n-2) =
\Ric_{\SS^{n-1}}$.  Using the sharper estimate \eqr{e:cm15a}
introduces the volume of $N$ into these eigenvalue estimates (as
predicted by Weyl's asymptotic formula).

An interesting  feature of these dimension estimates is that they
follow from ``rough'' properties of $M$ and are therefore
surprisingly stable under perturbation. For instance, in
\cite{CM14}, we actually proved Theorem \ref{t:cm14} for manifolds
with a volume doubling and a Neumann Poincar\'e inequality; unlike
a Ricci curvature bound, these properties are stable under
bi--Lipschitz transformations.

This finite dimensionality was not previously known even for
manifolds bi--Lipschitz to $\RR^n$ (except under additional
hypotheses, cf. results of Avellenada--Lin, \cite{AvLn}, and
Moser--Struwe, \cite{MrSt}).

\vskip2mm There are two particularly interesting directions which
have not been adequately explored.  The first is to develop
machinery to produce harmonic functions of polynomial growth.

\begin{question}
Suppose that $M^n$ has non--negative Ricci curvature and Euclidean
volume growth.  When can we produce harmonic functions of
polynomial growth on $M$?  The most interesting would be to solve
a  ``Dirichlet problem at infinity,'' where  polynomially growing
harmonic functions on a tangent cone at infinity give rise to
harmonic functions on $M$.  In complete generality, this is likely
to be rather delicate since these tangent cones need not be
unique.
\end{question}

The second direction is to get sharper dimension estimates for
holomorphic functions of polynomial growth.

\begin{question}
If $M^n$ is K\"ahler, then each holomorphic function is harmonic so
Theorem \ref{t:cm15} bounds the dimension  of the space of
polynomial growth holomorphic functions by $C \, d^{n-1}$.
However, on $\CC^{n/2}$, one get only $C \, d^{n/2}$ holomorphic
functions of degree $d$.   Does the sharper bound $C \, d^{n/2}$
hold?  A stronger curvature condition may be necessary (cf.
\cite{Ni} for one such result).
\end{question}

This is just a very brief overview (omitting many interesting
results), but we hope that it gives something of the flavor of the
subject; the interested reader may consult
 \cite{CM13} and references
therein for more information.

\section{Minimal surfaces and a generalized Bernstein theorem}  \label{s:28}

We next describe a similar finite dimensionality result for
minimal submanifolds and an application of this -- a ``generalized
Bernstein theorem'' -- proven in \cite{CM16}.  Recall that the
Bernstein theorem implies that, through dimension seven,
area--minimizing hypersurfaces are affine. A weaker form of this
is true in all dimensions by the Allard regularity theorem
\cite{Al}:
\begin{quotation}
There exists $\delta = \delta (k , n) > 0$ such that if
$\Sigma^k \subset \RR^{n}$ is a complete immersed minimal
submanifold with
\begin{equation}
        \frac{\Vol (B_r   \cap \Sigma)}{ \Vol (B_r \subset \RR^k) } \leq (1+
        \delta)
\end{equation}
for all $r$, then $\Sigma$ is  affine.
\end{quotation}
The
generalized Bernstein theorem, which should perhaps be called a
generalized Allard theorem instead, shows that {\it any} upper
 bound on the density gives a
corresponding upper bound for the dimension of the smallest affine
subspace containing the minimal surface.

The results of this section apply to a large class of generalized
minimal submanifolds   $\Sigma^k \subset \RR^{n}$: stationary
rectifiable $k$--varifolds with density   at least $1$ a.e. on the
support.  This
 includes the case of embedded minimal submanifolds and, for
 simplicity, we will focus on this case below.

\begin{theorem}  \cite{CM16}    \label{c:minimal}
If $\Sigma^k \subset \RR^{n}$ has density  bounded by $\V_{\Sigma}
$, then $\Sigma$ must be contained in an affine subspace of
dimension at most $C_k \, \V_{\Sigma} $.
\end{theorem}

Another way to think of Theorem \ref{c:minimal} is that it bounds
the number of linearly independent coordinate functions on
$\Sigma$ in terms of its volume.  The linear dependence  on
$\V_{\Sigma}$ in Theorem \ref{c:minimal} is sharp; namely, any
bound of the form $C_k \V_{\Sigma}^{\alpha}$ must have $\alpha
\geq 1$.

\medskip
  For a submanifold $\Sigma \subset \RR^n$, we will define the
  spaces of harmonic functions of polynomial
growth using the extrinsic distance; it will be clear from the
context which definition we are using. Namely,
  given $\Sigma \subset \RR^n$ and $d>0$,
 we define the  vector spaces $\cH_d (\Sigma)$ of harmonic functions of
polynomial
growth  by:

 \begin{definition}   \rm
A function $u$ is in  $\cH_d (\Sigma)$ if $u$ is harmonic on
$\Sigma$ and
\begin{equation}
    |u(x)| \leq C (1 +|x|^d) \, ,
\end{equation}
 for some $C$. Thus, the coordinate functions $x_i$ are in
$\cH_1 (\Sigma)$.
\end{definition}

Since the coordinate functions are   harmonic on $\Sigma$ (cf.
Proposition \ref{p:haco}), Theorem \ref{c:minimal} follows from a
bound for the dimensions of the spaces of harmonic functions on
$\Sigma$  of polynomial growth:

\begin{theorem}  \cite{CM16}     \label{t:minimal}
Let   $\Sigma^k \subset \RR^{n}$ have density  bounded by
$\V_{\Sigma} $.  For any $d \geq 1$,
\begin{equation}
        \dim \cH_d (\Sigma) \leq C_k \, \V_{\Sigma} \, d^{k-1}  \,
        .
\end{equation}
\end{theorem}

(The spectral properties of spherical minimal submanifolds have
been studied in their own right; see, for instance,
 Cheng--Li--Yau \cite{CgLiYa} or Choi--Wang \cite{CiWa}.)

\subsection{Other applications of function theory}
 We have just seen an application of function theory to
describe the geometry of the underlying space (in this case, a
bound on the dimension of the space of linear growth harmonic
functions controlled the complexity of the minimal submanifold).
Another example is the Bochner technique and resulting topological
restrictions of curvature.  There are  many other examples and,
indeed, often these sorts of applications motivate developments in
function theory.
  This is perhaps most evident in  (one variable) complex analysis, where
  function theory has played a major role.

  Function theory has also played an important role in  the theory of  quasi--regular
  maps, see \cite{G}. Recall that a map $F: \RR^n \to \RR^n$ is
  $K$--quasi--regular if $F$ and $dF$ are in $L^n$
  and
at almost every point we have
  \begin{equation}
    |dF| \leq K \, \det (dF)  \, .
  \end{equation}
   For instance, M. Bonk and J. Heinonen used function theoretic
  arguments to prove:

\begin{theorem}   \cite{BoHj}  \label{t:bonkhein}
If $M$ is a compact $n$--dimensional manifold and $F:\RR^n \to M$
is a (non--trivial) $K$--quasi--regular map, then the dimension of
the de Rham cohomology ring of $M$ is at most $C=C(n,K)$.
\end{theorem}

Finally, we note that in the theory of quasi--regular maps, the
most natural functions to study are no longer the harmonic ones.
Rather, one is interested in $\mathcal{A}$--harmonic functions,
i.e., functions $u$ satisfying
\begin{equation}
    \dv ( \mathcal{A} (\nabla u)) = 0 \, ,
\end{equation}
where $\mathcal{A}$ is a {\it nonlinear} map on the tangent space
satisfying several natural conditions (e.g., taking
$\mathcal{A}(x) = |x|^{p-2} \, x$ gives the so--called
$p$--Laplacian).

\section{Volumes for eigensections}     \label{s:29}

Let $M^n$ be a closed $n$-dimensional Riemannian manifold
and $V$ a vector bundle over $M$.
Suppose that $\{f_i\}$ is an $L^2$-orthonormal set of eigensections
of $V$ of the Laplacian with eigenvalues $\lambda_i$
(where $0=\lambda_0\leq \lambda_1\leq \lambda_2\leq \cdots$), that is
\begin{equation}
\int_M f_i\,f_j=\delta_{i,j}\text{ and }\Delta\,f_i+\lambda_i\,f_i=0\, .
\end{equation}
Given
$a = (a_1 , \cdots , a_i) \in \SS^{i-1}$,
we define a function $F^a_i (\epsilon)$ by
\begin{equation}        \label{e:fai}
        F^a_i (\epsilon) = \Vol (
        \{ x \in M \, | \,
        | \sum_{j=1}^i a_j \, f_j (x) | < \epsilon \} ) \, .
\end{equation}
Set $F_i(\epsilon)=F_i^{e_i}(\epsilon)$ where $e_i=(\delta_{i,j})_j$
and define $F^A_i(\epsilon)$ by averaging
$F_i^a(\epsilon)$ over $a \in \SS^{i-1}$
\begin{equation}        \label{e:fAi}
        F^A_i (\epsilon) = \mint_{ a \in \SS^{i-1}}  F^a_i (\epsilon)
         \, .
\end{equation}

In this section we will discuss the answer to
the following question of S.T. Yau:

\begin{quotation}
Let $M^n$ be closed and $V=\Omega^1(M)$ the bundle of
one forms on $M$.  Then $\limsup_{i \to \infty} F_i (\epsilon)$ and
$\liminf_{i \to \infty} F_i (\epsilon)$ are interesting
functions of $\epsilon$.  Are they positive?
Can one estimate the behavior of $\epsilon^{-n}
\liminf_{i \to \infty} F_i (\epsilon)$ as $\epsilon \to 0$?
The problem may be easier if we replace
$F_i(\epsilon)$ by $F^A_i(\epsilon)$.  We can of course consider problems
for $p$-forms with $p>1$.
\end{quotation}

It turns out that on a flat square torus
$\liminf_{i\to\infty}F_i=0$ on $\Omega^1$.
However as the next theorem shows then for eigenfunctions on any manifold
$F_i$ is positive.  We will also see below in Theorem \ref{c:averageN}
that the average $F^A_i$ for $p$-forms
is positive and we give a sharp lower bound.

\begin{theorem}   \label{t:lower}
\cite{CM23}
Let $M^n$ be closed
with
\begin{equation}  \label{e:l1}
\Ric_{M}\geq -(n-1)\, .
\end{equation}
There exists $C=C(n)>0$ and $\Lambda=\Lambda (n)>0$ such that
if $f$ is a
eigenfunction of $\Delta$ on $M$ with eigenvalue $\lambda\geq \Lambda$
and $0<\epsilon\leq 1$, then
\begin{equation}  \label{e:l22}
\Vol \, \left( \left\{x\in M|\,|f|^2
<\epsilon^2 \mint_M|f|^2\right\} \right)
\geq C\,\epsilon^{n}\,\Vol\, (M)\, .
\end{equation}
\end{theorem}

It is possible to generalize Theorem \ref{t:lower} to the case
where $M$ is assumed to have the doubling property and satisfy the
Neumann Poincar\'e inequality for $r\leq 1$. In this case,
however, the exponent in $\epsilon$ may not be $n$ but rather will
depend on the doubling constant $C_D$ and the constant $C_p$ in
the Poincar\'e inequality.

In contrast to $F_i$ for eigenforms, then the next theorem shows that the
average $F_i^A$ has always a positive lower bound.

\begin{theorem}   \label{c:averageN}
\cite{CM23}
There exists $C=C(q)$ so that if $M^n$ is closed and $V^{q}$ is a rank
$q$
vector bundle over $M$
with Laplace-type operator $\Delta_V$, then
for $0<\epsilon\leq 1$ and $i> q$
\begin{equation}  \label{e:Nc}
        F_i^A(\epsilon)
        \geq C\,\epsilon^q  \, [\Vol \, (M)]^{1 + q/2} \, .
\end{equation}
\end{theorem}

In contrast to
eigenfunctions, $F_i (\epsilon)$ need not
be positive for eigenforms:

Let $T^2$ be a flat square torus with side lengths $2\,\pi$ and
define one forms by
\begin{equation}
    \sigma_m =\cos (m\,x_1)\,dx_1+ \sin (m\,x_1)\,dx_2 \, ,
\end{equation}
 then for all $m$
\begin{equation}
        \{x  \in T^2\,|\, |\sigma_m |<1\}  =  \emptyset\text{ and
        hence }
\liminf_{i\to \infty} F_i(\epsilon)=0 \text{ for $\epsilon<1/ (4\pi^2)$}\, .
\end{equation}

\bibliographystyle{plain}

\end{document}